\documentclass[11pt]{article}
\usepackage{amsfonts,amssymb,amsmath,latexsym}

\newcommand{\file}{$\ti{\ }$/wisk/klaus/ren.tex\quad}
\renewcommand{\file}{}

\newcommand{\detail}[1]{\par\noi{\bf [Proof detail\ }{#1}
\hfill{\bf ]}\par\noi\hspace{-4pt}}
\renewcommand{\detail}[1]{}

\usepackage[dvips]{graphicx}

\newcommand{\si}{\ensuremath{\sigma}}

\newcommand{\di}{\mathrm{d}}
\newcommand{\dis}{\displaystyle}
\newcommand{\txt}{\textstyle}

\newcommand{\vc}{\vspace{\baselineskip}}

\newcommand{\med}{\medskip}
\newcommand{\noi}{\noindent}

\newcommand{\half}{{[0,\infty)}}
\newcommand{\R}{{\mathbb R}}
\newcommand{\N}{{\mathbb N}}
\newcommand{\Z}{{\mathbb Z}}

\newcommand{\Q}{{\mathbb Q}}

\newcommand{\Ai}{{\cal A}}

\newcommand{\Ci}{{\cal C}}
\newcommand{\Di}{{\cal D}}

\newcommand{\Fi}{{\cal F}}
\newcommand{\Gi}{{\cal G}}
\newcommand{\Hi}{{\cal H}}

\newcommand{\Li}{{\cal L}}
\newcommand{\Mi}{{\cal M}}
\newcommand{\Ni}{{\cal N}}

\newcommand{\Pc}{{\cal P}}
\newcommand{\Qi}{{\cal Q}}
\newcommand{\Ri}{{\cal R}}
\newcommand{\Si}{{\cal S}}

\newcommand{\Ui}{{\cal U}}
\newcommand{\Vi}{{\cal V}}
\newcommand{\Wi}{{\cal W}}
\newcommand{\Xc}{{\cal X}}
\newcommand{\Yi}{{\cal Y}}
\newcommand{\Zi}{{\cal Z}}

\newcommand{\al}{\alpha}

\newcommand{\bet}{\beta}
\newcommand{\ga}{\gamma}
\newcommand{\Ga}{\Gamma}
\newcommand{\de}{\delta}

\newcommand{\eps}{\varepsilon}
\newcommand{\la}{\lambda}

\newcommand{\sig}{\sigma}

\newcommand{\tet}{\theta}
\newcommand{\oo}{\omega}
\newcommand{\om}{\Omega}

\newcommand{\expo}{\mbox{\large\it e}}
\newcommand{\ex}[1]{\expo^{\,\textstyle{#1}}}

\newcommand{\as}{\mbox{a.s.}}

\newcommand{\volgt}{\ensuremath{\Rightarrow}}

\newcommand{\up}{\uparrow}
\newcommand{\down}{\downarrow}

\newcommand{\li}{\langle}

\newcommand{\re}{\rangle}

\newcommand{\ti}{\tilde}

\newcommand{\ov}{\overline}
\newcommand{\un}{\underline}

\newcommand{\pa}{\partial}

\newcommand{\tr}{{\rm tr}}

\newcommand{\sub}{\subset}
\newcommand{\beh}{\backslash}

\newcommand{\ffrac}[2]{{\textstyle\frac{{#1}}{{#2}}}}
\newcommand{\dif}[1]{\ffrac{\partial}{\partial{#1}}}
\newcommand{\diff}[1]{\ffrac{\partial^2}{{\partial{#1}}^2}}
\newcommand{\difif}[2]{\ffrac{\partial^2}{\partial{#1}\partial{#2}}}

\newcommand{\Diff}[1]{\frac{\partial^2}{{\partial{#1}}^2}}

\newcommand{\be}{\begin{equation}}
\newcommand{\ee}{\end{equation}}
\newcommand{\ba}{\begin{array}}
\newcommand{\ea}{\end{array}}
\newcommand{\bc}{\be\begin{array}{r@{\,}c@{\,}l}}
\newcommand{\ec}{\end{array}\ee}

\newtheorem{theorem}{Theorem}[section]
\newtheorem{proposition}[theorem]{Proposition}
\newtheorem{corollary}[theorem]{Corollary}
\newtheorem{conjecture}[theorem]{Conjecture}
\newtheorem{lemma}[theorem]{Lemma}

\newtheorem{remark}[theorem]{Remark}

\newtheorem{defi}[theorem]{Definition}
\newtheorem{note}[theorem]{Remark}
\newtheorem{question}[theorem]{Question}

\newcommand{\bt}{\begin{theorem}}
\newcommand{\et}{\end{theorem}}
\newcommand{\bl}{\begin{lemma}}
\newcommand{\el}{\end{lemma}}
\newcommand{\bp}{\begin{proposition}}
\newcommand{\ep}{\end{proposition}}
\newcommand{\bcor}{\begin{corollary}}
\newcommand{\ecor}{\end{corollary}}
\newcommand{\br}{\begin{remark}\rm}
\newcommand{\er}{\end{remark}}
\newcommand{\bdf}{\begin{defi}\rm}
\newcommand{\edf}{\hfill$\Diamond$\end{defi}}
\newcommand{\brm}{\begin{note}\rm}
\newcommand{\erm}{\hfill$\Diamond$\end{note}}
\newcommand{\bqu}{\begin{question}\rm}
\newcommand{\equ}{\hfill$\Diamond$\end{question}}

\newcommand{\halmos}{\rule{1ex}{1.4ex}}
\def \qed {\nopagebreak{\hspace*{\fill}$\halmos$\medskip}}

\begin{document}
\newcommand{\Widg}{\Wi_{\rm\scriptscriptstyle DG}}
\newcommand{\nul}{\makebox[7pt]{$\scriptstyle 0$}}
\newcommand{\mm}{\hspace{-1pt}-\hspace{-1pt}}
\newcommand{\bra}{\bet}
\newcommand{\isd}{\stackrel{\scriptscriptstyle\Di}{=}}
\newcommand{\vagto}{\stackrel{\mbox{\scriptsize vague}}{\longrightarrow}}
\newcommand{\Pois}{{\rm Pois}}
\newcommand{\Thin}{{\rm Thin}}
\newcommand{\Pss}{{\cal P}\!{\it ois}}
\newcommand{\bplim}{\mbox{\rm bp-}\lim}
\newcommand{\exlim}{\mbox{\rm ex-}\!\!\lim}
\newcommand{\subbp}[2]{_{\ba{c}\scriptstyle{#1}\\[-.05cm]\scriptstyle{#2}\ea}}
\newcommand{\x}{\mathbf x}
\newcommand{\y}{\mathbf y}
\newcommand{\vb}{\mathbf v}

\newcommand{\asto}[1]{\underset{{#1}\to\infty}{\longrightarrow}}
\newcommand{\Asto}[1]{\underset{{#1}\to\infty}{\Longrightarrow}}
\newcommand{\asim}[1]{\underset{{#1}\to\infty}{\sim}}
\newcommand{\astoo}[1]{\underset{{#1}\to 0}{\longrightarrow}}

\newcommand{\hut}{\ov}

\renewcommand{\labelenumi}{(\roman{enumi})}

\makeatletter\@addtoreset{equation}{section}\makeatother\def\theequation{\thesection.\arabic{equation}} 

\title{\vspace{-2.5cm} Renormalization analysis of catalytic Wright-Fisher diffusions}
\author{Klaus Fleischmann\vspace{6pt}\\
{\small Weierstrass Institute for Applied}\\
{\small Analysis and Stochastics}\\
{\small Mohrenstr.\ 39}\\
{\small D--10117 Berlin, Germany\vspace{3pt}}\\
{\small e-mail: fleischm@wias-berlin.de}\vspace{8pt}
\and Jan M. Swart\vspace{6pt}\\
{\small Mathematical Institute}\\
{\small University T\"ubingen}\\
{\small Auf der Morgenstelle 10}\\
{\small D--72076 T\"ubingen, Germany\vspace{3pt}}\\
{\small e-mail: jan.swart@uni-tuebingen.de}\vspace{4pt}}
\date{{\file\today}}
\maketitle
\vspace{-1cm}
\begin{abstract}\noindent
Recently, several authors have studied maps where a function,
describing the local diffusion matrix of a diffusion process with a
linear drift towards an attraction point, is mapped into the average
of that function with respect to the unique invariant measure of the
diffusion process, as a function of the attraction point. Such
mappings arise in the analysis of infinite systems of diffusions
indexed by the hierarchical group, with a linear attractive
interaction between the components. In this context, the mappings
are called renormalization transformations. We consider such maps
for catalytic Wright-Fisher diffusions. These are diffusions on the
unit square where the first component (the catalyst) performs an
autonomous Wright-Fisher diffusion, while the second component (the
reactant) performs a Wright-Fisher diffusion with a rate depending
on the first component through a catalyzing function. We determine
the limit of rescaled iterates of renormalization transformations
acting on the diffusion matrices of such catalytic Wright-Fisher
diffusions.
\end{abstract}

{\setlength{\parskip}{-2pt}
\footnotesize\tableofcontents}

\vspace{15pt}

\noi
{\bf\Large Part I}
\vspace{-5pt}

\section{Introduction and main result}\label{intro}

Several authors \cite{BCGH95,BCGH97,HS98,Sch98,CDG04} have studied
maps where a function, describing the local diffusion matrix of a
diffusion process, is mapped into the average of that function with
respect to the unique invariant measure of the diffusion process
itself. Such mappings arise in the analysis of infinite systems of
diffusion processes indexed by the hierarchical group, with a linear
attractive interaction between the components \cite{DG93a,DG96,DGV95}.
In this context, the mappings are called renormalization
transformations. We follow this terminology. For more on the relation
between hierarchically interacting diffusions and renormalization
transformations, see Appendix~\ref{hier}.

Formally, such renormalization transformations can be defined as follows.
\bdf\label{defi}
{\bf (Renormalization class and transformation)} Let $D\sub\R^d$ be
nonempty, convex, and open. Let $\Wi$ be a collection of continuous
functions $w$ from the closure $\ov D$ into the space $M^d_+$ of
symmetric non-negative definite $d\times d$ real matrices, such that
$\la w\in\Wi$ for every $\la>0$, $w\in\Wi$. We call $\Wi$ a 
{\em prerenormalization class} on $\ov D$ if the following three
conditions are satisfied:
\begin{enumerate}
\item For each constant $c>0$, $w\in\Wi$, and $x\in\ov D$, the
martingale problem for the operator $A^{c,w}_x$ is well-posed, where
\be\label{Adef}
A^{c,w}_x f(y):=\sum_{i=1}^d c\,(x_i-y_i)\dif{y_i}f(y)
+\sum_{i,j=1}^dw_{ij}(y)\difif{y_i}{y_j}f(y)\qquad(y\in\ov D),
\ee
and the domain of $A^{c,w}_x$ is the space of real functions on $\ov
D$ that can be extended to a twice continuously differentiable
function on $\R^d$ with compact support.
\item For each $c>0$, $w\in\Wi$, and $x\in\ov D$, the martingale
problem for $A^{c,w}_x$ has a unique stationary solution with
invariant law denoted by $\nu^{c,w}_x$.
\item For each $c>0$, $w\in\Wi$, $x\in\ov D$, and $i,j=1,\ldots,d$,
one has $\dis\int_{\ov D}\nu^{c,w}_x(\di y)|w_{ij}(y)|<\infty$.
\end{enumerate}
If $\Wi$ is a prerenormalization class, then we define for each $c>0$
and $w\in\Wi$ a matrix-valued function $F_cw$ on $\ov D$ by
\be\label{Fc}
F_cw(x):=\int_{\ov D}\nu^{c,w}_x(dy)w(y)\qquad(x\in\ov D).
\ee
We say that $\Wi$ is a {\em renormalization class} on $\ov D$ if in addition:
\begin{enumerate}\addtocounter{enumi}{3}
\item For each $c>0$ and $w\in\Wi$, the function $F_cw$ is an element of $\Wi$.
\end{enumerate}

If $\Wi$ is a renormalization class and $c>0$, then the map
$F_c:\Wi\to\Wi$ defined by (\ref{Fc}) is called the 
{\em renormalization transformation} on $\Wi$ with {\em migration
constant} $c$. In (\ref{Adef}), $w$ is called the {\em diffusion
matrix} and $x$ the {\em attraction point.}
\edf
\brm{\bf (Associated SDE)}\label{sderem}\hspace{-.4pt}
It is well-known that $\ov D$-valued (weak) solutions
$\y=(\y^1,\ldots,\y^d)$ to the stochastic differential equation (SDE)
\be\label{sde}
\di\y^i_t=c\,(x_i-\y^i_t)\di t
+\sqrt 2\sum_{j=1}^{n}\sig_{ij}(\y_t)\di B^j_t\qquad(t\geq 0,\ i=1,\ldots,d),
\ee
where $B=(B^1,\ldots,B^n)$ is $n$-dimensional (standard) Brownian
motion ($n\geq 1$), solve the martingale problem for $A^{c,w}_x$ if
the $d\times n$ matrix-valued function $\sig$ is continuous and
satisfies $\sum_k\sig_{ik}\sig_{jk}=w_{ij}$. Conversely
\cite[Theorem~5.3.3]{EK}, every solution to the martingale problem for
$A^{c,w}_x$ can be represented as a solution to the SDE (\ref{sde}),
where there is some freedom in the choice of the root $\sig$ of the
diffusion matrix $w$.
\erm

\noindent
In the present paper, we concern ourselves with the following
renormalization class on $[0,1]^2$.
\bdf\label{Wcatdef}
{\bf (Renormalization class of catalytic Wright-Fisher diffusions)}
We set $\Wi_{\rm cat}:=\{w^{\al,p}:\al>0,\ p\in\Hi\}$, where
\be\label{walp}
w^{\al,p}(x):=\left(\ba{@{}cc@{}}\al x_1(1-x_1)&0\\
0&p(x_1)x_2(1-x_2)\ea\right)\qquad(x=(x_1,x_2)\in[0,1]^2),
\ee
and
\be\label{Hdef}
\Hi:=\{p:p\mbox{ a real function on }[0,1],\ p\geq 0,
\ p\mbox{ Lipschitz continuous}\}.
\ee
Moreover, we put
\be\label{Hlr}
\Hi_{l,r}:=\{p\in\Hi:\; 1_{\{p(0)>0\}}=l,\ 1_{\{p(1)>0\}}=r\}\qquad(l,r=0,1),
\ee
and set $\Wi^{l,r}_{\rm cat}:=\{w^{\al,p}:\al>0,
\ p\in\Hi_{l,r}\}\quad(l,r=0,1)$.
\edf
By Remark~\ref{sderem}, solutions $\y=(\y^1,\y^2)$ to the martingale
problem for $A^{c,w^{\al,p}}_x$ can be represented as solutions to the SDE
\be\ba{rr@{\,}c@{\,}l}\label{catsde}
{\rm (i)}&\dis\di\y^1_t&=&\dis c\,(x_1-\y^1_t)\di t
+\sqrt{2\al\y^1_t(1-\y^1_t)}\di B^1_t,\\[5pt]
{\rm (ii)}&\dis\di\y^2_t&=&\dis c\,(x_2-\y^2_t)\di t
+\sqrt{2p(\y^1_t)\y^2_t(1-\y^2_t)}\di B^2_t.
\ec
We call $\y^1$ the Wright-Fisher {\em catalyst} with {\em resampling
rate} $\al$ and $\y^2$ the Wright-Fisher {\em reactant} with {\em
catalyzing function} $p$.\med

\noi
For any renormalization class $\Wi$ and any sequence of (strictly)
positive migration constants $(c_k)_{k\geq 0}$, we define {\em
iterated renormalization transformations} $F^{(n)}:\Wi\to\Wi$, as
follows:
\be\label{Fn}
F^{(n+1)}w:=F_{c_n}(F^{(n)}w)\quad(n\geq 0)\quad\mbox{with}
\quad F^{(0)}w:=w\qquad(w\in\Wi_{\rm cat}).
\ee
We set $s_0:=0$ and
\be\label{sn}
s_n:=\sum_{k=0}^{n-1}\frac{1}{c_k}\qquad(1\leq n\leq\infty).
\ee
Here is our main result:
\bt{\bf (Main result)}\label{main}\newline
{\bf (a)} The set $\Wi_{\rm cat}$ is a renormalization class on
$[0,1]^2$ and $F_c(\Wi^{l,r}_{\rm cat})\sub\Wi^{l,r}_{\rm cat}$
$(c>0,\ l,r=0,1)$.\vc

\noi
{\bf (b)} Fix (positive) migration constants $(c_k)_{k\geq 0}$ such that
\be\label{sga}
{\rm (i)}\quad s_n\asto{n}\infty\qquad\mbox{and}\qquad{\rm (ii)}
\quad\frac{s_{n+1}}{s_n}\asto{n}1+\ga^\ast
\ee
for some $\ga^\ast\geq 0$. If $w\in\Wi^{l,r}_{\rm cat}$ $(l,r=0,1)$,
then uniformly on $[0,1]^2$,
\be\label{wast}
s_nF^{(n)}w\asto{n}w^\ast,
\ee
where the limit $w^\ast$ is the unique solution in $\Wi^{l,r}_{\rm cat}$
to the equation
\be\ba{rr@{\,}c@{\,}ll}\label{wiga}
{\rm (i)}&\dis(1+\ga^\ast)F_{1/\ga^\ast}w^\ast&=&\dis w^\ast\qquad
&\mbox{if}\ \ \ga^\ast>0,\\[5pt]
{\rm (ii)}&
\dis\ffrac{1}{2}\sum_{i,j=1}^2w^\ast_{ij}(x)\difif{x_i}{x_j}w^\ast(x)+w^\ast(x)
&=&\dis 0\qquad(x\in[0,1]^2)\quad&\mbox{if}\ \ \ga^\ast=0.
\ec
{\bf (c)} The matrix $w^\ast$ is of the form $w^\ast=w^{1,p^\ast}$,
where $p^\ast=p^\ast_{l,r,\ga^\ast}\in\Hi_{l,r}$
depends on $l,r,$ and~$\ga^\ast$. One has
\be
p^\ast_{0,0,\ga^\ast}\equiv 0\quad\mbox{and}
\quad p^\ast_{1,1,\ga^\ast}\equiv 1\qquad\mbox{ for all }\ga^\ast\geq 0.
\ee
For each $\ga^\ast\geq 0$, the function $p^{\ast}_{0,1,\ga^\ast}$ is
concave, nondecreasing, and satisfies $p^\ast_{0,1,\ga^\ast}(0)=0$,
$p^\ast_{0,1,\ga^\ast}(1)=1$. By symmetry, analoguous statements hold
for $p^\ast_{1,0,\ga^\ast}$.
\et
Conditions (\ref{sga})~(i) and (ii) are satisfied, for example, for
$c_k=(1+\ga^\ast)^{-k}$. Note that the functions
$p^\ast_{0,0,\ga^\ast}$ and $p^\ast_{1,1,\ga^\ast}$ are independent of
$\ga^\ast\geq 0$. We believe that on the other hand,
$p^\ast_{0,1,\ga^\ast}$ is not constant as a function of $\ga^\ast$,
but we have not proved this.

The function $p^\ast_{0,1,0}$ is the unique nonnegative
solution to the equation
\be
\ffrac{1}{2}x(1-x)\diff{x}p(x)+p(x)(1-p(x))=0\qquad(x\in[0,1])
\ee
with boundary conditions $p(0)=0$ and $p(1)>0$. This function occurred
before in the work of Greven, Klenke, and Wakolbinger
\cite[formulas~(1.10)--(1.11)]{GKW01}. In Section~\ref{discus} we
discuss the relation between their work and ours.\med

\noi
{\bf Outline} In Part~I of the paper
(Sections~\ref{intro}--\ref{dispro}) we present our results and our
main techniques for proving them. Part~II
(Sections~\ref{Wcatsec}--\ref{final}) contains detailed proofs. Since
the motivation for studying renormalization classes comes from the
study of linearly interacting diffusions on the hierarchical group, we
explain this connection in Appendix~\ref{app}.\med

\noi
{\bf Outline of Part~I} In the next section, we place our main result
in a broader context. We give a more thorough introduction to the
theory of renormalization classes on compact sets and discuss earlier
results on this topic. In Section~\ref{bracon}, we discuss special
properties of the renormalization class $\Wi_{\rm cat}$ from
Definition~\ref{Wcatdef}. In particular, we show how techniques from
the theory of spatial branching processes can be used to prove
Theorem~\ref{main}. In Section~\ref{dispro} we discuss the relation of
our work with that in \cite{GKW01} and mention some open problems.\med

\noi
{\bf Notation} If $E$ is a separable, locally compact, metrizable
space, then $\Ci(E)$ denotes the space of continuous real functions on
$E$. If $E$ is compact then we equip $\Ci(E)$ with the supremumnorm
$\|\cdot\|_\infty$. We let $B(E)$ denote the space of all bounded
Borel measurable real functions on $E$. We write $\Ci_+(E)$ and
$\Ci_{[0,1]}(E)$ for the spaces of all $f\in\Ci(E)$ with $f\geq 0$ and
$0\leq f\leq 1$, respectively, and define $B_+(E)$ and $B_{[0,1]}(E)$
analogously. We let $\Mi(E)$ denote the space of all finite measures
on $E$, equipped with the topology of weak convergence. The subspaces
of probability measures is denoted by $\Mi_1(E)$. We write $\Ni(E)$
for the space of finite counting measures, i.e., measures of the form
$\nu=\sum_{i=1}^m\de_{x_i}$ with $x_1,\ldots,x_m\in E$ ($m\geq 0$). We
interpret $\nu$ as a collection of particles, situated at positions
$x_1,\ldots,x_m$. For $\mu\in\Mi(E)$ and $f\in B(E)$ we use the
notation $\li\mu,f\re:=\int_E f\,\di\mu$ and $|\mu|:=\mu(E)$.  By
definition, $\Di_E\half$ is the space of cadlag functions $w:\half\to
E$, equipped with the Skorohod topology. We denote the law of a random
variable $y$ by $\Li(y)$. If $\y=(\y_t)_{t\geq 0}$ is a Markov process
in $E$ and $x\in E$, then $P^x$ denotes the law of $\y$ started in
$\y_0=x$. If $\mu$ is a probability law on $E$ then $P^\mu$ denotes
the law of $\y$ started with initial law $\Li(\y_0)=\mu$. For
time-inhomogeneous processes, we use the notation $P^{t,x}$ or
$P^{t,\mu}$ to denote the law of the process started at time $t$ with
initial state $\y_t=x$ or initial law $\Li(\y_t)=\mu$, respectively.
We let $E^x,E^\mu,\ldots$ etc.\ denote expectation with respect to
$P^x,P^\mu,\ldots$, respectively.

\section{Renormalization classes on compact sets}\label{gensec}

\subsection{Some general facts and heuristics}

In this section, we explain that our main result is a special case of
a type of theorem that we believe holds for many more renormalization
classes on compact sets in $\R^d$. Moreover, we describe some
elementary properties that hold generally for such renormalization
classes. The proofs of Lemmas~\ref{contlem}--\ref{L:shape} can be
found in Section~\ref{gensub} below.

Fix a prerenormalization class $\Wi$ on a set $\ov D$ where
$D\sub\R^d$ is open, bounded, and convex. Then $\Wi$ is a subset of
the cone $\Ci(\ov D,M^d_+)$ of continuous $M^d_+$-valued functions on
$\ov D$. We equip $\Ci(\ov D,M^d_+)$ with the topology of uniform
convergence. Our first lemma says that the equilibrium measures
$\nu^{c,w}_x$ and the renormalized diffusion matrices $F_cw(x)$ are
continuous in their parameters.
\bl{\bf(Continuity in parameters)}\label{contlem}
\begin{itemize}
\item[{\bf (a)}] The map $(x,c,w)\mapsto\nu^{c,w}_x$ from
$\ov D\times(0,\infty)\times\Wi$ into $\Mi_1(\ov D)$ is continuous.
\item[{\bf (b)}] The map $(x,c,w)\mapsto F_cw(x)$ from
$\ov D\times(0,\infty)\times\Wi$ into $M^d_+$ is continuous.
\end{itemize}
\el
In particular, $x\mapsto\nu^{c,w}_x$ is a continuous probability
kernel on $\ov D$, and $F_cw\in\Ci(\ov D,M^d_+)$ for all $c>0$ and
$w\in\Wi$. Recall from Definition~\ref{defi} that $\la w\in\Wi$ for
all $w\in\Wi$ and $\la>0$. The reason why we have included this
assumption is that it is convenient to have the next scaling lemma
around, which is a consequence of time scaling.
\bl{\bf(Scaling property of renormalization transformations)}\label{schaal}
One has
\be\left.\ba{rr@{\,}c@{\,}l}\label{schafo}
{\rm (i)}&\nu^{\la c,\la w}_x&=&\nu^{c,w}_x\\
{\rm (ii)}&F_{\la c}(\la w)&=&\la F_cw\\
\ea\qquad\right\}\quad(\la,c>0,\ w\in\Wi,\ x\in\ov D).
\ee
\el
The following simple lemma will play a crucial role in what follows.
\bl{\bf (Mean and covariance matrix)}\label{numom}
For all $x\in\ov D$ and $i,j=1,\ldots,d$, the mean and covariances of
$\nu^{c,w}_x$ are given by
\be\ba{rr@{\;}c@{\;}l}\label{moments2}
{\rm (i)}&\dis\int_{\ov D}\nu^{c,w}_x(\di y)(y_i-x_i)&=&0,\\
{\rm (ii)}&\dis\int_{\ov D}\nu^{c,w}_x(\di y)(y_i-x_i)(y_j-x_j)
&=&\frac{1}{c}F_cw_{ij}(x).
\ec 
\el
For any $w\in\Ci(\ov D,M^d_+)$, we call
\be\label{eff}
\pa_w D:=\{x\in\ov D:w_{ij}(x)=0\ \forall i,j=1,\ldots,d\}
\ee
the {\em effective boundary} of $D$ (associated with $w$). If $\y$ is
a solution to the martingale problem for the operator
$\sum_{i,j=1}^dw_{ij}(y)\difif{y_i}{y_j}$ (i.e., the operator in
(\ref{Adef}) without the drift), then, by martingale convergence,
$\y_t$ converges a.s.\ to a limit $\y_\infty$; it is not hard to see
that $\y_\infty\in\pa_w D$ a.s. The next lemma says that the effective
boundary is invariant under renormalization.
\bl[Invariance of effective boundary]\label{effinv}
One has $\pa_{F_cw}D=\pa_wD$ for all $w\in\Wi$, $c>0$.
\el

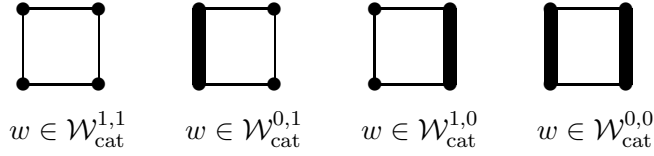
\begin{figure}
\begin{center}
\setlength{\unitlength}{1cm}
\begin{tabular}{c@{\qquad}c@{\qquad}c@{\qquad}c}
   \begin{picture}(1,1.1)(-.1,0)
   \put(-.2,-.2){\framebox(1,1)[c]{}}
   \put(-.2,-.2){\circle*{.17}}
   \put(-.2,.8){\circle*{.17}}
   \put(.8,-.2){\circle*{.17}}
   \put(.8,.8){\circle*{.17}}  
   \end{picture}
& \begin{picture}(1,1.1)(-.1,0)
   \put(-.2,-.2){\framebox(1,1)[c]{}}
   \put(-.2,-.2){\circle*{.17}}
   \put(-.2,.8){\circle*{.17}}
   \put(.8,-.2){\circle*{.17}}
   \put(.8,.8){\circle*{.17}}  
   \linethickness{.17cm}
   \put(-.2,-.2){\line(0,1){1}}
   \end{picture}
& \begin{picture}(1,1.1)(-.1,0)
   \put(-.2,-.2){\framebox(1,1)[c]{}}
   \put(-.2,-.2){\circle*{.17}}
   \put(-.2,.8){\circle*{.17}}
   \put(.8,-.2){\circle*{.17}}
   \put(.8,.8){\circle*{.17}}  
   \linethickness{.17cm}
   \put(.8,-.2){\line(0,1){1}}
   \end{picture}
& \begin{picture}(1,1.1)(-.1,0)
   \put(-.2,-.2){\framebox(1,1)[c]{}}
   \put(-.2,-.2){\circle*{.17}}
   \put(-.2,.8){\circle*{.17}}
   \put(.8,-.2){\circle*{.17}}
   \put(.8,.8){\circle*{.17}}  
   \linethickness{.17cm}
   \put(-.2,-.2){\line(0,1){1}}
   \put(.8,-.2){\line(0,1){1}}
   \end{picture}
\\[12pt]
$w\in\Wi^{1,1}_{\rm cat}$ & $w\in\Wi^{0,1}_{\rm cat}$ 
& $w\in\Wi^{1,0}_{\rm cat}$ &  $w\in\Wi^{0,0}_{\rm cat}$
\end{tabular}
\caption[Effective boundaries for $w\in\Wi_{\rm cat}$.]
{Effective boundaries for $w\in\Wi_{\rm cat}$.}\label{catef}
\end{center}
\end{figure}

For example, for diffusion matrices $w$ from the renormalization class
$\Wi=\Wi_{\rm cat}$, there occur four different effective boundaries,
depending on whether $w\in\Wi^{1,1}_{\rm cat}$, $\Wi^{0,1}_{\rm cat}$,
$\Wi^{1,0}_{\rm cat}$, or $\Wi^{0,0}_{\rm cat}$. These effective
boundaries are depicted in Figure~\ref{catef}. The statement from
Theorem~\ref{main}~(a) that $F_c(\Wi^{l,r}_{\rm cat})\sub\Wi^{l,r}_{\rm cat}$
is just the translation of Lemma~\ref{effinv} to the special set-up there.

{F}rom now on, let $\Wi$ be a renormalization class, i.e., $\Wi$
satisfies also condition~(iv) from Definition~\ref{defi}. Fix a
sequence of (positive) migration constants $(c_k)_{k\geq 0}$. By
definition, the {\em iterated probability kernels} $K^{w,(n)}$
associated with a diffusion matrix $w\in\Wi$ (and the constants
$(c_k)_{k\geq 0}$) are the probability kernels on $\ov D$ defined
inductively by
\be\label{Kdef}
K^{w,(n+1)}_x(\di z):=
\int_{\ov D}\nu^{c_n,F^{(n)}w}_x(\di y)K^{w,(n)}_y(\di z)
\quad(n\geq 0)\quad\mbox{with}\quad K^{w,(0)}_x(\di y):=\de_x(\di y),
\ee
with $F^{(n)}$ as in (\ref{Fn}). Note that
\be\label{KF}
F^{(n)}w(x)=\int_{\ov D} K^{w,(n)}_x(\di y)w(y)\qquad(x\in\ov D,\ n\geq 0).
\ee
The next lemma follows by iteration from Lemmas~\ref{contlem} and
\ref{numom}. It their essence, this lemma
and Lemma~\ref{clusK} below go back to \cite{BCGH95}.
\bl[Basic properties of iterated probability kernels]\label{basic}
For each $w\in\Wi$, the $K^{w,(n)}$ are continuous probability kernels
on $\ov D$. Moreover, for all $x\in\ov D$, $i,j=1,\ldots,d$, and $n\geq 0$,
the mean and covariance matrix of $K^{w,(n)}_x$ are given by
\be\ba{rr@{\;}c@{\;}l}\label{moments}
{\rm (i)}&\dis\int_{\ov D} K^{w,(n)}_x(\di y)(y_i-x_i)&=&0,\\
{\rm (ii)}&\dis\int_{\ov D} K^{w,(n)}_x(\di y)(y_i-x_i)(y_j-x_j)
&=&s_nF^{(n)}w_{ij}(x).
\ec
\el
We equip the space $\Ci(\ov D,\Mi_1(\ov D))$ of continuous probability
kernels on $\ov D$ with the topology of uniform convergence (since
$\Mi_1(\ov D)$ is compact, there is a unique uniform structure on
$\Mi_1(\ov D)$ generating the topology). For `nice' renormalization
classes, it seems reasonable to conjecture that the kernels
$K^{w,(n)}$ converge as $n\to\infty$ to some limit $K^{w,\ast}$ in
$\Ci(\ov D,\Mi_1(\ov D))$. If this happens, then formula
(\ref{moments})~(ii) tells us that the rescaled renormalized diffusion
matrices $s_nF^{(n)}w$ converge uniformly on $\ov D$ to the covariance
matrix of $K^{w,\ast}$. This gives a heuristic explanation why we need
to rescale the iterates $F^{(n)}w$ with the scaling constants $s_n$
from (\ref{sn}) to get a nontrivial limit in (\ref{wast}).

We now explain the relevance of the conditions (\ref{sga})~(i) and
(ii) in the present more general context. If the iterated kernels
converge to a limit $K^{w,\ast}$, then condition~(\ref{sga})~(i)
guarantees that this limit is concentrated on the effective boundary:
\bl[Concentration on the effective boundary]\label{clusK}
If $s_n\asto{n}\infty$, then for any $f\in\Ci(\ov D)$ such that $f=0$ on
$\pa_w D$:
\be\label{concen}
\lim_{n\to\infty}\sup_{x\in\ov D}
\Big|\int_{\ov D}K^{w,(n)}_x(\di y)f(y)\Big|=0.
\ee
\el
In fact, the condition $s_n\to\infty$ guarantees that the
corresponding system of hierarchically interacting diffusions with
migration constants $(c_k)_{k\geq 0}$ {\em clusters in the local mean
field limit}, see \cite[Theorem~3]{DG93a} or Appendix~\ref{hier}
below.

To explain also the relevance of condition~(\ref{sga})~(ii), we
observe that using Lemma~\ref{schaal}, we can convert the rescaled
iterates $s_nF^{(n)}$ into (usual, not rescaled) iterates of another
transformation. For this purpose, it will be convenient to modify the
definition of our scaling constants $s_n$ a little bit. Fix some
$\bet>0$ and put
\be\label{hutsn}
\hut s_n:=\bet+s_n\qquad(n\geq 0).
\ee
Define {\em rescaled renormalization transformations} $\hut F_\ga:\Wi\to\Wi$ by
\be\label{hatFc}
\hut F_\ga w:=(1+\ga)F_{1/\ga}w\qquad(\ga>0,\ w\in\Wi).
\ee
Using (\ref{schafo})~(ii), one easily deduces that
\be\label{sFw}
\hut s_nF^{(n)}w=\hut F_{\ga_{n-1}}\circ\cdots\circ\hut F_{\ga_0}(\bet w)
\qquad(w\in\Wi,\ n\geq 1),
\ee
where
\be\label{gan}
\ga_n:=\frac{1}{\hut s_nc_n}\qquad(n\geq 0).
\ee

We can reformulate the conditions~(\ref{sga})~(i) and (ii) in terms of
the constants $(\ga_n)_{n\geq 0}$. Indeed, it is not hard to
check\footnote{To see this, let $\ov s_\infty\in(0,\infty]$ denote the
limit of the $\ov s_n$ and note that on the one
hand,\label{footnote} $\sum_n1/(\hut s_nc_n)\geq\sum_n\log(1+1/(\hut
s_nc_n))=\log(\prod_n\hut s_{n+1}/\hut s_n)=\log(\hut s_\infty/\hut
s_1)$, while on the other hand $\sum_n 1/(\hut
s_nc_n)\leq\prod_n(1+1/(\hut s_nc_n))=\prod_n\hut s_{n+1}/\hut
s_n=\hut s_\infty/\hut s_1$.} that equivalent formulations of
condition~(\ref{sga})~(i) are:

\be\label{sumga}
{\rm (i)}\quad s_n\asto{n}\infty,\qquad{\rm (ii)}
\quad\hut s_n\asto{n}\infty,\qquad{\rm (iii)}\quad\sum_n\ga_n=\infty.
\ee
Since $\hut s_{n+1}/\hut s_n=1+\ga_n$ we see moreover that, for
any $\ga^\ast\in[0,\infty]$, equivalent formulations of
condition~(\ref{sga})~(ii) are:
\be\label{gamma}
{\rm (i)}\quad\frac{s_{n+1}}{s_n}\asto{n}1+\ga^\ast,
\qquad{\rm (ii)}\quad\frac{\hut s_{n+1}}{\hut s_n}\asto{n}1+\ga^\ast,
\qquad{\rm (iii)}\quad\ga_n\asto{n}\ga^\ast.
\ee
If $0<\ga^\ast<\infty$, then, in the light of (\ref{sFw}), we expect
$\hut s_nF^{(n)}w$ to converge to a fixed point of the transformation
$\hut F_{\ga^\ast}$. If $\ga^\ast=0$, the situation is more complex.
In this case, we expect the orbit $\hut s_nF^{(n)}w\mapsto\hut
s_{n+1}F^{(n+1)}w\mapsto\cdots$, for large $n$, to approximate a
continuous flow, the generator of which is
\be\label{gener}
\lim_{\ga\to 0}\ga^{-1}\Big(\hut F_\ga w-w\Big)(x)
=\ffrac{1}{2}\sum_{i,j=1}^dw_{ij}(x)\difif{x_i}{x_j}w(x)+w(x)\qquad(x\in\ov D).
\ee
To see that the right-hand side of this equation equals the left-hand
side if $w$ is twice continuously differentiable, one needs a Taylor
expansion of $w$ together with the moment formulas (\ref{moments2})
for $\nu^{1/\ga,w}_x$. Under condition condition~(\ref{sumga})~(iii),
we expect this continuous flow to reach equilibrium.

In the light if these considerations, we are led to at the following
general conjecture.
\begin{conjecture}\label{con}{\bf (Limits of rescaled renormalized
diffusion matrices)}
Assume that $s_n\to\infty$ and $s_{n+1}/s_n\to 1+\ga^\ast$ for some
$\ga^\ast\in[0,\infty]$. Then, for any $w\in\Wi$,
\be\label{concon}
s_nF^{(n)}w\asto{n}w^\ast,
\ee
where $w^\ast$ satisfies
\be\ba{rr@{\,}c@{\,}ll}\label{afp}
{\rm (i)}&\dis\hut F_{\ga^\ast}w^\ast&=&\dis w^\ast
\qquad&\mbox{if}\ \ 0<\ga^\ast<\infty,\\[5pt]
{\rm (ii)}&\dis\ffrac{1}{2}\sum_{i,j=1}^d
w^\ast_{ij}(x)\difif{x_i}{x_j}w^\ast(x)+w^\ast(x)&=&\dis 0
\qquad(x\in\ov D)\quad&\mbox{if}\ \ \ga^\ast=0,\\[5pt]
{\rm (iii)}&\dis\lim_{\ga\to\infty}\ov F_\ga w^\ast&=&w^\ast
\qquad&\mbox{if}\ \ \ga^\ast=\infty.
\ec
\end{conjecture}
We call (\ref{afp})~(ii), which is in some sense the $\ga^\ast\to 0$
limit of the fixed point equation (\ref{afp})~(i), the {\em asymptotic
fixed point equation}. A version of formula (\ref{afp})~(ii) occured
in \cite[formula~(1.3.5)]{Swa99} (a minus sign is missing there).

In particular, one may hope that for a given effective boundary, the
equations in (\ref{afp}) have a unique solution. Our main result
(Theorem~\ref{main}) confirms this conjecture for the renormalization
class $\Wi_{\rm cat}$ and for $\ga^\ast<\infty$. In the next section,
we discuss numerical evidence that supports Conjecture~\ref{con} in
the case $\ga^\ast=0$ for other renormalization classes on compacta as
well.

In previous work on renormalization classes, fixed shapes have played
an important role. By definition, for any prerenormalization class
$\Wi$, a {\em fixed shape} is a subclass $\hat\Wi\sub\Wi$ of the form
$\hat\Wi=\{\la w:\la>0\}$ with $0\neq w\in\Wi$, such that
$F_c(\hat\Wi)\sub\hat\Wi$ for all $c>0$. The next lemma describes how
fixed shapes for renormalization classes on compact sets typically
arise.
\bl{\bf(Fixed shapes)}\label{L:shape}
Assume that for each $0<\ga^\ast<\infty$, there is a $0\neq
w^\ast=w^\ast_{\ga^\ast}\in\Wi$ such that
$s_nF^{(n)}w\asto{n}w^\ast_{\ga^\ast}$ whenever $w\in\Wi$,
$s_n\to\infty$, and $s_{n+1}/s_n\to 1+\ga^\ast$. Then:\med

\noi
{\bf (a)} $w^\ast_{\ga^\ast}$ is the unique solution in $\Wi$ of
equation (\ref{afp})~(i).\med

\noi
{\bf (b)} If $w^\ast=w^\ast_{\ga^\ast}$ does not depend on $\ga^\ast$, then
\be\label{fix}
F_c(\la w^\ast)=(\ffrac{1}{\la}+\ffrac{1}{c})^{-1}w^\ast\qquad(\la,c>0).
\ee
Moreover, $\{\la w^\ast:\la>0\}$ is the unique fixed shape in $\Wi$.\med

\noi
{\bf (c)} If the $w^\ast_{\ga^\ast}$ for different values of $\ga^\ast$
are not constant multiples of each other, then $\Wi$ contains no fixed shapes.
\el
Note that by Theorem~\ref{main}, $\Wi^{0,1}_{\rm cat}$ is a
renormalization class satisfying the general assumptions of
Lemma~\ref{L:shape}. The unique solution of (\ref{afp})~(i) in
$\Wi^{0,1}_{\rm cat}$ is of the form $w^\ast=w^{1,p^\ast}$ where
$p^\ast=p^\ast_{0,1,\ga^\ast}$. We conjecture that the
$p^\ast_{0,1,\ga^\ast}$ for different values of $\ga^\ast$ are not
constant multiples of each other, and, as a consequence, that
$\Wi^{0,1}_{\rm cat}$ contains no fixed shapes.

Many facts and conjectures that we have discussed can be generalized
to renormalization classes on unbounded $D$, but in this case, the
second moments of the iterated kernels $K^{w,(n)}$ may diverge as
$n\to\infty$. As a result, because of formula (\ref{moments})~(ii),
the $s_n$ may no longer be the right scaling factors to find a
nontrivial limit of the renormalized diffusion matrices; see, for
example, \cite{BCGH97}.

\subsection{Numerical solutions to the asymptotic fixed point equation}
\label{numer}

Let $t\mapsto w(t,\,\cdot\,)$ be a solution to the continuous flow
with the generator in (\ref{gener}), i.e., $w$ is an $M^d_+$-valued
solution to the nonlinear partial differential equation
\be\label{flow}
\dif{t}w(t,x)=\ffrac{1}{2}\sum_{i,j=1}^dw_{ij}(t,x)\difif{x_i}{x_j}w(t,x)
+w(t,x)\qquad(t\geq 0,\ x\in\ov D).
\ee
Solutions to (\ref{flow}) are quite easy to simulate on a computer. We
have simulated solutions for all kind of diffusion matrices (including
nondiagonal ones) on the unit square $[0,1]^2$, with the effective
boundaries 1--6 depicted in Figure~\ref{fixfig}. For all initial
diffusion matrices $w(0,\,\cdot\,)$ we tried, the solution converged
as $t\to\infty$ to a fixed point $w^\ast$. In all cases except case~6,
the fixed point was unique. The fixed points are listed in
Figure~\ref{fixfig}. The functions $p^\ast_{0,1,0}$ and $q^\ast$ from
Figure~\ref{fixfig} are plotted in Figure~\ref{graph}. Here
$p^\ast_{0,1,0}$ is the function from Theorem~\ref{main}~(c).

The fixed points for the effective boundaries in cases 1,2, and 4 are
the unique solutions of equation (\ref{wiga})~(ii) from
Theorem~\ref{main} in the classes $\Wi^{1,1}_{\rm cat}$,
$\Wi^{0,1}_{\rm cat}$, and $\Wi^{0,0}_{\rm cat}$, respectively. The
simulations suggest that the domain of attraction of these fixed
points (within the class of ``all'' diffusion matrices on $[0,1]^2$)
is actually a lot larger than the classes $\Wi^{1,1}_{\rm cat}$,
$\Wi^{0,1}_{\rm cat}$, and $\Wi^{0,0}_{\rm cat}$.

The function $q^\ast$ from case~3 satisfies $q^\ast(x_1,1)=x_1(1-x_1)$
and is zero on the other parts of the boundary. In contrast to what
one might perhaps guess in view of case~2, $q^\ast$ is {\em not} of
the form $q^\ast(x_1,x_2)=f(x_2)x_1(1-x_1)$ for some function $f$.

Case~5 is somewhat degenerate since in this case the fixed point is
not continuous.

The only case where the fixed point is not unique is case~6. Here, $m$
can be any positive definite matrix, while $g^\ast$, depending on $m$,
is the unique solution on $(0,1)^2$ of the equation
$1+\frac{1}{2}\sum_{i,j=1}^2m_{ij}\difif{x_i}{x_i}g^\ast(x)=0$, with
zero boundary conditions.

\begin{figure}
\begin{center}
\setlength{\unitlength}{1cm}
\begin{tabular}{|c|c|c|}
\hline
case & effective boundary & fixed points $w^\ast$ of (\ref{flow})\\
\hline
1& \begin{picture}(1,1.1)(.2,.2)
   \put(-.2,-.2){\framebox(1,1)[c]{}}
   \put(-.2,-.2){\circle*{.17}}
   \put(-.2,.8){\circle*{.17}}
   \put(.8,-.2){\circle*{.17}}
   \put(.8,.8){\circle*{.17}}  
   \end{picture}
& $\left(\ba{@{}cc@{}}x_1(1-x_1)&0\\0&x_2(1-x_2)\ea\right)$\\
2& \begin{picture}(1,1.1)(.2,.2)
   \put(-.2,-.2){\framebox(1,1)[c]{}}
   \put(-.2,-.2){\circle*{.17}}
   \put(-.2,.8){\circle*{.17}}
   \put(.8,-.2){\circle*{.17}}
   \put(.8,.8){\circle*{.17}}  
   \linethickness{.17cm}
   \put(-.2,-.2){\line(0,1){1}}
   \end{picture}
 & $\left(\ba{@{}cc@{}}x_1(1-x_1)&0\\0&p^\ast_{0,1,0}(x_1)x_2(1-x_2)\ea\right)$\\
3& \begin{picture}(1,1.1)(.2,.2)
   \put(-.2,-.2){\framebox(1,1)[c]{}}
   \put(-.2,-.2){\circle*{.17}}
   \put(-.2,.8){\circle*{.17}}
   \put(.8,-.2){\circle*{.17}}
   \put(.8,.8){\circle*{.17}}  
   \linethickness{.17cm}
   \put(-.2,-.2){\line(0,1){1}}
   \put(-.2,-.2){\line(1,0){1}}
   \end{picture}
 & $\left(\ba{@{}cc@{}}q^\ast(x_1,x_2)&0\\0&q^\ast(x_2,x_1)\ea\right)$\\
4& \begin{picture}(1,1.1)(.2,.2)
   \put(-.2,-.2){\framebox(1,1)[c]{}}
   \put(-.2,-.2){\circle*{.17}}
   \put(-.2,.8){\circle*{.17}}
   \put(.8,-.2){\circle*{.17}}
   \put(.8,.8){\circle*{.17}}  
   \linethickness{.17cm}
   \put(-.2,-.2){\line(0,1){1}}
   \put(.8,-.2){\line(0,1){1}}
   \end{picture}
 & $\left(\ba{@{}cc@{}}x_1(1-x_1)&0\\0&0\ea\right)$\\
5& \begin{picture}(1,1.1)(.2,.2)
   \put(-.2,-.2){\framebox(1,1)[c]{}}
   \put(-.2,-.2){\circle*{.17}}
   \put(-.2,.8){\circle*{.17}}
   \put(.8,-.2){\circle*{.17}}
   \put(.8,.8){\circle*{.17}}  
   \linethickness{.17cm}
   \put(-.2,-.2){\line(0,1){1}}
   \put(-.2,-.2){\line(1,0){1}}
   \put(.8,-.2){\line(0,1){1}}
   \end{picture}
 & $\left(\ba{@{}cc@{}}x_1(1-x_1)1_{\{x_2>0\}}&0\\0&0\ea\right)$\\
6& \begin{picture}(1,1.1)(.2,.2)
   \put(-.2,-.2){\framebox(1,1)[c]{}}
   \put(-.2,-.2){\circle*{.17}}
   \put(-.2,.8){\circle*{.17}}
   \put(.8,-.2){\circle*{.17}}
   \put(.8,.8){\circle*{.17}}  
   \linethickness{.17cm}
   \put(-.2,-.2){\line(0,1){1}}
   \put(-.2,-.2){\line(1,0){1}}
   \put(.8,-.2){\line(0,1){1}}
   \put(-.2,.8){\line(1,0){1}}
   \end{picture}
 & $g^\ast(x_1,x_2)\left(\ba{@{}cc@{}} m_{11}&m_{12}\\m_{21}&m_{22}\ea\right)$\\[14pt]
\hline
\end{tabular}
\caption[Fixed points of the flow (\ref{flow}).]
{Fixed points of the flow (\ref{flow}).}\label{fixfig}
\end{center}
\end{figure}

\begin{figure}
\begin{center}
\includegraphics[width=4cm,height=4cm]{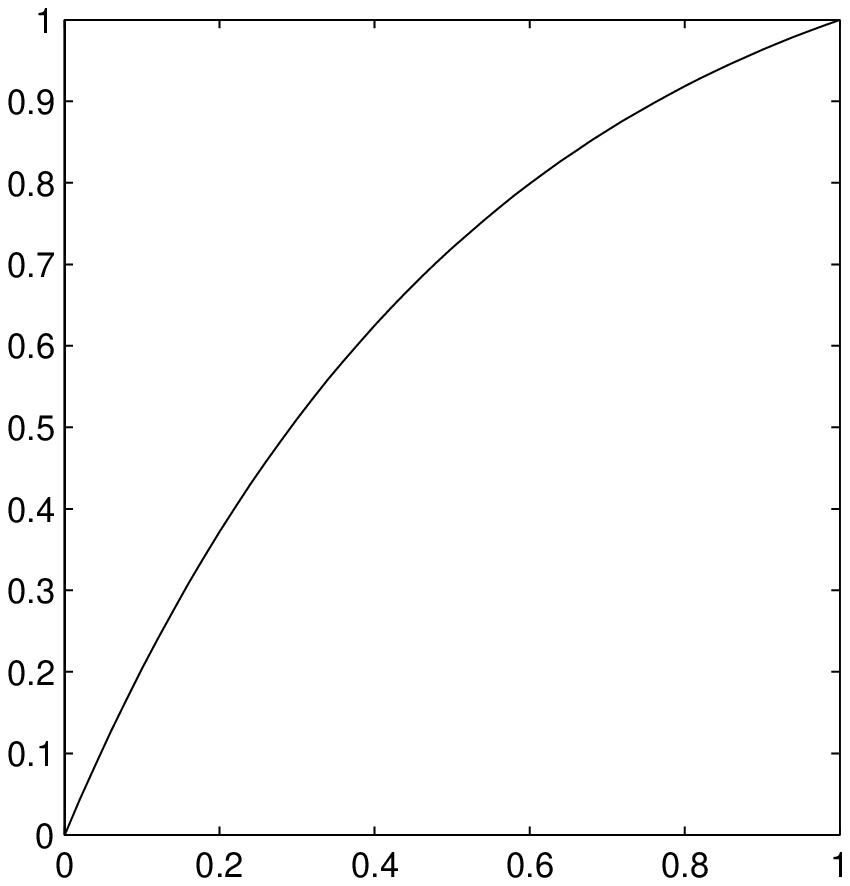}
\hspace{1cm}
\includegraphics[width=5cm,height=4cm]{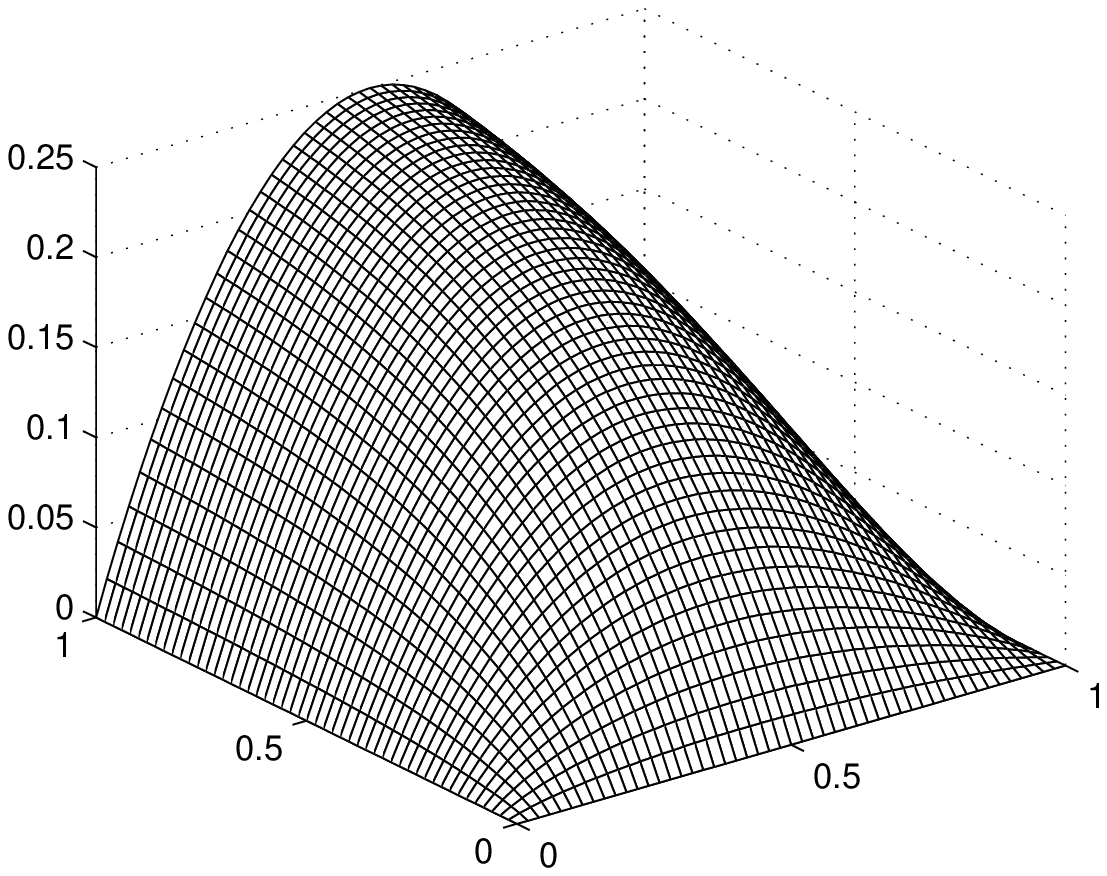}
\caption[The functions $p^\ast_{0,1,0}$ and $q^\ast$ from cases~2 and 3
of Figure~\ref{fixfig}.]{The functions $p^\ast_{0,1,0}$ and $q^\ast$ from
cases~2 and 3 of Figure~\ref{fixfig}.}\label{graph}
\end{center}
\end{figure}

\subsection{Previous rigorous results}\label{prevsec}

In this section we discuss some results that have been derived previously
for renormalization classes on compact sets.
\bt\label{WFTheorem}
{\bf\cite{BCGH95,DGV95} (Universality class of Wright-Fisher models)}
Let $D:=\{x\in\R^d:x_i>0\ \forall i,\ \sum_{i=1}^dx_i<1\}$, and let
$\{e_0,\ldots,e_d\}$, with $e_0:=(0,\ldots,0)$ and
$e_1:=(1,0,\ldots,0),\ldots,\ e_d:=(0,\ldots,0,1)$ be the extremal
points of $\ov D$. Let $w^\ast_{ij}(x):=x_i(\de_{ij}-x_j)$ $(x\in\ov
D$, $i,j=1,\ldots,d)$ denote the standard Wright-Fisher diffusion
matrix, and assume that $\Wi$ is a renormalization class on $\ov D$
such that $w^\ast\in\Wi$ and $\pa_w\ov D=\{e_0,\ldots,e_d\}$ for all
$w\in\Wi$. Let $(c_k)_{k\geq 0}$ be migration constants such that
$s_n\to\infty$ as $n\to\infty$. Then, for all $w\in\Wi$, uniformly on
$\ov D$,
\be\label{toWF}
s_n F^{(n)}w\asto{n} w^\ast.
\ee
\et
The convergence in (\ref{toWF}) is a consequence of Lemmas~\ref{basic}
and \ref{clusK}: The first moment formula (\ref{moments})~(i) and
(\ref{concen}) show that $K^{w,(n)}_x$ converges to the unique
distribution on $\{e_0,\ldots,e_d\}$ with mean $x$, and by the second
moment formula (\ref{moments})~(ii) this implies the convergence of
$s_nF^{(n)}w$.

In order for the iterates in (\ref{toWF}) to be well-defined,
Theorem~\ref{WFTheorem} {\em assumes} that a renormalization class
$\Wi$ of diffusion matrices $w$ on $\ov D$ with effective boundary
$\{e_0,\ldots,e_d\}$ is given. The problem of finding a nontrivial
example of such a renormalization class is open in dimensions greater
than one. In the one-dimensional case, however, the following result
is known.
\bl{\bf \cite{DG93b} (Renormalization class on the unit interval)}
The set
\be\label{WWF}
\Widg:=\{w\in\Ci[0,1]:w=0\mbox{ on }\{0,1\},\ w>0\mbox{ on }(0,1),
\ w\mbox{ Lipschitz}\}
\ee
is a renormalization class on $[0,1]$.
\el
About renormalization of isotropic diffusions, the following result
is known. Below, $\pa D:=\ov D\beh D$ denotes the topological boundary of $D$.
\bt\label{isoTheorem}
{\bf\cite{HS98} (Universality class of isotropic models)}
Let $D\sub\R^d$ be open, bounded, and convex and let $m\in M^d_+$ be
fixed and (strictly) positive definite. Set
$w^\ast_{ij}(x):=m_{ij}g^\ast(x)$, where $g^\ast$ is the unique
solution of $1+\ffrac{1}{2}\sum_{ij}m_{ij}\difif{x_i}{x_j}g^\ast(x)=0$
for $x\in D$ and $g^\ast(x)=0$ for $x\in\pa D$. Assume that $\Wi$ is a
renormalization class on $\ov D$ such that $w^\ast\in\Wi$ and such
that each $w\in\Wi$ is of the form
\be
w_{ij}(x)=m_{ij}g(x)\qquad(x\in\ov D,\ i,j=1,\ldots,d),
\ee
for some $g\in\Ci(\ov D)$ satisfying $g>0$ on $D$ and $g=0$ on $\pa
D$. Let $(c_k)_{k\geq 0}$ be migration constants such that
$s_n\to\infty$ as $n\to\infty$. Then, for all $w\in\Wi$, uniformly on
$\ov D$,
\be\label{togast}
s_n F^{(n)}w\asto{n}w^\ast.
\ee
\et
The proof of Theorem~\ref{isoTheorem} follows the same lines as the
proof of Theorem~\ref{WFTheorem}, with the difference that in this
case one needs to generalize the first moment formula
(\ref{moments})~(i) in the sense that $\int_{\ov D} K^{w,(n)}_x(\di
y)h(y)=h(x)$ for any $m$-harmonic function $h$, i.e., $h\in\Ci(\ov D)$
satisfying $\sum_{ij}m_{ij}\difif{x_i}{x_j}h(x)=0$ for $x\in D$. The
kernel $K^{w,(n)}_x$ now converges to the $m$-harmonic measure on $\pa
D$ with mean $x$, and this implies (\ref{togast}).

Again, in dimensions $d\geq 2$, the problem of finding a `reasonable'
class $\Wi$ satisfying the assumptions of Theorem~\ref{isoTheorem} is
so far unresolved. The problem with verifying conditions (i)--(iv)
from Definition~\ref{defi} in an explicit set-up is that (i) and (ii)
usually require some smoothness of $w$, while (iv) requires that one
can prove the same smoothness for $F_cw$, which is difficult.

The proofs of Theorems~\ref{WFTheorem} and \ref{isoTheorem} are based
on the same principle. For any diffusion matrix $w$, let $H_w$ denote
the class of $w$-harmonic functions, i.e., functions $h\in\Ci(\ov D)$
satisfying $\sum_{ij}w_{ij}(x)\difif{x_i}{x_j}h(x)=0$ on $D$. If $w$
belongs to one of the renormalization classes in
Theorems~\ref{WFTheorem} and \ref{isoTheorem}, then $H_w$ has the
property that $T^c_{x,t}h(H_w)\sub H_w$ for all $c>0$, $x\in\ov D$,
and $t\geq 0$, where $T^c_{x,t}h(y):= h(x+(y-x)e^{-ct})$ is the
semigroup with generator $\sum_{i=1}^dc(x_i-y_i)\dif{y_i}$, i.e., the
operator in (\ref{Adef}) without the diffusion part.  In this case we
say that $w$ has {\em invariant harmonics}; see \cite{Swa00}. As a
consequence, one can prove that the iterated kernels satisfy
$\int_{\ov D} K^{w,(n)}_x(\di y)h(y)=h(x)$ for all $h\in H_w$ and
$x\in\ov D$. If $s_n\to\infty$, then this implies that $K^{w,(n)}_x$
converges to the unique $H_w$-harmonic measure on $\pa_wD$ with mean
$x$. Diffusion matrices from $\Wi_{\rm cat}$ do not in general have
invariant harmonics. Therefore, to prove Theorem~\ref{main}, we need
new techniques.

Note that in the renormalization classes from Theorems~\ref{WFTheorem}
and \ref{isoTheorem}, the unique attraction point $w^\ast$ does not
depend on $\ga^\ast$. Therefore, by Lemma~\ref{L:shape}, these
renormalization classes contain a unique fixed shape, which is given
by $\{\la w^\ast:\la>0\}$.

\section{Connection with branching theory}\label{bracon}

{F}rom now on, we focuss on the renormalization class $\Wi_{\rm cat}$.
We will show that for this renormalization class, the rescaled
renormalization transformations $\ov F_\ga$ from (\ref{hatFc}) can be
expressed in terms of the log-Laplace operators of a discrete time
branching process on $[0,1]$. This will allow us to use techniques
from the theory of spatial branching processes to verify
Conjecture~\ref{con} for the renormalization class $\Wi_{\rm cat}$ in
the case $\ga^\ast<\infty$.

\subsection{Poisson-cluster branching processes}

We first need some concepts and facts from branching theory. Finite
measure-valued branching processes (on $\R$) in discrete time have
been introduced by Ji\v rina \cite{Jir64}. We need to consider only a
special class. Let $E$ be a separable, locally compact, and metrizable
space. We call a continuous map $\Qi$ from $E$ into $\Mi_1(\Mi(E))$ a
{\em continuous cluster mechanism}. By definition, an $\Mi(E)$-valued
random variable $\Xc$ is a {\em Poisson cluster measure} on $E$ with
locally finite {\em intensity measure} $\mu$ and continuous cluster
mechanism $\Qi$, if its log-Laplace transform satisfies
\be\label{randclust}
-\log E\big[\ex{-\li\Xc,f\re}\big]=\int_E\!\mu(\di x)\Big(1-\int_{\Mi(E)}\!\!\Qi(x,\di\chi)\ex{-\li\chi,f\re}\Big)\quad(f\in B_+(E)).
\ee
For given $\mu$ and $\Qi$, such a Poisson cluster measure exists, and
is unique in distribution, provided that the right-hand side of
(\ref{randclust}) is finite for $f=1$. It may be constructed as
$\Xc=\sum_i\chi_{x_i}$, where $\sum_i\de_{x_i}$ is a (possibly
infinite) Poisson point measure with intensity $\mu$, and given
$x_1,x_2,\ldots$, the $\chi_{x_1},\chi_{x_2},\ldots$ are independent
random variables with laws
$\Qi(x_1,\,\cdot\,),\Qi(x_2,\,\cdot\,),\ldots$, respectively.

\detail{Choose compact sets $(C_n)_{n\geq 0}$ with $C_n\up E$ and set
$B_0:=C_0$ and $B_n:=C_n\beh C_{n-1}$ $(n\geq 1)$. Then the measures
$\mu_n(\,\cdot\,):=\mu(\,\cdot\,\cap C_n)$ are finite, and we can
construct independent Poisson cluster measures $\Xc_n$ with
intensity measure $\mu_n$ and cluster mechanism $\Qi$. Set
$\Xc:=\sum_{n=0}^\infty\Xc_n$. Since for each $n$,
\[
-\log E\big[\ex{-\li\Xc_n,f\re}\big]
=\int_{B_n}\!\mu(\di x)\Big(1-\int_{\Mi(E)}\!\!\Qi(x,\di\chi)
\ex{-\li\chi,f\re}\Big)\quad(f\in B_+(E)),
\]
we see that
\[\ba{l}
\dis-\log E\big[\ex{-\li\Xc,f\re}\big]
=-\log\prod_n E\big[\ex{-\li\Xc_n,f\re}\big]\\[5pt]
\dis\quad=-\sum_n\log E\big[\ex{-\li\Xc_n,f\re}\big]
=\sum_n\int_{B_n}\!\mu(\di x)\Big(1-\int_{\Mi(E)}\!\!\Qi(x,\di\chi)
\ex{-\li\chi,f\re}\Big)\\[5pt]
\dis\quad=\int_E\!\mu(\di x)\Big(1-\int_{\Mi(E)}\!\!\Qi(x,\di\chi)
\ex{-\li\chi,f\re}\Big)\qquad(f\in B_+(E)),
\ea\]
so $\Xc$ satisfies (\ref{randclust}). To show that $\Xc$ is a.s.\
finite it suffices to show that $\Pois(\Xc)$ is a.s.\ finite. Since
$\Pois(\Xc)=\sum_n\Pois(\Xc_n)$ where the $\Pois(\Xc_n)$ are
independent, it suffices to show that a.s.\ only finitely many of the
$\Pois(\Xc_n)$ are nonzero. Set
\[
K_n:=\int_{B_n}\!\mu(\di x)\Big(1-\int_{\Mi(E)}\!\!\Qi(x,\di\chi)
\ex{-\li\chi,1\re}\Big).
\]
Since we are assuming that the right-hand side of (\ref{randclust})
is finite for $f=1$, we have
\[
\sum_nK_n<\infty.
\]
Now
\[
P[\Pois(\Xc_n)\neq 0]=1-E\big[\ex{-|\Xc_n|}\big]=1-\ex{-K_n}\leq K_n,
\]
so
\[
\sum_nP[\Pois(\Xc_n)\neq 0]<\infty,
\]
which by Borel-Cantelli implies that only finitely many of the
$\Pois(\Xc_n)$ are nonzero.}

Now fix a finite sequence of functions $q_k\in\Ci_+(E)$ and continuous
cluster mechanisms $\Qi_k$ ($k=1,\ldots,n$), define
\be\label{Vk}
\Ui_kf(x):=q_k(x)\Big(1-\int_{\Mi(E)}\!\!\Qi_k(x,\di\chi)\ex{-\li\chi,f\re}\Big)\qquad(x\in E,\ f\in B_+(E),\ k=1,\ldots,n),
\ee
and assume that
\be\label{fincon}
\sup_{x\in E}\Ui_k1(x)<\infty\qquad(k=1,\ldots,n).
\ee
Then $\Ui_k$ maps $B_+(E)$ into $B_+(E)$ for each $k$, and for each
$\Mi(E)$-valued initial state $\Xc_0$, there exists a
(time-inhomogeneous) Markov chain $(\Xc_0,\ldots,\Xc_n)$ in $\Mi(E)$,
such that $\Xc_k$, given $\Xc_{k-1}$, is a Poisson cluster measure
with intensity $q_k\Xc_{k-1}$ and cluster mechanism $\Qi_k$.
\detail{We claim that condition (\ref{fincon}) implies that a Poisson
cluster measure with intensity $q_k\mu$ and cluster mechanism
$\Qi_k$ exists and is a.s.\ finite for each finite $\mu$. Note that
we do not insist that $q_k$ is bounded and therefore $q_k\mu$ need
not be finite. By our earlier remarks, it suffices to check that
\[
\int_E\!q_k(x)\mu(\di x)\Big(1-\int_{\Mi(E)}\!\!\Qi(x,\di\chi)
\ex{-\li\chi,1\re}\Big)<\infty
\]
for each finite measure $\mu$. For this, it clearly suffices that
\[
x\mapsto q_k(x)\Big(1-\int_{\Mi(E)}\!\!\Qi(x,\di\chi)\ex{-\li\chi,1\re}\Big)
\]
is bounded, which is just condition (\ref{fincon}).}
It is not hard to see that
\be\label{Vi}
E^\mu\big[\ex{-\li\Xc_n,f\re}\big]
=\ex{-\li\mu,\Ui_1\circ\cdots\circ\Ui_n f\re}
\qquad(\mu\in\Mi(E),\ f\in B_+(E)).
\ee
\detail{It is not hard to check that $\Ui_k(rf)\leq r\Ui(f)$ for all
$r\geq 1$) (compare Lemma~\ref{rhlem}), so $\Ui_k(r)\leq r\Ui_k(1)$.
Together with the monotonicity of $\Ui_k$ this shows that under the
condition (\ref{fincon}), $\Ui_k$ maps $B_+(E)$ into $B_+(E)$. The
proof of (\ref{Vi}) is standard.}
We call $\Xc=(\Xc_0,\ldots,\Xc_n)$ the {\em Poisson-cluster branching
process} on $E$ with {\em weight functions} $q_1,\ldots,q_n$ and
cluster mechanisms $\Qi_1,\ldots,\Qi_n$. The operator $\Ui_k$ is
called the {\em log-Laplace operator} of the transition law from
$\Xc_{k-1}$ to $\Xc_k$. Note that we can write (\ref{Vi}) in the
suggestive form
\be\label{Vi2}
P^\mu\big[\Pois(f\Xc_n)=0\big]
=P\big[\Pois\big((\Ui_1\circ\cdots\circ\Ui_n f)\mu\big)=0\big].
\ee
Here, if $\mu$ is an $\Mi(E)$-valued random variable, then
$\Pois(\mu)$ denotes an $\Ni(E)$-valued random variable such that
conditioned on $\mu$, $\Pois(\mu)$ is a Poisson point measure with
intensity $\mu$.

\subsection{The renormalization branching process}\label{RBP}

We will now construct a Poisson-cluster branching process on $[0,1]$
of a special kind, and show that the rescaled renormalization
transformations on $\Wi_{\rm cat}$ can be expressed in terms of the
log-Laplace operators of this branching process.

By Lemma~\ref{ergo} below, for each $\ga>0$ and $x\in[0,1]$, the SDE
\be\label{Yclx}
\di\y(t)=\ffrac{1}{\ga}\,(x-\y(t))\di t+\sqrt{2\y(t)(1-\y(t))}\di B(t),
\ee
has a unique (in law) stationary solution. We denote this solution by
$(\y^\ga_x(t))_{t\in\R}$. Let $\tau_\ga$ be an independent exponentially
distributed random variable with mean $\ga$, and set
\be\label{Zidef}
\Zi^\ga_x:=\int_0^{\tau_\ga}\de_{\y^\ga_x(-t/2)}\di t\qquad(\ga>0,\ x\in[0,1]).
\ee
Define constants $q_\ga$ and continuous (by Corollary~\ref{Ugacont} below)
cluster mechanisms $\Qi_\ga$ by
\be\label{Qqdef}
q_\ga:=\ffrac{1}{\ga}+1\qquad\mbox{and}\qquad\Qi_\ga(x,\,\cdot\,)
:=\Li(\Zi^\ga_x)\qquad(\ga>0,\ x\in[0,1]),
\ee
and let $\Ui_\ga$ denote the log-Laplace operator with (constant)
weight function $q_\ga$ and cluster mechanism $\Qi_\ga$, i.e.,
\be\label{Uiga}
\Ui_\ga f(x):=q_\ga\Big(1-\int_{\Mi([0,1])}\!\!\Qi_\ga(x,\di\chi)
\ex{-\li\chi,f\re}\Big)\qquad(x\in[0,1],\ f\in B_+[0,1],\ \ga>0).
\ee
We now establish the connection between renormalization transformations
on $\Wi_{\rm cat}$ and log-Laplace operators.
\bp{\bf (Identification of the renormalization transformation)}\label{identprop}
Let $\hut F_\ga$ be the rescaled renormalization transformation on
$\Wi_{\rm cat}$ defined in (\ref{hatFc}). Then
\be\label{hutform}
\hut F_\ga w^{\txt 1,p}=w^{\txt 1,\Ui_\ga p}\qquad(p\in\Hi,\ \ga>0).
\ee
\ep
Fix a diffusion matrix $w^{\al,p}\in\Wi_{\rm cat}$ and migration constants
$(c_k)_{k\geq 0}$. Define constants $\hut s_n$ and $\ga_n$ as in (\ref{hutsn})
and (\ref{gan}), respectively, where $\bet:=1/\al$. Then
Proposition~\ref{identprop} and formula (\ref{sFw}) show that
\be\label{renbra}
\hut s_nF^{(n)}w^{\al,p}=w^{\txt 1,\Ui_{\ga_{n-1}}\circ\cdots
\circ\Ui_{\ga_0}(\frac{p}{\al})}.
\ee
Here $\Ui_{\ga_{n-1}},\ldots,\Ui_{\ga_0}$ are the log-Laplace
operators of the Poisson-cluster branching process
$\Xc=(\Xc_{-n},\ldots,\Xc_0)$ with weight functions
$q_{\ga_{n-1}},\ldots,q_{\ga_0}$ and cluster mechanisms
$\Qi_{\ga_{n-1}},\ldots,\Qi_{\ga_0}$. We call $\Xc$ (started at some
time $-n$ in an initial law $\Li(\Xc_{-n})$) the {\em renormalization
  branching process}. By formulas (\ref{Vi}) and (\ref{renbra}), the
study of the limiting behavior of rescaled iterated renormalization
transformations on $\Wi_{\rm cat}$ reduces to the study of the
renormalization branching process $\Xc$ in the limit $n\to\infty$.

\subsection{Convergence to a time-homogeneous process}

Let $\Xc=(\Xc_{-n},\ldots,\Xc_0)$ be the renormalization branching
process introduced in the last section. If the constants
$(\ga_k)_{k\geq 0}$ satisfy $\sum_n\ga_n=\infty$ and
$\ga_n\to\ga^\ast$ for some $\ga^\ast\in\half$, then $\Xc$ is almost
time-homo\-geneous for large $n$. More precisely, we will prove the
following convergence result.
\bt{\bf(Convergence to a time-homogenous branching process)}\label{supercon}
Assume that $\Li(\Xc_{-n})\Asto{n}\mu$ for some probability law $\mu$
on $\Mi([0,1])$.\med

\noi
{\bf (a)} If $0<\ga^\ast<\infty$, then
\be\label{disco}
\Li(\Xc_{-n},\Xc_{-n+1},\ldots)\Asto{n}
\Li(\Yi^{\ga^\ast}_0,\Yi^{\ga^\ast}_1,\ldots),
\ee
where $\Yi^{\ga^\ast}$ is the time-homogenous branching process with
log-Laplace operator $\Ui_{\ga^\ast}$ in each step and initial law
$\Li(\Yi^{\ga^\ast}_0)=\mu$.\med

\noi
{\bf (b)} If $\ga^\ast=0$, then
\be\label{coco}
\Li\Big(\big(\Xc_{-k_n(t)}\big)_{t\geq 0}\Big)\Asto{n}
\Li\Big(\big(\Yi^0_t\big)_{t\geq 0}\Big),
\ee
where $\Rightarrow$ denotes weak convergence of laws on path space,
$k_n(t):=\min\{k:0\leq k\leq n$, $\sum_{l=k}^{n-1}\ga_l\leq t\}$, and
$\Yi^0$ is the super-Wright-Fisher diffusion with activity and growth
parameter both identically $1$ and initial law $\Li(\Yi^0_0)=\mu$.
\et
The super-Wright-Fisher diffusion was studied in \cite{FSsup}. By
definition, $\Yi^0$ is the time-homo\-ge\-neous Markov process in
$\Mi[0,1]$ with continuous sample paths, whose Laplace functionals are
given by
\be
E^\mu\big[\ex{-\li\Yi^0_t,f\re}\big]=\ex{-\li\mu,\Ui^0_tf\re}
\qquad(\mu\in\Mi[0,1],\ f\in B_+[0,1],\ t\geq 0).
\ee
Here $\Ui^0_tf=u_t$ is the unique mild solution of the semilinear
Cauchy equation
\be
\left\{\ba{r@{\,}c@{\,}l}\label{cau}
\dif{t}u_t(x)&=&\ffrac{1}{2}x(1-x)\diff{x}u_t(x)+u_t(x)(1-u_t(x))
\quad(t\geq 0,\ x\in[0,1]),\\
u_0&=&f.\ea\right.
\ee
For a further study of the renormalization branching process $\Xc$
and its limiting processes $\Yi^{\ga^\ast}$ ($\ga^\ast\geq 0$) we
will use the technique of embedded particle systems, which we explain
in the next section.

\subsection{Weighted and Poissonized branching processes}

In this section, we explain how from a Poisson-cluster branching
process it is possible to construct other branching processes by
weighting and Poissonization. We first need to introduce spatial
branching particle systems in some generality.

Let $E$ again be separable, locally compact, and metrizable. For
$\nu\in\Ni(E)$ and $f\in B_{[0,1]}(E)$, we adopt the notation
\be
f^{\txt 0}:=1\quad\mbox{and}\quad f^{\txt\nu}:=\prod_{i=1}^mf(x_i)
\quad\mbox{when}\quad\nu=\sum_{i=1}^m\de_{x_i}\quad(m\geq 1).
\ee
We call a continuous map $x\mapsto Q(x,\,\cdot\,)$ from $E$ into
$\Mi_1(\Ni(E))$ a {\em continuous offspring mechanism.}

Fix continuous offspring mechanisms $Q_k$ ($1\leq k\leq n$), and
let $(X_0,\ldots,X_n)$ be a Markov chain in $\Ni(E)$ such that,
given that $X_{k-1}=\sum_{i=1}^m\de_{x_i}$, the next step of the
chain $X_k$ is a sum of independent random variables with laws
$Q_k(x_i,\,\cdot\,)$ ($i=1,\ldots,m$). Then
\be\label{V}
E^\nu\big[(1-f)^{\txt X_n}\big]=(1-U_1\circ\cdots\circ U_nf)^{\txt\nu}
\qquad(\nu\in\Ni(E),\ f\in B_{[0,1]}(E)),
\ee
where $U_k:B_{[0,1]}(E)\to B_{[0,1]}(E)$ is defined as
\be\label{VV}
U_kf(x):=1-\int_{\Ni(E)}\!\!Q^k(x,\di\nu)(1-f)^{\txt\nu}
\qquad(1\leq k\leq n,\ x\in E,\ f\in B_{[0,1]}(E)).
\ee
We call $U_k$ the {\em generating operator} of the transition law from
$X_{k-1}$ to $X_k$, and we call $X=(X_0,\ldots,X_n)$ the {\em
branching particle system} on $E$ with generating operators
$U_1,\ldots,U_n$. It is often useful to write (\ref{V}) in the
suggestive form
\be\label{V2}
P^\nu\big[\Thin_f(X_n)=0\big]=P\big[\Thin_{U_1\circ\cdots\circ U_n f}(\nu)
=0\big]\qquad(\nu\in\Ni(E),\ f\in B_{[0,1]}(E)).
\ee
Here, if $\nu$ is an $\Ni(E)$-valued random variable and $f\in
B_{[0,1]}(E)$, then $\Thin_f(\nu)$ denotes an $\Ni(E)$-valued random
variable such that conditioned on $\nu$, $\Thin_f(\nu)$ is obtained
from $\nu$ by independently throwing away particles from $\nu$, where
a particle at $x$ is kept with probability $f(x)$. One has the
elementary relations
\be\label{elrel}
\Thin_f(\Thin_g(\nu))\isd\Thin_{fg}(\nu)\quad\mbox{and}
\quad\Thin_f(\Pois(\mu))\isd\Pois(f\mu),
\ee
where $\isd$ denotes equality in distribution. 

We are now ready to describe weighted and Poissonized branching
processes. Let $\Xc=(\Xc_0,\ldots,\Xc_n)$ be a Poisson-cluster
branching process on $E$, with continuous weight functions
$q_1,\ldots,q_n$, continuous cluster mechanisms $\Qi_1,\ldots,\Qi_n$,
and log-Laplace operators $\Ui_1,\ldots,\Ui_n$ given by (\ref{Vk}) and
satisfying (\ref{fincon}). Let $\Zi^k_x$ denote an $\Mi(E)$-valued
random variable with law $\Qi_k(x,\,\cdot\,)$. Let $h\in\Ci_+(E)$ be
bounded, $h\neq 0$, and put $E^h:=\{x\in E:h(x)>0\}$. For $f\in
B_+(E^h)$, define $hf\in B_+(E)$ by $hf(x):=h(x)f(x)$ if $x\in E^h$
and $hf(x):=0$ otherwise.
\bp{\bf(Weighting of Poisson-cluster branching processes)}\label{weightprop}
Assume that there exists a constant $K<\infty$ such that $\Ui_kh\leq
Kh$ for all $k=1,\ldots,n$. Then there exists a Poisson-cluster
branching process $\Xc^h=(\Xc^h_0,\ldots,\Xc^h_n)$ on $E^h$ with
weight functions $(q_1^h,\ldots,q^h_n)$ given by $q^h_k:=q_k/h$,
continuous cluster mechanisms $\Qi^h_1,\ldots,\Qi^h_n$ given by
\be\label{Qih}
\Qi^h_k(x,\,\cdot\,):=\Li(h\Zi^k_x)\qquad(x\in E^h),
\ee
and log-Laplace operators $\Ui^h_1,\ldots,\Ui^h_n$ satisfying
\be\label{htrafo1}
h\,\Ui^h_kf:=\Ui_k(hf)\qquad(f\in B_+(E^h)).
\ee
The processes $\Xc$ and $X^h$ are related by
\be\label{weighting}
\Li(\Xc^h_0)=\Li(h\Xc_0)\quad\mbox{implies}\quad\Li(\Xc^h_k)
=\Li(h\Xc_k)\qquad(0\leq k\leq n).
\ee
\ep
\bp{\bf(Poissonization of Poisson-cluster branching processes)}\label{Poisprop}
Assume that $\Ui_kh\leq h$ for all $k=1,\ldots,n$. Then there exists
a branching particle system $X^h=(X^h_0,\ldots,X^h_n)$ on $E^h$ with
continuous offspring mechanisms $Q^h_1,\ldots,Q^h_n$ given by
\be\label{Zxdef}
Q^h_k(x,\,\cdot\,):=\frac{q_k(x)}{h(x)}P\big[\Pois(h\Zi^k_x)\in\cdot\,\big]
+\Big(1-\frac{q_k(x)}{h(x)}\Big)\de_0(\,\cdot\,)
\qquad(x\in E^h),
\ee
and generating operators $U^h_1,\ldots,U^h_n$ satisfying
\be\label{htrafo}
hU^h_kf:=\Ui_k(hf)\qquad(f\in B_{[0,1]}(E^h)).
\ee
The processes $\Xc$ and  $X^h$ are related by
\be\label{Poissonization}
\Li(X^h_0)=\Li(\Pois(h\Xc_0))\quad\mbox{implies}\quad\Li(X^h_k)
=\Li(\Pois(h\Xc_k))\qquad(0\leq k\leq n).
\ee
\ep
Here, the right-hand side of (\ref{Zxdef}) is always a probability
measure, despite that it may happen that $q_k(x)/h(x)>1$. The
(straightforward) proofs of Propositions~\ref{weightprop} and
\ref{Poisprop} can be found in Section~\ref{Poissec} below. If
(\ref{weighting}) holds then we say that $\Xc^h$ is obtained from
$\Xc$ by {\em weighting} with density $h$. If (\ref{Poissonization})
holds then we say that $X^h$ is obtained from $\Xc$ by {\em
Poissonization} with density $h$. Proposition~\ref{Poisprop} says
that a Poisson-cluster branching process $\Xc$ contains, in a way,
certain `embedded' branching particle systems $X^h$. Poissonization
relations for superprocesses and embedded particle systems have
enjoyed considerable attention, see \cite{FStrim} and references
therein.

A function $h\in B_+(E)$ such that $\Ui_kh\leq h$ is called {\em
$\Ui_k$-super\-harmonic}. If the reverse inequality holds we say
that $h$ is {\em $\Ui_k$-sub\-harmonic}. If $\Ui_kh=h$ then $h$ is
called {\em $\Ui_k$-harmonic}.

\subsection{Extinction versus unbounded growth for embedded particle
systems}\label{exex}

In this section we explain how embedded particle systems can be used
to prove Theorem~\ref{main}. Throughout this section $(\ga_k)_{k\geq 0}$
are positive constants such that $\sum_n\ga_n=\infty$ and
$\ga_n\to\ga^\ast$ for some $\ga^\ast\in\half$, and
$\Xc=(\Xc_{-n},\ldots,\Xc_0)$ is the renormalization branching process
on $[0,1]$ defined in Section~\ref{RBP}. We write
\be
\Ui^{(n)}:=\Ui_{\ga_{n-1}}\circ\cdots\circ\Ui_{\ga_0}.
\ee
In view of formula (\ref{renbra}), in order to prove Theorem~\ref{main},
we need the following result.
\bp{\bf (Limits of iterated log-Laplace operators)}\label{Uit}
Uniformly on $[0,1]$,
\be\ba{rr@{\,}c@{\,}ll}\label{Ucon}
{\rm (i)}&\dis\lim_{n\to\infty}\Ui^{(n)}p&=&1&
\dis\qquad(p\in\Hi_{1,1}),\\[5pt]
{\rm (ii)}&\dis\lim_{n\to\infty}\Ui^{(n)}p&=&0&
\dis\qquad(p\in\Hi_{0,0}),\\[5pt]
{\rm (iii)}&\dis\lim_{n\to\infty}\Ui^{(n)}p&=&p^{\ast}_{0,1,\ga^\ast}&
\dis\qquad(p\in\Hi_{0,1}),
\ec
where $p^{\ast}_{0,1,\ga^\ast}:[0,1]\to[0,1]$ is a function depending
on $\ga^\ast$ but not on $p\in\Hi_{0,1}$.
\ep
In our proof of Proposition~\ref{Uit}, we will use embedded particle
systems $X^h=(X^h_{-n},\ldots,X^h_0)$ obtained from $\Xc$ by Poissonization
with certain $h$ taken from the classes $\Hi_{1,1}$, $\Hi_{0,0}$, and
$\Hi_{0,1}$.
\bl{\bf(Embedded particle system with $h_{1,1}$)}\label{11lem}
The constant function $h_{1,1}(x):=1$ is $\Ui_\ga$-harmonic for each
$\ga>0$. The corresponding embedded particle system $X^{h_{1,1}}$ on
$[0,1]$ satisfies
\be\label{expl}
P^{-n,\de_x}\big[|X^{h_{1,1}}_0|\in\cdot\,\big]\Asto{n}\de_\infty
\ee
uniformly\footnote{Since $\Mi_1[0,\infty]$ is compact in the topology
of weak convergence, there is a unique uniform structure compatible
with the topology, and therefore it makes sense to talk about uniform
convergence of $\Mi_1[0,\infty]$-valued functions (in this case,
$x\mapsto P^{-n,\de_x}\big[|X^{h_{1,1}}_0|\in\cdot\,\big]$).} for
all $x\in[0,1]$.
\el
In (\ref{expl}) and similar formulas below, $\Rightarrow$ denotes
weak convergence of probability measures on $[0,\infty]$. Thus,
(\ref{expl}) says that for processes started with one particle
on the position $x$ at times $-n$, the number of particles at
time zero converges to infinity as $n\to\infty$.
\bl{\bf(Embedded particle system with $h_{0,0}$)}\label{00lem}
The function $h_{0,0}(x):=x(1-x)$ $(x\in[0,1])$ is
$\Ui_\ga$-super\-harmonic for each $\ga>0$. The corresponding
embedded particle system $X^{h_{0,0}}$ on $(0,1)$ is critical
and satisfies
\be\label{ext}
P^{-n,\de_x}\big[|X^{h_{0,0}}_0|\in\cdot\,\big]\Asto{n}\de_0
\ee
locally uniformly for all $x\in(0,1)$.
\el
Here, a branching particle system $X$ is called {\em critical} if each
particle produces on average one offspring (in each time step and
independent of its position). Formula (\ref{ext}) says that the
embedded particle system $X^{h_{0,0}}$ gets extinct during the time
interval $\{-n,\ldots,0\}$ with probability tending to one as
$n\to\infty$. We can summarize Lemmas~\ref{11lem} and \ref{00lem} by
saying that the embedded particle system associated with $h_{1,1}$
grows unboundedly while the embedded particle system associated with
$h_{0,0}$ becomes extinct as $n\to\infty$.

We will also consider an embedded particle systems $X^{h_{0,1}}$ for a
certain $h_{0,1}$ taken from $\Hi_{0,1}$. It turns out that this
system either gets extinct or grows unboundedly, each with a positive
probability. In order to determine these probabilities, we need to
consider embedded particle systems for the time-homogeneous processes
$\Yi^{\ga^\ast}$ ($\ga^\ast\in\half$) from (\ref{disco}) and
(\ref{coco}). If $h\in\Hi_{0,1}$ is $\Ui_{\ga^\ast}$-superharmonic for
some $\ga^\ast>0$, then Poissonizing the process $\Yi^{\ga^\ast}$ with
$h$ yields a branching particle system on $(0,1]$ which we denote by
$Y^{\ga^\ast,h}=(Y^{\ga^\ast,h}_0,Y^{\ga^\ast,h}_1,\ldots)$.
Likewise, if $h\in\Hi_{0,1}$ is twice continuously differentiable and
satisfies
\be\label{contsup}
\ffrac{1}{2}x(1-x)\diff{x}h(x)-h(x)(1-h(x))\leq 0,
\ee
then Poissonizing the super-Wright-Fisher diffusion $\Yi^0$ with $h$
yields a continuous-time branching particle system on $(0,1]$, which
we denote by $Y^{0,h}=(Y^{0,h}_t)_{t\geq 0}$. For example, for $m\geq 4$,
the function $h(x):=1-(1-x)^m$ satisfies (\ref{contsup}).
\bl{\bf(Embedded particle system with $h_{0,1}$)}\label{01lem}
The function $h_{0,1}(x):=1-(1-x)^7$ is $\Ui_\ga$-super\-harmonic for
each $\ga>0$. The corresponding embedded particle system $X^{h_{0,1}}$
on $(0,1]$ satisfies
\be\label{expext}
P^{-n,\de_x}\big[|X^{h_{0,1}}_0|\in\cdot\,\big]\Asto{n}
\rho_{\ga^\ast}(x)\de_\infty+(1-\rho_{\ga^\ast}(x))\de_0,
\ee
locally uniformly for all $x\in(0,1]$, where
\be\label{rhodef}
\rho_{\ga^\ast}(x):=\left\{\ba{ll}\dis P^{\de_x}[Y^{\ga^\ast,h_{0,1}}_k\neq 0
\ \forall k\geq 0]\quad&(0<\ga^\ast<\infty),\\[5pt]
\dis P^{\de_x}[Y^{0,h_{0,1}}_t\neq 0\ \forall t\geq 0]\quad
&(\ga^\ast=0).\ea\right.
\ee
\el
We now explain how Lemmas~\ref{11lem}--\ref{01lem} imply
Proposition~\ref{Uit}. In doing so, it will be more convenient to work
with weighted branching processes than with Poissonized branching
processes. A little argument (which can be found in Lemma~\ref{exgro}
below) shows that Lemmas~\ref{11lem}--\ref{01lem} are equivalent to
the next proposition.
\bp{\bf(Extinction versus unbounded growth)}\label{P:exgr}
Let $h_{1,1}$, $h_{0,0}$, and $h_{0,1}$ be as in
Lemmas~\ref{11lem}--\ref{01lem}. For $\ga^\ast\in\half$,
put $p^\ast_{1,1,\ga^\ast}(x):=1$, $p^\ast_{0,0,\ga^\ast}(x):=0$
$(x\in[0,1])$, and
\be\label{pdef}
p^\ast_{0,1,\ga^\ast}(0):=0\quad\mbox{and}
\quad p^\ast_{0,1,\ga^\ast}(x):=h_{0,1}(x)\rho_{\ga^\ast}(x)\qquad(x\in(0,1]),
\ee
with $\rho_{\ga^\ast}$ as in (\ref{rhodef}). Then, for $(l,r)=(1,1),(0,0)$,
and $(0,1)$,
\be\label{intexp}
P^{-n,\de_x}\big[\li\Xc_0,h_{l,r}\re\in\cdot\,\big]\Asto{n}
\ex{-p^\ast_{l,r,\ga^\ast}(x)}\de_0
+\big(1-\ex{-p^\ast_{l,r,\ga^\ast}(x)}\big)\de_\infty,
\ee
uniformly for all $x\in[0,1]$.
\ep
Formula (\ref{intexp}) says that the weighted branching process
$\Xc^{h_{l,r}}$ exhibits a form of `extinction versus unbounded
growth'. More precisely, for large $n$ the total mass of
$h_{l,r}\Xc_0$ is close to $0$ or $\infty$ with high probability.\vc

\noi
{\bf Proof of Proposition~\ref{Uit}} By (\ref{Vi}),
\be\label{Vi3}
\Ui^{(n)}p(x)=-\log E^{-n,\de_x}\big[\ex{-\li\Xc_0,p\re}\big]
\qquad(p\in B_+[0,1],\ x\in[0,1]).
\ee
We first prove formula (\ref{Ucon})~(ii). For $(l,r)=(0,0)$,
formula (\ref{intexp}) says that
\be\label{intext}
P^{-n,\de_x}[\li\Xc_0,h_{0,0}\re\in\cdot\,]\Asto{n}\de_0
\ee
uniformly for all $x\in[0,1]$. If $p\in\Hi_{0,0}$, then we can
find $r>0$ such that $p\leq rh_{0,0}$. Therefore, (\ref{intext})
implies that for any $p\in\Hi_{0,0}$,
\be\label{intext2}
P^{-n,\de_x}[\li\Xc_0,p\re\in\cdot\,]\Asto{n}\de_0.
\ee
By (\ref{Vi3}) it follows that
\be\label{innie}
\Ui^{(n)}p(x)=-\log E^{-n,\de_x}\big[\ex{-\li\Xc_0,p\re}\big]\asto{n}0,
\ee
where the limits in (\ref{intext2}) and (\ref{innie}) are uniform in
$x\in[0,1]$. This proves formula (\ref{Ucon})~(ii). To prove formula
(\ref{Ucon})~(iii), note that for any $p\in\Hi_{0,1}$ we can choose
$0<r_-<r_+$ such that $r_-h_{0,1}\leq p+h_{0,0}\leq r_+h_{0,1}$.
Therefore, (\ref{intexp}) implies that
\be\label{intext3}
P^{-n,\de_x}[\li\Xc_0,p\re+\li\Xc_0,h_{0,0}\re\in\cdot\,]\Asto{n}
\ex{-p^\ast_{0,1,\ga^\ast}(x)}\de_0
+\big(1-\ex{-p^\ast_{0,1,\ga^\ast}(x)}\big)\de_\infty.
\ee
Using moreover (\ref{intext}), we see that
\be\label{intext4}
P^{-n,\de_x}[\li\Xc_0,p\re\in\cdot\,]\Asto{n}
\ex{-p^\ast_{0,1,\ga^\ast}(x)}\de_0
+\big(1-\ex{-p^\ast_{0,1,\ga^\ast}(x)}\big)\de_\infty.
\ee
By (\ref{Vi3}), it follows that
\be
\Ui^{(n)}p(x)=-\log E^{-n,\de_x}\big[\ex{-\li\Xc_0,p\re}\big]\asto{n}
p^\ast_{0,1,\ga^\ast}(x)
\ee
where all limits are uniform in $x\in[0,1]$. This proves
(\ref{Ucon})~(iii). The proof of (\ref{Ucon})~(i) is similar but easier.\qed

\section{Discussion, open problems}\label{dispro}

\subsection{Discussion}\label{discus}

Consider a $([0,1]^2)^{\Z^2}$-valued process
$\x=(\x_\xi)_{\xi\in\Z^2}=(\x^1_\xi,\x^2_\xi)_{\xi\in\Z^2}$,
solving a system of SDE's of the form
\bc\label{Z2}
\di\x^1_\xi(t)&=&\dis\sum_{\eta:\,|\eta-\xi|=1}\!
\big(\x^1_\eta(t)-\x^1_\xi(t)\big)\,\di t
+\sqrt{2\al\x^1_\xi(t)(1-\x^1_\xi(t))}\,\di B^1_\xi(t),\\[5pt]
\di\x^2_\xi(t)&=&\dis\sum_{\eta:\,|\eta-\xi|=1}\!
\big(\x^2_\eta(t)-\x^2_\xi(t)\big)\,\di t
+\sqrt{2p(\x^1_\xi(t))\x^2_\xi(t)(1-\x^2_\xi(t))}\,\di B^2_\xi(t),
\ec
where $\al>0$ is a constant, $p$ is a nonnegative function on $[0,1]$
satisfying $p(0)=0$ and $p(1)>0$, and $(B^i_\xi)^{i=1,2}_{\xi\in\Z^2}$
is a collection of independent Brownian motions. We call $\x$ a system
of linearly interacting catalytic Wright-Fisher diffusions with
catalyzation function $p$. It is expected that $\x$ {\em clusters},
i.e., $\x(t)$ converges in distribution as $t\to\infty$ to a limit
$(\x_\xi(\infty))_{\xi\in\Z^2}$ such that
$\x_\xi(\infty)=\x_0(\infty)$ for all $\xi\in\Z^2$ and $\x_0(\infty)$
takes values in the effective boundary associated with the diffusion
matrix $w^{\al,p}$ (see (\ref{eff})). Heuristic arguments, based on
renormalization, yield a formula for the {\em clustering distribution}
$\Li(\x_0(\infty))$ in terms of the diffusion matrix $w^\ast$ which is
the unique solution of the asymptotic fixed point equation
(\ref{afp})~(ii) in the renormalization class $\Wi^{0,1}_{\rm cat}$;
see Conjecture~\ref{Z2con} in Appendix~\ref{clustapp} below.

The present paper is inspired by the work of Greven, Klenke and
Wakolbinger \cite{GKW01}. They study a model that is closely related
to (\ref{Z2}), but where $\x^1$ is replaced by a voter model. They
show that their model clusters and determine its clustering
distribution $\Li(\x_0(\infty))$, which turns out to coincide with the
mentioned prediction for (\ref{Z2}) based on renormalization theory.
In fact, they believe their results to hold for the model in
(\ref{Z2}) too, but they could not prove this due to certain technical
difficulties that a $[0,1]$-valued catalyst would create, compared to
the simpler $\{0,1\}$-valued voter model.

The work in \cite{GKW01} not only provides the main motivation for the
present paper, but also inspired some of our techniques for proving
Theorem~\ref{main}. This concerns in particular the proof of
Proposition~\ref{identprop}, which makes the connection between
renormalization transformations and a branching process. We hope that
conversely, our techniques may shed some light on the problems left
open by \cite{GKW01}, in particular, the question whether their
results stay true if the voter model catalyst is replaced by a
Wright-Fisher catalyst. It seems plausible that their results may not
hold for the model in (\ref{Z2}) if the catalyzing function $p$ grows
too fast at $0$. On the other hand, our proofs suggest that $p$ with a
finite slope at $0$ should be OK. (In particular, while deriving
formula (\ref{intext3}), we use that $p$ can be bounded from above by
$r_+h_{0,1}$ for some $r_+>0$, which requires that $p$ has a finite
slope at $0$.)

Our results are also interesting in the wider program of studying
renormalization classes in the sense of Definition~\ref{defi}. We
conjecture that the class $\Wi^{0,1}_{\rm cat}$, unlike all
renormalization classes studied previously, contains no fixed shapes
(see the discussion following Lemma~\ref{L:shape}). In fact, we expect
this to be the usual situation. In this sense, the renormalization
classes studied so far were all of a special type.

\subsection{Open problems}\label{probsec}

The general program of studying renormalization classes in the sense
of Definition~\ref{defi} contains a wealth of open problems. In our
proofs, we make heavy use of the single-way nature of the catalyzation
in (\ref{catsde}), in particular, the fact that $\y^1$ is an
autonomous process which allows one to condition on $\y^1$ and
consider $\y^2$ as a process in a random environment created by
$\y^1$.  As soon as one leaves the single-way catalytic regime one
runs into several difficulties, both technically (it is hard to prove
that a given class of matrices is a renormalization class in the sense
of Definition~\ref{defi}) and conceptually (it is not clear when
solutions to the asymptotic fixed shape equation (\ref{afp})~(ii) are
unique). Therefore, it seems at present hard to verify the complete
picture for renormalization classes on the unit square that arises
from the numerical simulations described in Section~\ref{numer} and
Figures~\ref{fixfig} and~\ref{graph}, unless one or more essential new
ideas are added.

In this context, the study of the nonlinear partial differential
equation (\ref{flow}) and its fixed points seems to be a challenging
problem. This may be a hard problem from an analytic point of view,
since the equation is degenerate and not in divergence form. For the
renormalization class $\Wi_{\rm cat}$, the quasilinear equation
(\ref{flow}) reduces to the semilinear equation (\ref{cau}), which is
analytically easier to treat and moreover has a probabilistic
interpretation in terms of a superprocess. For a study of the
semilinear equation (\ref{cau}) we refer to \cite{FSsup}. We do not
know whether solutions to equation (\ref{flow}) can in general be
represented in terms of a stochastic process of some sort.

Even for the renormalization class $\Wi_{\rm cat}$, several
interesting problems are left open. One of the most urgent ones is to
prove that the functions $p^\ast_{0,1,\ga^\ast}$ are not constant in
$\ga^\ast$, and therefore, by Lemma~\ref{L:shape}~(c), $\Wi^{0,1}_{\rm
  cat}$ contains no fixed shapes. Moreover, we have not investigated
the iterated renormalization transformations in the regime
$\ga^\ast=\infty$.  Also, we believe that the convergence in
(\ref{Ucon})~(ii) does not hold if the condition that $p$ is Lipschitz
is dropped, in particular, if $p$ has an infinite slope at $0$ or an
infinite negative slope at $1$. For $p\in\Hi_{0,0}$, it seems
plausible that a properly rescaled version of the iterates
$\Ui^{(n)}p$ converges to a universal limit, but we have not
investigated this either. Finally, we have not investigated the
convergence of the iterated kernels $K^{w,(n)}$ from (\ref{Kdef}) (in
particular, we have not verified Conjecture~\ref{Kcon}) for the
renormalization class $\Wi_{\rm cat}$.

Our methods, combined with those in \cite{BCGH95}, can probably be
extended to study the action of iterated renormalization
transformations on diffusion matrices of the following more general
form (compared to (\ref{walp})):
\be\label{walp2}
w(x)=\left(\ba{@{}cc@{}}g(x_1)&0\\0&p(x_1)x_2(1-x_2)\ea\right)
\qquad(x=\in[0,1]^2),
\ee
where $g:[0,1]\to\R$ is Lipschitz, $g(0)=g(1)=0$, $g>0$ on $(0,1)$,
and $p\in\Hi$ as before. This would, however, require a lot of extra
technical work and probably not generate much new insight. The
numerical simulations mentioned in Section~\ref{numer} suggest that
many diffusion matrices of an even more general form than
(\ref{walp2}) also converge under renormalization to the limit points
$w^\ast$ from Theorem~\ref{main}, but we don't know how to prove this.

\vspace{20pt}
\noi
{\bf\Large Part II}
\vspace{5pt}

\noi
{\bf Outline of Part~II} In Section~\ref{Wcatsec}, we verify that
$\Wi_{\rm cat}$ is a renormalization class, we prove
Proposition~\ref{identprop}, which connects the renormalization
transformations $F_c$ to the log-Laplace operators $\Ui_\ga$, and we
collect a number of technical properties of the operators $\Ui_\ga$
that will be needed later on. In Section~\ref{convsec} we prove
Theorem~\ref{supercon} about the convergence of the renormalization
branching process to a time-homogeneous limit. In
Section~\ref{embsec}, we prove the statements from Section~\ref{exex}
about extinction versus unbounded growth of embedded particle systems,
with the exception of Lemma~\ref{00lem}, which is proved in
Section~\ref{00sec}. In Section~\ref{final}, finally, we combine the
results derived by that point to prove our main theorem.

\section{The renormalization class $\Wi_{\rm cat}$}\label{Wcatsec}

In this section we prove Theorem~\ref{main}~(a) and
Proposition~\ref{identprop}, as well as
Lemmas~\ref{contlem}--\ref{L:shape} from Section~\ref{gensec}. The
section is organized according to the techniques used.
Section~\ref{gensub} collects some facts that hold for general
renormalization classes on compact sets. In Section~\ref{coupsec} we
use the SDE (\ref{catsde}) to couple catalytic Wright-Fisher
diffusions. In Section~\ref{dualsec} we apply the moment duality for
the Wright-Fisher diffusion to the catalyst and to the reactant
conditioned on the catalyst. In Section~\ref{concavesec} we prove that
monotone concave catalyzing functions form a preserved class under
renormalization.

\subsection{Renormalization classes on compact sets}\label{gensub}

In this section, we prove the lemmas stated in Section~\ref{gensec}.
Recall that $D\sub\R^d$ is open, bounded, and convex, and that $\Wi$
is a prerenormalization class on $\ov D$, equipped with the topology
of uniform convergence.\med

\noi
{\bf Proof of Lemma~\ref{contlem}} To see that
$(x,c,w)\mapsto\nu^{c,w}_x$ is continuous, let $(x_n,c_n,w_n)$ be a
sequence converging in $\ov D\times(0,\infty)\times\Wi$ to a limit
$(x,c,w)$. By the compactness of $\ov D$, the sequence
$(\nu_{x_n}^{c_n,w_n})_{n\geq 0}$ is tight, and each limit point
$\nu^\ast$ satisfies
\be\label{invar}
\li\nu^\ast,A^{c,w}_xf\re=0\qquad(f\in\Ci^{(2)}(D)).
\ee
Therefore, by \cite[Theorem~4.9.17]{EK}, $\nu^\ast$ is an invariant
law for the martingale problem associated with $A^{c,w}_x$. Since we
are assuming uniqueness of the invariant law, $\nu^\ast=\nu^{c,w}_x$
and therefore $\nu^{c_n,w_n}_{x_n}\Rightarrow\nu^{c,w}_x$. The
continuity of $F_cw(x)$ is a simple consequence of the continuity of
$\nu^{c,w}_x$.\qed

\noi
{\bf Proof of Lemma~\ref{schaal}} Formula (\ref{schafo})~(i) follows
from the fact that rescaling the time in solutions $(\y_t)_{t\geq 0}$
to the martingale problem for $A^{c,w}_x$ by a factor $\la$ has no
influence on the invariant law. Formula (\ref{schafo})~(ii) is a
direct consequence of formula (\ref{schafo})~(i).\qed

\noi
{\bf Proof of Lemma~\ref{numom}} This follows by inserting the
functions $f(x)=x_i$ and $f(x)=x_ix_j$ into the equilibrium equation
(\ref{invar}).\qed

\noi
{\bf Proof of Lemma~\ref{effinv}} If $x\in\pa_wD$, then $\y_t:=x$
($t\geq 0$) is a stationary solution to the martingale problem for
$A^{c,w}_x$, and therefore $\nu^{c,w}_x=\de_x$ and $F_cw(x)=w(x)=0$.
On the other hand, if $x\not\in\pa_wD$, then $\y_t:=x$ ($t\geq 0$) is
not a stationary solution to the martingale problem for $A^{c,w}_x$
and therefore $\int_{\ov D}\nu^{c,w}_x(\di y)|y-x|^2>0$. Let
$\tr(w(y)):=\sum_iw_{ii}(y)$ denote the trace of $w(y)$. By
(\ref{moments2})~(ii), $\frac{1}{c}\tr(F_cw)(x)=\frac{1}{c}\int_{\ov
  D}\nu^{c,w}_x(\di y)\tr(w(y))=\int_{\ov D}\nu^{c,w}_x(\di
y)|y-x|^2>0$ and therefore $F_cw(x)\neq 0$.\qed

\noi
{F}rom now on assume that $\Wi$ is a renormalization class. Note that
\be\label{iterK}
K^{w,(n)}=\nu^{c_{n-1},F^{(n-1)}w}\cdots\nu^{c_0,w}\qquad(n\geq 1),
\ee
where we denote the composition of two probability kernels $K,L$ on $\ov D$ by
\be
(K L)_x(\di z):=\int_{\ov D} K_x(\di y)L_y(\di z).
\ee
{\bf Proof of Lemma~\ref{basic}} This is a direct consequence of
Lemmas~\ref{contlem} and \ref{numom}. In particular, the relations
(\ref{moments}) follow by iterating the relations (\ref{moments2}).\qed

\noi
{\bf Proof of Lemma~\ref{clusK}} Recall that $\tr(w(y))$ denotes the
trace of $w(y)$. Formulas (\ref{KF}) and (\ref{moments})~(ii) show that
\be
\int_{\ov D}\! K^{w,(n)}_x(\di y)\,|y-x|^2
=s_n\!\int_{\ov D}\! K^{w,(n)}_x(\di y)\,\tr(w(y)).
\ee
Since $\ov D$ is compact, the left-hand side of this equation is
bounded uniformly in $x\in\ov D$ and $n\geq 1$, and therefore,
since we are assuming $s_n\to\infty$,
\be\label{tranul}
\lim_{n\to\infty}\sup_{x\in D}\int_{\ov D} K^{w,(n)}_x(\di y)\tr(w(y))=0.
\ee
Since $w$ is symmetric and nonnegative definite, $\tr(w(y))$ is
nonnegative, and zero if and only if $y\in\pa_w D$. If $f\in\Ci(\ov
D)$ satisfies $f=0$ on $\pa_w D$, then, for every $\eps>0$, the sets
$C_m:=\{x\in\ov D:|f(x)|\geq\eps+m\,\tr(w(x))\}$ are compact with
$C_m\down\emptyset$ as $m\up\infty$, so there exists an $m$ (depending
on $\eps$) such that $|f|<\eps+m\,\tr(w)$. Therefore,
\be\ba{l}
\dis\limsup_{n\to\infty}\:\sup_{x\in\ov D}
\Big|\int_{\ov D}K^{w,(n)}_x(\di y)f(y)\Big|
\leq\limsup_{n\to\infty}\:\sup_{x\in\ov D}
\int_{\ov D}K^{w,(n)}_x(\di y)|f(y)|\\
\dis\quad\leq\eps+m\limsup_{n\to\infty}\:\sup_{x\in\ov D}
\int_{\ov D}K^{w,(n)}_x(\di y)\tr(w(y))=\eps.
\ec
Since $\eps>0$ is arbitrary, (\ref{concen}) follows.\qed

\noi
{\bf Proof of Lemma~\ref{L:shape}} By (\ref{sFw}), (\ref{sumga}), and
(\ref{gamma}), $w^\ast_{\ga^\ast}=\lim_{n\to\infty}(\ov
F_{\ga^\ast})^nw$ for each $w\in\Wi$. By Lemma~\ref{contlem}~(b), $\ov
F_{\ga^\ast}:\Wi\to\Wi$ is continuous, so $w^\ast_{\ga^\ast}$ is the
unique fixed point of $\ov F_{\ga^\ast}$. This proves part~(a).

Now let $0\neq w\in\Wi$ and assume that $\hat\Wi=\{\la w:\la>0\}$ is a
fixed shape. Then $\hat\Wi\ni s_nF^{(n)}w\asto{n}w^\ast_{\ga^\ast}$
whenever $s_n\to\infty$ and $s_{n+1}/s_n\to1+\ga^\ast$ for some
$0<\ga^\ast<\infty$, which shows that $\hat\Wi=\{\la
w^\ast_{\ga^\ast}:\la>0\}$. Thus, $\Wi$ can contain at most one fixed
shape, and if it does, then the $w^\ast_{\ga^\ast}$ for different
values of $\ga^\ast$ must be constant multiples of each other. This
proves part~(c) and the uniqueness statement in part~(b).

To complete the proof of part~(b), note that if $w^\ast=w^\ast_{\ga^\ast}$ does not depend on $\ga^\ast$, then $w^\ast\in\Wi$ solves (\ref{afp})~(i) for all $0<\ga^\ast<\infty$, hence $F_cw^\ast=(1+\frac{1}{c})^{-1}w^\ast$ for all $c>0$, and therefore, by scaling (Lemma~\ref{schaal}), $F_c(\la w^\ast)=\la F_{c/\la}(w^\ast)=\la(1+\frac{\la}{c})^{-1}w^\ast=(\frac{1}{\la}+\frac{1}{c})^{-1}w^\ast$.\qed

\subsection{Coupling of catalytic Wright-Fisher diffusions}\label{coupsec}

In this section we verify condition (i) of Definition~\ref{defi} for
the class $\Wi_{\rm cat}$, and we prepare for the verification of
conditions (ii)--(iv) in Section~\ref{dualsec}. In fact, we will show
that the larger class $\ov\Wi_{\rm cat}:=\{w^{\al,p}:\al>0,\
p\in\Ci_+[0,1]\}$ is also a renormalization class, and the equivalents
of Theorem~\ref{main}~(a) and Proposition~\ref{identprop} remain true
for this larger class. (We do not know, however, if the convergence
statements in Theorem~\ref{main}~(b) also hold in this larger class;
see the discussion in Section~\ref{probsec}.)

For each $c\geq 0$, $w\in\ov\Wi_{\rm cat}$ and $x\in[0,1]^2$, the
operator $A^{c,w}_x$ is a densely defined linear operator on
$\Ci([0,1]^2)$ that maps the identity function into zero and, as one
easily verifies, satisfies the positive maximum principle. Since
$[0,1]^2$ is compact, the existence of a solution to the martingale
problem for $A^{c,w}_x$, for each $[0,1]^2$-valued initial condition,
now follows from general theory (see \cite{RW87}, Theorem~5.23.5, or
\cite[Theorem~4.5.4 and Remark~4.5.5]{EK}).

We are therefore left with the task of verifying uniqueness of
solutions to the martingale problem for $A^{c,w}_x$. By
\cite[Problem~4.19, Corollary~5.3.4, and Theorem 5.3.6]{EK}, it
suffices to show that solutions to (\ref{catsde}) are pathwise unique.
\bl[Monotone coupling of Wright-Fisher diffusions]\label{leq}
Assume that $0\leq x\leq \ti x\leq 1$, $c\geq 0$ and that
$(P_t)_{t\geq 0}$ is a progressively measurable, nonnegative
process such that $\sup_{t\geq 0,\oo\in\om}P_t(\oo)<\infty$.
Let $\y,\ti\y$ be $[0,1]$-valued solutions to the SDE's
\bc
\di\y_t&=&\dis c\,(x-\y_t)\di t+\sqrt{2P_t\y_t(1-\y_t)}\di B_t,\\
\di\ti\y_t&=&\dis c\,(\ti x-\ti\y_t)\di t+\sqrt{2P_t\ti\y_t(1-\ti\y_t)}\di B_t,
\ec
where in both equations $B$ is the same Brownian motion. If
$\y_0\leq \ti\y_0$ a.s., then
\be
\y_t\leq \ti\y_t\quad\forall t\geq 0\quad\as
\ee
\el
{\bf Proof} This is an easy adaptation of a technique due to
Yamada and Watanabe \cite{YW71}. Since $\int_{0+}\frac{\di x}{x}=\infty$,
it is possible to choose $\rho_n\in\Ci\half$ such that
$\int_0^\infty\rho_n(x)\di x=1$ and
\be
0\leq\rho_n(x)\leq\frac{1}{nx}1_{(0,1]}(x)\qquad\qquad(x\geq 0).
\ee
Define $\phi_n\in\Ci^{(2)}(\R)$ by
\be
\phi_n(x):=\int_0^{x\vee 0}\!\!\di y\int_0^y\!\!\di z\,\rho_n(z).
\ee
One easily verifies that $\phi_n(x)$, $x\phi'_n(x)$, and $x\phi''_n(x)$
are nonnegative and converge, as $n\to\infty$, to $x\vee 0$, $x\vee 0$,
and $0$, respectively. By It\^o's formula:
\be\ba{r@{\,}c@{\,}l@{\qquad}r}
\dis E[\phi_n(\y_t-\ti\y_t)]&=&\dis E[\phi_n(\y_0-\ti\y_0)]&{\rm(i)}\\
&&\dis+c\,(x-\ti x)\int_0^t E[\phi'_n(\y_s-\ti\y_s)]\di s
-c\int_0^tE[(\y_s-\ti\y_s)\phi'_n(\y_s-\ti\y_s)]\di s&{\rm(ii)}\\
&&\dis+\int_0^tE\Big[P_s\Big(\sqrt{\y_s(1-\y_s)}
-\sqrt{\ti\y_s(1-\ti\y_s)}\Big)^2\phi''_n(\y_s-\ti\y_s)\Big]\di s.&{\rm(iii)}
\ec
Here the terms in (ii) are nonpositive, and hence, letting
$n\to\infty$ and using the elementary estimate
\be\label{Hoelder}
|\sqrt{y(1-y)}-\sqrt{\ti y(1-\ti y)}|\leq |y-\ti y|^{\frac{1}{2}}
\qquad(y,\ti y\in[0,1]),
\ee
the properties of $\phi_n$, and the fact that the process $P$ is
uniformly bounded, we find that
\be
E[0\vee(\y_t-\ti\y_t)]\leq E[0\vee(\y_0-\ti\y_0)]=0,
\ee
by our assumption that $\y_0\leq\ti\y_0$. This shows that
$\y_t\leq\ti\y_t$ a.s. for each fixed $t\geq 0$, and by the
continuity of sample paths the statement holds for all
$t\geq 0$ almost surely.\qed

\begin{corollary}[Pathwise uniqueness]\label{uni}
For all $c\geq 0$, $\al>0$, $p\in\Ci_+[0,1]$ and $x\in[0,1]$,
solutions to the SDE (\ref{catsde}) are pathwise unique.
\end{corollary}
{\bf Proof} Let $(\y^1,\y^2)$ and $(\ti\y^1,\ti\y^2)$ be solutions to
(\ref{catsde}) relative to the same pair $(B^1,B^2)$ of Brownian
motions, with $(\y^1_0,\y^2_0)=(\ti\y^1_0,\ti\y^2_0)$. Applying
Lemma~\ref{leq}, with inequality in both directions, we see that
$\y^1=\ti\y^1$ a.s. Applying Lemma~\ref{leq} two more times, this time
using that $\y^1=\ti\y^1$ a.s., we see that also $\y^2=\ti\y^2$
a.s.\qed

\begin{corollary}[Exponential coupling]\label{excoup}
Assume that $x\in[0,1]$, $c\geq 0$, and $\al>0$. Let $\y,\ti\y$
be solutions to the SDE
\be\label{alWF}
\di\y_t=c\,(x-\y_t)\di t+\sqrt{2\al\y_t(1-\y_t)}\di B_t,
\ee
relative to the same Brownian motion $B$. Then
\be\label{exco}
E\big[|\ti\y_t-\y_t|\big]=e^{-ct}E\big[|\ti\y_0-\y_0|\big].
\ee
\end{corollary}
{\bf Proof} If $\y_0=y$ and $\ti\y_0=\ti y$ are deterministic and
$y\leq\ti y$, then by Lemma~\ref{leq} and a simple moment calculation
\be\label{absmom}
E\big[|\ti\y_t-\y_t|\big]=E[\ti\y_t-\y_t]=e^{-ct}|\ti y-y|.
\ee
The same argument applies when $y\geq\ti y$. The general case where
$\y_0$ and $\ti\y_0$ are random follows by conditioning on
$(\y_0,\ti\y_0)$.\qed
\begin{corollary}[Ergodicity]\label{ergo}
The Markov process defined by the SDE (\ref{Yclx}) has a unique
invariant law $\Ga^\ga_x$ and is ergodic, i.e, solutions to
(\ref{Yclx}) started in an arbitrary initial law $\Li(\y_0)$
satisfy $\Li(\y_t)\Asto{t}\Ga^\ga_x$.
\end{corollary}
{\bf Proof} Since our process is a Feller diffusion on a compactum,
the existence of an invariant law follows from a simple time averaging
argument. Now start one solution $\ti\y$ of (\ref{Yclx}) in this
invariant law and let $\y$ be any other solution, relative to the same
Brownian motion. Corollary~\ref{excoup} then gives ergodicity and, in
particular, uniqueness of the invariant law.\qed

\brm{\bf (Density of invariant law)}\label{R:betadis}
It is well-known (see, for example \cite[formula~(5.70)]{Ewe04}) that
$\Ga^\ga_x$ is a $\bet(\al_1,\al_2)$-distribution, where
$\al_1:=x/\ga$ and $\al_2:=(1-x)/\ga$, i.e., $\Ga^\ga_x=\de_x$
$(x\in\{0,1\})$ and
\be\label{Gadef}
\Ga^\ga_x(\di y)=\frac{\Ga(\al_1+\al_2)}{\Ga(\al_1)\Ga(\al_2)}
\,y^{\al_1-1}(1-y)^{\al_2-1}\di y\qquad(x\in(0,1)).
\ee
\erm

\noi
We conclude this section with a lemma that prepares for the verification
of condition (iv) in Definition~\ref{defi} for the class $\Wi_{\rm cat}$.
\bl\label{monot}
{\bf (Monotone coupling of stationary Wright-Fisher diffusions)} Assume
that $c>0$, $\al>0$ and $0\leq x\leq\ti x\leq 1$. Then the pair of equations
\bc\label{kopwf}
\di\y_t&=&\dis c\,(x-\y_t)\di t+\sqrt{2\al\y_t(1-\y_t)}\di B_t,\\
\di\ti\y_t&=&\dis c\,(\ti x-\ti\y_t)\di t+\sqrt{2\al\ti\y_t(1-\ti\y_t)}\di B_t
\ec
has a unique stationary solution $(\y_t,\ti\y_t)_{t\in\R}$. This
stationary solution satisfies
\be\label{leqt}
\y_t\leq\ti\y_t\quad\forall t\in\R\quad\as
\ee
\el
{\bf Proof} Let $(\y_t,\ti\y_t)_{t\geq 0}$ be a solution of
(\ref{kopwf}) and let $(\y'_t,\ti\y'_t)_{t\geq 0}$ be another one,
relative to the same Brownian motion $B$. Then, by Lemma~\ref{excoup},
$E[|\y_t-\y'_t|]\to 0$ and also $E[|\ti\y_t-\ti\y'_t|]\to 0$ as
$t\to\infty$. Hence we may argue as in the proof of
Corollary~\ref{ergo} that (\ref{kopwf}) has a unique invariant law and
is ergodic. Now start a solution of (\ref{kopwf}) in an initial
condition such that $\y_0\leq\ti\y_0$. By ergodicity, the law of this
solution converges as $t\to\infty$ to the invariant law of
(\ref{kopwf}) and using Lemma~\ref{leq} we see that this invariant law
is concentrated on $\{(y,\ti y)\in[0,1]^2:y\leq\ti y\}$. Now consider,
on the whole real time axis, the stationary solution to (\ref{kopwf})
with this invariant law. Applying Lemma~\ref{leq} once more, we see
that (\ref{leqt}) holds.\qed

\subsection{Duality for catalytic Wright-Fisher diffusions}\label{dualsec}

In this section we prove Theorem~\ref{main}~(a) and
Proposition~\ref{identprop}. Moreover, we will show that their
statements remain true if the renormalization class $\Wi_{\rm cat}$ is
replaced by the larger class $\ov\Wi_{\rm cat}:=\{w^{\al,p}:\al>0,\
p\in\Ci_+[0,1]\}$. We begin by recalling the usual moment duality for
Wright-Fisher diffusions.

For $\ga>0$ and $x\in[0,1]$, let $\y$ be a solution to the SDE
\be\label{WFsde}
\di\y(t)=\ffrac{1}{\ga}\,(x-\y(t))\di t+\sqrt{2\y(t)(1-\y(t))}\di B(t),
\ee
i.e., $\y$ is a Wright-Fisher diffusion with a linear drift towards $x$.
It is well-known that $\y$ has a moment dual. To be precise, let
$(\phi,\psi)$ be a Markov process in $\N^2=\{0,1,\ldots\}^2$ that jumps as:
\be\ba{r@{\,}c@{\,}l@{\qquad}l}\label{phipsi}
(\phi_t,\psi_t)&\to&(\phi_t-1,\psi_t)
&\mbox{with rate }\ \phi_t(\phi_t-1)\\
(\phi_t,\psi_t)&\to&(\phi_t-1,\psi_t+1)
&\mbox{with rate }\ \ffrac{1}{\ga}\phi_t.
\ec
Then one has the following {\em duality relation} (see for example
Lemma 2.3 in \cite{Shi80} or Proposition~1.5 in \cite{GKW01})
\be\label{WFdual}
E^y\big[\y_t^nx^m\big]=E^{(n,m)}\big[y^{\phi_t}x^{\psi_t}\big]
\qquad(y\in[0,1],\ (n,m)\in\N^2),
\ee
where $0^0:=1$. The duality in (\ref{WFdual}) has the following
heuristic explanation. Consider a population containing a fixed, large
number of organisms, that come in two genetic types, say I and II.
Each pair of organisms in the population is {\em resampled} with rate
$2$. This means that one organism of the pair (chosen at random) dies,
while the other organism produces one child of its own genetic type.
Moreover, each organism is replaced with rate $\frac{1}{\ga}$ by an
organism chosen from an infinite reservoir where the frequency of type
I has the fixed value $x$. In the limit that the number of organisms
in the population is large, the relative frequency $\y_t$ of type I
organisms follows the SDE (\ref{WFsde}). Now $E[\y_t^n]$ is the
probability that $n$ organisms sampled from the population at time $t$
are all of type I. In order to find this probability, we follow the
ancestors of these organisms back in time. Viewed backwards in time,
these ancestors live for a while in the population, until, with rate
$\frac{1}{\ga}$, they jump to the infinite reservoir. Moreover, due to
resampling, each pair of ancestors coalesces with rate $2$ to one
common ancestor. Denoting the number of ancestors that lived at time
$t-s$ in the population and in the reservoir by $\phi_s$ and $\psi_s$,
respectively, we see that the probability that all ancestors are of
type~I is $E^y[\y_t^n]=E^{(n,0)}[y^{\phi_t}x^{\psi_t}]$. This gives a
heuristic explanation of (\ref{WFdual}).

Since eventually all ancestors of the process $(\phi,\psi)$ end up in
the reservoir, we have $(\phi_t,\psi_t)\to(0,\psi_\infty)$ as
$t\to\infty$ a.s.\ for some $\N$-valued random variable $\psi_\infty$.
Taking the limit $t\to\infty$ in (\ref{WFdual}), we see that the
moments of the invariant law $\Ga^\ga_x$ from Corollary~\ref{ergo} are
given by:
\be\label{Gamom}
\int\Ga^\ga_x(\di y)y^n=E^{(n,0)}[x^{\psi_\infty}]\qquad(n\geq 0).
\ee
It is not hard to obtain an inductive formula for the moments of
$\Ga^\ga_x$, which can then be solved to yield the formula
\be\label{WFmoments}
\int\Ga^\ga_x(\di y)y^n=\prod_{k=0}^{n-1}\frac{x+k\ga}{1+k\ga}\qquad(n\geq 1).
\ee
In particular, it follows that
\be\label{WFfix}
\int\Ga^\ga_x(\di y)y(1-y)=\frac{1}{1+\ga}x(1-x).
\ee
This is the important {\em fixed shape property} of the Wright-Fisher
diffusion (see formula (\ref{fix})).

We now consider catalytic Wright-Fisher diffusions $(\y^1,\y^2)$ as in
(\ref{catsde}) with $p\in\Ci_+[0,1]$ and apply duality to the catalyst
$\y^2$ conditioned on the reactant $\y^1$. Let
$(\y^1_t,\y^2_t)_{t\in\R}$ be a stationary solution to the SDE
(\ref{catsde}) with $c=1/\ga$. Let $(\ti\phi,\ti\psi)$ be a
$\N^2$-valued process, defined on the same probability space as
$(\y^1,\y^2)$, such that conditioned on the past path
$(\y^1_{-t})_{t\leq 0}$, the process $(\ti\phi,\ti\psi)$ is a
(time-inhomogeneous) Markov process that jumps as:
\be\ba{r@{\,}c@{\,}l@{\qquad}l}\label{yphipsi}
(\ti\phi_t,\ti\psi_t)&\to&(\ti\phi_t-1,\ti\psi_t)
&\mbox{with rate }\ p(\y^1_{-t})\ti\phi_t(\ti\phi_t-1),\\
(\ti\phi_t,\ti\psi_t)&\to&(\ti\phi_t-1,\ti\psi_t+1)
&\mbox{with rate }\ \ffrac{1}{\ga}\ti\phi_t.
\ec
Then, in analogy with (\ref{WFdual}),
\be\label{reactdual}
E[(\y^2_0)^nx_2^m|(\y^1_{-t})_{t\leq 0}]
=E^{(n,m)}[(\y^2_{-t})^{\ti\phi_t}x_2^{\ti\psi_t}|(\y^1_{-t})_{t\leq 0}]
\qquad((n,m)\in\N^2,\ t\geq 0).
\ee
We may interpret (\ref{yphipsi}) by saying that pairs of ancestors in
a finite population coalesce with time-dependent rate $2p(\y^1_{-t})$
and ancestors jump to an infinite reservoir with constant rate
$\frac{1}{\ga}$.  Again, eventualy all ancestors end up in the
reservoir, and therefore $(\ti\phi_t,\ti\psi_t)\to(0,\ti\psi_\infty)$
as $t\to\infty$ a.s.\ for some $\N$-valued random variable
$\ti\psi_\infty$. Taking the limit $t\to\infty$ in (\ref{reactdual})
we find that
\be\label{reactdual2}
E[(\y^2_0)^nx_2^m|(\y^1_{-t})_{t\leq 0}]
=E^{(n,m)}[x_2^{\ti\psi_\infty}|(\y^1_{-t})_{t\leq 0}]
\qquad((n,m)\in\N^2,\ t\geq 0).
\ee
\bl{\bf(Uniqueness of invariant law)}\label{unin}
For each $c>0$, $w\in\ov\Wi_{\rm cat}$, and $x\in[0,1]^2$, there exists
a unique invariant law $\nu^{c,w}_x$ for the martingale problem for
$A^{c,w}_x$.
\el
{\bf Proof} Our process being a Feller diffusion on a compactum, the
existence of an invariant law follows from time averaging. We need to
show uniqueness. If $(\y^1,\y^2)=\y^1_t,\y^2_t)_{t\in\R}$ is a
stationary solution, then $\y^1$ is an autonomous process, and
$\Li(\y^1_0)=\Ga^{1/c}_x$, the unique invariant law from
Corollary~\ref{ergo}. Therefore, $\Li((\y^1_t)_{t\in\R})$ is
determined uniquely by the requirement that $(\y^1,\y^2)$ be
stationary. By (\ref{reactdual2}), the conditional distribution of
$\y^2_0$ given $(\y^1_t)_{t\leq 0}$ is determined uniquely, and
therefore the joint distribution of $\y^2_0$ and $(\y^1_t)_{t\leq 0}$
is determined uniquely. In particular,
$\Li(\y^1_0,\y^2_0)=\nu^{c,w}_x$ is determined uniquely.\qed

\brm{\bf (Reversibility)} 
It seems that the invariant law $\nu^{c,w}_x$ from Lemma~\ref{unin} is
reversible. In many cases (densities of) reversible invariant measures
can be obtained in closed form by solving the equations of detailed
balance. This is the case, for example, for the one-dimensional
Wright-Fisher diffusion. We have not attempted this for the catalytic
Wright-Fisher diffusion.
\erm
The next proposition implies Proposition~\ref{identprop} and prepares
for the proof of Theorem~\ref{main}~(a).
\bp{\bf(Extended renormalization class)}\label{exclass}
The set $\ov\Wi_{\rm cat}$ is a renormalization class on $[0,1]^2$, and
\be\label{hutfo}
\hut F_\ga w^{\txt 1,p}=w^{\txt 1,\Ui_\ga p}\qquad(p\in\Ci_+[0,1],\ \ga>0).
\ee
\ep
{\bf Proof} To see that $\ov\Wi_{\rm cat}$ is a renormalization class
we need to check conditions~(i)--(iv) from Definition~\ref{defi}. By
Lemma~\ref{uni}, the martingale problem for $A^{c,w}_x$ is well-posed
for all $c\geq 0$, $w\in\Wi_{\rm cat}$ and $x\in[0,1]^2$. By
Lemma~\ref{unin}, the corresponding Feller process on $[0,1]^2$ has a
unique invariant law $\nu^{c,w}_x$. This shows that conditions~(i) and
(ii) from Definition~\ref{defi} are satisfied. Note that by the
compactness of $[0,1]^2$, any continuous function on $[0,1]^2$ is
bounded, so condition~(iii) is automatically satisfied. Hence $\Wi$ is
a prerenormalization class. As a consequence, for any
$p\in\Ci_+[0,1]$, $\hut F_\ga w^{1,p}$ is well-defined by (\ref{Fc})
and (\ref{hatFc}). We will now first prove (\ref{hutfo}) and then show
that $\ov\Wi_{\rm cat}$ is a renormalization class.

Fix $\ga>0$, $p\in\Ci_+[0,1]$, and $x\in[0,1]^2$. Let
$(\y^1_t,\y^2_t)_{t\in\R}$ be a stationary solution to the SDE
(\ref{catsde}) with $\al=1$ and $c=1/\ga$. Then
\be
\hut F_\ga w^{1,p}_{ij}(x)=(1+\ga)E[w^{1,p}_{ij}(\y^1_0,\y^2_0)]
\qquad(i,j=1,2).
\ee
Since $w^{1,p}_{ij}=0$ if $i\neq j$, it is clear that $\hut F_\ga
w^{1,p}_{ij}(x)=0$ if $i\neq j$. Since $\Li(\y^1_0)=\Ga^\ga_x$ it
follows from (\ref{WFfix}) that $\hut F_\ga
w^{1,p}_{11}(x)=x_1(1-x_1)$. We are left with the task of showing that
\be
\hut F_\ga w^{1,p}_{22}(x)=\Ui_\ga p(x_1)x_2(1-x_2).
\ee
Here, by (\ref{moments2})~(ii),
\bc\label{momtruc}
\dis\hut F_\ga w^{1,p}_{22}(x)&=&\dis(1+\ga)E[p(\y^1_0)\y^2_0(1-\y^2_0)]\\[5pt]
&=&\dis(\ffrac{1}{\ga}+1)E[(\y^2_0-x_2)^2].
\ec
By (\ref{reactdual2}), using the fact that $E[\y^2_0]=x_2$ (which
follows from (\ref{reactdual}) or more elementary from (\ref{moments})~(i)),
we find that
\be\label{coal}
E[(\y^2_0-x_2)^2]=E[(\y^2_0)^2]-(x_2)^2
=E^{(2,0)}[x_2^{\ti\psi_\infty}]-(x_2)^2=P^{(2,0)}[\ti\psi_\infty=1]x_2(1-x_2)
\qquad(t\geq 0).
\ee
Note that $P^{(2,0)}[\ti\psi_\infty=1]$ is the probability that the
two ancestors coalesce before one of them leaves the population. The
probability of {\em noncoalescence} is given by
\be\label{noncoa}
P^{(2,0)}[\ti\psi_\infty=2]
=E\big[\ex{-\int_0^{\frac{1}{2}\tau_\ga}2p(y^1_{-t})\di t}\big],
\ee
where $\tau_\ga$ is an exponentially distributed random variable
with mean $\ga$. Combining this with (\ref{momtruc}) and (\ref{coal})
we find
that
\bc
\dis\hut F_\ga w^{1,p}_{22}(x)&=&\dis(\ffrac{1}{\ga}+1)
E\big[1-\ex{-\int_0^{\tau_\ga}p(y^1_{-t/2})\di t}\big]x_2(1-x_2)\\[5pt]
&=&\dis q_\ga E\big[1-\ex{-\li\Zi^\ga_x,p\re}\big]x_2(1-x_2)\\[5pt]
&=&\dis\Ui_\ga p(x_1)x_2(1-x_2),
\ec
where we have used the definition of $\Ui_\ga$.

We still have to show that $\ov\Wi_{\rm cat}$ satisfies condition~(iv)
from Definition~\ref{defi}. For any $\al>0$ and $p\in\Ci_+[0,1]$, by
scaling (Lemma~\ref{schaal}) and (\ref{hutfo}),
\be\label{wp}
F_cw^{\txt\al,p}=\al F_{\frac{c}{\al}}w^{\txt 1,\frac{p}{\al}}
=\al(1+\frac{\al}{c})^{-1}\ov F_{\frac{c}{\al}}w^{\txt 1,\frac{p}{\al}}
=w^{\txt(\frac{1}{\al}+\frac{1}{c})^{-1},
(\frac{1}{\al}+\frac{1}{c})^{-1}\Ui_{\frac{c}{\al}}(\frac{p}{\al})}.
\ee
By Lemma~\ref{contlem}, this diffusion matrix is continuous, which
implies that $\Ui_{\frac{c}{\al}}(\frac{p}{\al})$ is continuous.\qed

\noi
Our proof of Propostion~\ref{exclass} has a corollary.
\bcor{\bf(Continuity in parameters)}\label{Ugacont}
The map $(x,\ga)\mapsto\Qi_\ga(x,\cdot)$ from $[0,1]\times(0,\infty)$ to
$\Mi_1(\Mi[0,1])$ and the map $(x,\ga,p)\mapsto\Ui_\ga p(x)$ from
$[0,1]\times(0,\infty)\times\Ci_+[0,1]$ to $\R$ are continuous.
\ecor
{\bf Proof} By Lemma~\ref{contlem}, the diffusion matrix in (\ref{wp})
is continuous in $x,\ga$, and $p$, which implies the continuity of
$\Ui_\ga p(x)$. It follows that the map
$(x,\ga)\mapsto\int\Qi_\ga(x,\di\chi)\ex{-\li\chi,f\re}$ is continuous
for all $f\in\Ci_+[0,1]$, so by \cite[Theorem~4.2]{Kal76},
$(x,\ga)\mapsto\Qi_\ga(x,\cdot)$ is continuous.\qed

\noi
{\bf Proof of Theorem~\ref{main}~(a)}\label{aproof} We need to show
that $\Wi_{\rm cat}$ is a renormalization class and that $F_c$ maps
the subclasses $\Wi^{l,r}_{\rm cat}$ into themselves. It has already
been explained in Section~\ref{gensec} that the latter fact is a
consequence of Lemma~\ref{effinv}. Since in Proposition~\ref{exclass}
it has been shown that $\ov\Wi_{\rm cat}$ is a renormalization class,
we are left with the task to show that $F_c$ maps $\Wi_{\rm cat}$ into
itself. By (\ref{hutfo}) and scaling, it suffices to show that
$\Ui_\ga$ maps $\Hi$ into itself.

Fix $0\leq x\leq\ti x\leq 1$. By Lemma~\ref{monot}, we can couple the
processes $\y^\ga_x$ and $\y^\ga_{\ti x}$ from (\ref{Yclx}) such that
\be\label{leqt2}
\y^\ga_x(t)\leq\y^\ga_{\ti x}(t)\quad\forall t\leq 0\quad\as
\ee
Since the function $z\mapsto 1-e^{-z}$ on $\half$ is Lipschitz continuous
with Lipschitz constant $1$,
\be\ba{l}\label{Lipone}
\dis\big|\Ui_\ga p(\ti x)-\Ui_\ga p(x)\big|\\[5pt]
\dis\quad=\Big|(\ffrac{1}{\ga}+1)
E\big[1-\ex{-\int_0^{\tau_\ga}p(\y^\ga_{\ti x}(-t/2))\di t}\big]
-(\ffrac{1}{\ga}+1)
E\big[1-\ex{-\int_0^{\tau_\ga}p(\y^\ga_x(-t/2))\di t}\big]\Big|\\[5pt]
\dis\quad\leq (\ffrac{1}{\ga}+1)E\Big[\int_0^{\tau_\ga}\big|
p(\y^\ga_{\ti x}(-t/2))-p(\y^\ga_x(-t/2))\big|\di t\Big]\\[5pt]
\dis\quad\leq (\ffrac{1}{\ga}+1)LE\Big[\int_0^{\tau_\ga}\big|
\y^\ga_{\ti x}(-t/2)-\y^\ga_x(-t/2)\big|\di t\Big]\\[5pt]
\dis\quad=(\ffrac{1}{\ga}+1)L\ga(\ti x-x)=L(1+\ga)|\ti x-x|,
\ec
where $L$ is the Lipschitz constant of $p$ and we have used the
same exponentially distributed $\tau_\ga$ for $\y^\ga_x$ and
$\y^\ga_{\ti x}$.\qed

\subsection{Monotone and concave catalyzing functions}\label{concavesec}

In this section we prove that the log-Laplace operators $\Ui_\ga$ from
(\ref{Uiga}) map monotone functions into monotone functions, and
monotone concave functions into monotone concave functions. We do not
know if in general $\Ui_\ga$ maps concave functions into concave
functions.
 \bp{\bf(Preservation of monotonicity and concavity)}\label{moncon}
Let $\ga>0$. Then:\smallskip

\noi
{\bf (a)} If $f\in\Ci_+[0,1]$ is nondecreasing, then $\Ui_\ga f$ is
nondecreasing.\smallskip

\noi
{\bf (b)} If $f\in\Ci_+[0,1]$ is nondecreasing and concave, then
$\Ui_\ga f$ is nondecreasing and concave.
\ep
{\bf Proof} Our proof of Proposition~\ref{moncon} is in part based on
ideas from \cite[Appendix~A]{BCGH97}. The proof is quite long and will
depend on several lemmas. We remark that part~(a) can be proved in a
more elementary way using Lemma~\ref{monot}.

We recall some facts from Hille-Yosida theory. A linear operator $A$
on a Banach space $V$ is closable and its closure $\ov A$ generates a
strongly continuous contraction semigroup $(S_t)_{t\geq 0}$ if and only if
\be\ba{rl}\label{HilYos1}
{\rm (i)}&\Di(A)\mbox{ is dense},\\
{\rm (ii)}&A\mbox{ is dissipative},\\
{\rm (iii)}&\Ri(1-\al A)\mbox{ is dense for some, and hence for all }\al>0.
\ec
\detail{To see that condition (iii) holds for all $\al>0$ if it holds
for some $\al>0$, one needs the elementary fact that $\ov{\Ri(1-\al A)}
=\Ri(1-\al\ov A)$.}
Here, for any linear operator $B$ on $V$, $\Di(B)$ and $\Ri(B)$ denote
the domain and range of $B$, respectively. For each $\al>0$, the operator
$(1-\al\ov A):\Di(\ov A)\to V$ is a bijection and its inverse
$(1-\al\ov A)^{-1}:V\to\Di(\ov A)$ is a bounded linear operator, given by
\be
(1-\al\ov A)^{-1}u=\int_0^\infty\!\! S_tu\;\al^{-1}e^{-t/\al}\di t
\qquad(u\in V,\ \al>0).
\ee
If $E$ is a compact metrizable space and $\Ci(E)$ is the Banach space
of continuous real functions on $E$, equipped with the supremumnorm,
then a linear operator $A$ on $\Ci(E)$ is closable and its closure
$\ov A$ generates a Feller semigroup if and only if (see
\cite[Theorem~4.2.2 and remarks on page~166]{EK})
\be\ba{rl}\label{HilYos2}
{\rm (i)}&1\in\Di(\ov A)\mbox{ and }\ov A1=0,\\
{\rm (ii)}&\Di(A)\mbox{ is dense},\\
{\rm (iii)}&A\mbox{ satisfies the positive maximum principle},\\
{\rm (iv)}&\Ri(1-\al A)\mbox{ is dense for some, and hence for all }\al>0.
\ec
If $\ov A$ generates a Feller semigroup and $g\in\Ci(E)$, then the
operator $\ov A+g$ (with domain $\Di(\ov A+g):=\Di(\ov A)$) generates
a strongly continuous semigroup $(S^g_t)_{t\geq 0}$ on $\Ci(E)$. If
$g\leq 0$ then $(S^g_t)_{t\geq 0}$ is contractive. If $(\xi_t)_{t\geq 0}$
is the Feller process with generator $\ov A$, then one has the
Feynman-Kac representation
\be\label{FeyKac}
S^g_tu(x)=E^x[u(\xi(t))\ex{\int_0^tg(\xi(s))\di s}\big]
\qquad(t\geq 0,\ x\in E,\ g,u\in\Ci(E)).
\ee
Let $\Ci^{(n)}([0,1]^2)$ denote the space of continuous real functions
on $[0,1]^2$ whose partial derivatives up to $n$-th order exist and
are continuous on $[0,1]^2$ (including the boundary), and put
$\Ci^{(\infty)}([0,1]^2):=\bigcap_n\Ci^{(n)}([0,1]^2)$. Define a
linear operator $B$ on $\Ci([0,1]^2)$ with domain
$\Di(B):=\Ci^{(\infty)}([0,1]^2)$ by
\be\label{Bdef}
Bu(x,y):=y(1-y)\diff{y}u(x,y)+\ffrac{1}{\ga}(x-y)\dif{y}u(x,y).
\ee
Below, we will prove:
\bl{\bf(Feller semigroup)}\label{Felgen}
The closure in $\Ci([0,1]^2)$ of the operator $B$ generates a Feller
semigroup on $\Ci([0,1]^2)$.
\el
Write
\be\ba{c@{\,}c@{\,}l}
\Ci_+&:=&\dis\big\{u\in\Ci([0,1]^2):u\geq 0\big\},\\[5pt]
\Ci_{1+}&:=&\dis
\big\{u\in\Ci^{(1)}([0,1]^2):\dif{y}u,\dif{x}u\geq 0\big\},\\[5pt]
\Ci_{2+}&:=&\dis
\big\{u\in\Ci^{(2)}([0,1]^2):\diff{y}u,\difif{x}{y}u,\diff{x}u\geq 0\big\}.
\ec
Let $\ov\Si$ denote the closure of a set $\Si\sub\Ci([0,1]^2)$.
We need the following lemma.
\bl{\bf (Preserved classes)}\label{Lapcon}
Let $g\in\Ci([0,1]^2)$ and let $(S^g_t)_{t\geq 0}$ be the strongly
continuous semigroup with generator $\ov B+g$. Then, for each
$t\geq 0$:\smallskip

\noi
{\bf (a)} If $g\in\ov{\Ci_{1+}}$, then $S^g_t$ maps
$\ov{\Ci_+\cap\Ci_{1+}}$ into itself.\smallskip

\noi
{\bf (b)} If $g\in\ov{\Ci_{1+}\cap\Ci_{2+}}$, then $S^g_t$ maps
$\ov{\Ci_+\cap\Ci_{1+}\cap\Ci_{2+}}$ into itself.
\el
To see why Lemma~\ref{Lapcon} implies Proposition~\ref{moncon}, let
$(\x(t),\y(t))_{t\geq 0}$ denote the Feller process in $[0,1]^2$
generated by $\ov B$. It is easy to see that $\x(t)=\x(0)$ a.s.\ for
all $t\geq 0$. For fixed $\x(0)=x$, the process $(\y(t))_{t\geq 0}$ is
the diffusion given by the SDE (\ref{WFsde}). Therefore, by
Feynman-Kac, for each $g\in\Ci([0,1]^2)$,
\be\label{FK}
E^y\big[\ex{\int_0^tg(x,\y(s))\di s}\big]=S^g_t1(x,y),
\ee
where $1$ denotes the constant function $1\in\Ci([0,1]^2)$. By (\ref{Uiga}),
\be\label{Uiexpr}
\Ui_\ga f(x)=(\ffrac{1}{\ga}+1)\Big(1-\int\Ga^\ga_x(\di y)
E^y\big[\ex{-\int_0^{\tau_\ga}f(\y_x(s))\di s}\big]\Big)
\qquad(f\in\Ci_+[0,1]),
\ee
where $\Ga^\ga_x$ is the invariant law of $(\y(t))_{t\geq 0}$ from
Corollary~\ref{ergo} and $\tau_\ga$ is an exponential time with mean
$\ga$, independent of $(\y(t))_{t\geq 0}$. Setting $g(x,y):=-f(y)$ in
(\ref{FK}), using the ergodicity of $(\y(t))_{t\geq 0}$ (see
Corollary~\ref{ergo}), we find that for each $z\in[0,1]$ and $t\geq
0$,
\bc\label{limla}
\dis\int\Ga^\ga_x(\di y)E^y\big[\ex{-\int_0^tf(\y(s))\di s}\big]
&=&\dis\lim_{r\to\infty}\int P^z[\y(r)\in\di y]\,
E^y\big[\ex{-\int_0^tg(x,\y(s))\di s}\big]\\[7pt]
&=&\dis\lim_{r\to\infty}S^0_rS^g_t1(x,z).
\ec
It follows from Lemma~\ref{Lapcon} that for each fixed $r,t$, and $z$,
the function $x\mapsto S^0_rS^g_t1(x,z)$ is nondecreasing if $f$ is
nonincreasing, and nondecreasing and convex if $f$ is nonincreasing
and concave. Therefore, taking the expectation over the randomness of
$\tau_\ga$, the claims follow from (\ref{Uiexpr}) and
(\ref{limla}).\qed

\noi
We still need to prove Lemmas~\ref{Felgen} and \ref{Lapcon}.\med

\noi
{\bf Proof of Lemma~\ref{Felgen}} It is easy to see that the operator
$B$ from (\ref{Bdef}) is densely defined, satisfies the positive
maximum principle, and maps the constant function $1$ into $0$.
Therefore, by Hille-Yosida (\ref{HilYos2}), we must show that the
range $\Ri(1-\al B)$ is dense in $\Ci([0,1]^2)$ for some, and hence
for all $\al>0$. Let $\Pc_n$ denote the space of polynomials on
$[0,1]^2$ of $n$-th and lower order, i.e., the space of functions
$f:[0,1]^2\to\R$ of the form
\be
f(x,y)=\sum_{k,l\geq 0}a_{kl}\,x^ky^l\quad\mbox{ with $a_{k,l}=0$ for $k+l>n$.}
\ee
Set $\Pc_\infty:=\bigcup_n\Pc_n$. It is easy to see that $B$ maps the
space $\Pc_n$ into itself, for each $n\geq 0$. Since each $\Pc_n$ is
finite-dimensional, a simple argument (see
\cite[Proposition~1.3.5]{EK}) shows that the image of $\Pc_\infty$
under $1-\al B$ is dense in $\Ci([0,1]^2)$ for all but countably many,
and hence for all $\al>0$.\qed

\noi
As a first step towards proving Lemma~\ref{Lapcon}, we prove:
\bl{\bf(Smooth solutions to Laplace equation)}\label{smoothlap}
Let $\al>0$, $g\in\Ci^{(2)}([0,1])$, $g\leq 0$, $v\in\Ci([0,1]^2)$,
and assume that $u\in\Ci^{(\infty)}([0,1]^2)$ solves the Laplace equation
\be\label{Lap}
(1-\al(B+g))u=v.
\ee
{\bf (a)} If $g\in\Ci_{1+}$, then $v\in\Ci_{+}\cap\Ci_{1+}$ implies
$u\in\Ci_{+}\cap\Ci_{1+}$.\smallskip

\noi
{\bf (b)} If $g\in\Ci_{1+}\cap\Ci_{2+}$, then
$v\in\Ci_{+}\cap\Ci_{1+}\cap\Ci_{2+}$ implies
$u\in\Ci_{+}\cap\Ci_{1+}\cap\Ci_{2+}$.
\el
{\bf Proof} Let $u^y:=\dif{y}u$, $u^{xy}:=\difif{x}{y}u$, etc.\ denote
the partial derivatives of $u$ and similarly for $v$ and $g$, whenever
they exist. Set $c:=\frac{1}{\ga}$. Define linear operators $B'$ and
$B''$ on $\Ci([0,1]^2)$ with domains
$\Di(B')=\Di(B''):=\Ci^{(\infty)}([0,1]^2)$ by
\bc
B'&:=&\dis y(1-y)\diff{y}+\big(c(x-y)+2(\ffrac{1}{2}-y)\big)\dif{y},\\[5pt]
B''&:=&\dis y(1-y)\diff{y}+\big(c(x-y)+4(\ffrac{1}{2}-y)\big)\dif{y}.
\ec
Then
\be\ba{r@{\,}c@{\,}l@{\quad}r@{\,}c@{\,}l}\label{Bcommut}
\dif{y}Bu&=&(B'-c)u^y,&\dif{y}B'u&=&(B''-c-2)u^y,\\[5pt]
\dif{x}Bu&=&Bu^x+cu^y,&\dif{x}B'u&=&B'u^x+cu^y.
\ec
Therefore, it is easy to see that
\be\ba{rr@{\,}c@{\,}l}\label{partdi}
{\rm (i)}&(1-\al(B'-c+g))u^y&=&v^y+\al g^yu,\\
{\rm (ii)}&(1-\al(B+g))u^x&=&v^x+\al(cu^y+g^xu),\\
{\rm (iii)}&(1-\al(B''-2c-2+g))u^{yy}&=&v^{yy}+\al(2g^yu^y+g^{yy}u),\\
{\rm (iv)}&(1-\al(B'-c+g))u^{xy}&=&v^{xy}+\al(cu^{yy}+g^yu^x+g^{xy}u+g^xu^y),\\
{\rm (v)}&(1-\al(B+g))u^{xx}&=&v^{xx}+\al(2cu^{xy}+2g^xu^x+g^{xx}u),
\ec
where in (i) and (ii) we assume that $v\in\Ci^{(1)}([0,1]^2)$ and in
(iii)--(v) we assume that $v\in\Ci^{(2)}([0,1]^2)$. By
Lemma~\ref{Felgen}, the closure of the operator $B$ generates a Feller
processes in $[0,1]^2$. Exactly the same proof shows that $B'$ and
$B''$ also generate Feller processes on $[0,1]^2$. Therefore, by
Feynman-Kac, $u$ is nonnegative if $v$ is nonnegative and
$u^y,\ldots,u^{xx}$ are nonnegative if the right-hand sides of the
equations (i)--(v) are well-defined and nonnegative. (Instead of using
Feynman-Kac, this follows more elementarily from the fact that $B,B'$,
and $B''$ satisfy the positive maximum principle.) In particular, if
$g^y,g^x\geq 0$ and $v\in\Ci^{(1)}([0,1]^2)$, $v,v^y,v^x\geq 0$, then
it follows that $u,u^y,u^x\geq 0$. If moreover
$g^{yy},g^{xy},g^{xx}\geq 0$ and $v\in\Ci^{(2)}([0,1]^2)$,
$v^{yy},v^{xy},v^{yy}\geq 0$, then also $u^{yy},u^{xy},u^{yy}\geq
0$.\qed

\noi
In order to prove Lemma~\ref{Lapcon}, based on Lemma~\ref{smoothlap},
we will show that the Laplace equation (\ref{Lap}) has smooth
solutions $u$ for sufficiently many functions $v$. Here `suffiently
many' will mean dense in the topology of uniform convergence of
functions and their derivatives up to second order. To this aim, we
make $\Ci^{(2)}([0,1]^2)$ into a Banach space by equipping it with the
norm
\be\label{2norm}
\|u\|_{(2)}:=\|u\|+\|u^y\|+\|u^x\|+\|u^{yy}\|+2\|u^{xy}\|+\|u^{xx}\|.
\ee
Here, to reduce notation, we denote the supremumnorm by
$\|f\|:=\|f\|_\infty$. Note the factor 2 in the second term
from the right in (\ref{2norm}), which is crucial for the next key lemma.
\bl{\bf(Semigroup on twice diffferentiable functions)}\label{twicegen}
The closure in $\Ci^{(2)}([0,1]^2)$ of the operator $B$ generates a
strongly continuous contraction semigroup on $\Ci^{(2)}([0,1]^2)$.
\el
{\bf Proof} We must check the conditions (i)--(iii) from
(\ref{HilYos1}). It is well-known (see for example
\cite[Proposition~7.1 from the appendix]{EK}) that the space
$\Pc_\infty$ of polynomials is dense in $\Ci^{(2)}([0,1]^2)$.
Therefore $\Di(B)=\Ci^{(\infty)}([0,1]^2)$ is dense, and copying the
proof of Lemma~\ref{Felgen} we see that $\Ri(1-\al B)$ is dense for
all but countably many $\al$. To complete the proof, we must show that
$B$ is dissipative, i.e., that
\be
\|(1-\eps B)u\|_{(2)}\geq\|u\|_{(2)}\qquad
(\eps>0,\ u\in\Ci^{(\infty)}([0,1]^2)).
\ee
Using (\ref{Bcommut}), we calculate
\bc\label{doordif}
\dis\dif{y}(1-\eps B)u&=&\dis(1-\eps(B'-c))u^y,\\[5pt]
\dis\dif{x}(1-\eps B)u&=&\dis(1-\eps B)u^x-\eps cu^y,\\[5pt]
\dis\diff{y}(1-\eps B)u&=&\dis(1-\eps(B''-2c-2))u^{yy},\\[5pt]
\dis\difif{x}{y}(1-\eps B)u&=&\dis(1-\eps(B'-c))u^{xy}-\eps cu^{yy},\\[5pt]
\dis\diff{x}(1-\eps B)u&=&\dis(1-\eps B)u^{xx}-2\eps cu^{xy}.
\ec
Using the disipativity of $B,B'$, and $B''$ with respect to the
supremumnorm (which follows from the positive maximum principle) we
see that $\|(1-\eps(B'-c))u^y\|=(1+\eps c)\|(1-\frac{\eps}{1+\eps c}B)u^y\|
\geq(1+\eps c)\|u^y\|$ etc. We conclude therefore from
(\ref{doordif}) that
\bc
\dis\|(1-\eps B)u\|_{(2)}&\geq &\dis\|(1-\eps B)u\|
+\|(1-\eps(B'-c))u^y\|+\|(1-\eps B)u^x\|-\eps c\|u^y\|\\[5pt]
&&\dis+\|(1-\eps(B''-2c-2))u^{yy}\|+2\|(1-\eps(B'-c))u^{xy}\|
-2\eps c\|u^{yy}\|\\[5pt]
&&\dis+\|(1-\eps B)u^{xx}\|-2\eps c\|u^{xy}\|\\[5pt]
&\geq&\dis\|u\|+(1+\eps c)\|u^y\|+\|u^x\|-\eps c\|u^y\|\\[5pt]
&&\dis+(1+\eps(2c+2))\|u^{yy}\|+2(1+\eps c)\|u^{xy}\|-2\eps c\|u^{yy}\|\\[5pt]
&&\dis+\|u^{xx}\|-2\eps c\|u^{xy}\|\geq\|u\|_{(2)}
\ec
for each $\eps>0$, which shows that $B$ is dissipative with respect to
the norm $\|\cdot\|_{(2)}$.\qed

\noi
{\bf Proof of Lemma~\ref{Lapcon}} Let $g\in\Ci^{(2)}([0,1]^2)$. Then
$u\mapsto gu$ is a bounded operator on both $\Ci([0,1]^2)$ and
$\Ci^{(2)}([0,1]^2)$, so we can choose a $\la>0$ such that
\be
\|gu\|\leq\la\|u\|\quad\mbox{and}\quad\|gu\|_{(2)}\leq\la\|u\|_{(2)}
\ee
for all $u$ in $\Ci([0,1]^2)$ and $\Ci^{(2)}([0,1]^2)$, respectively.
Put $\ti g:=g-\la$. By Lemma~\ref{Felgen}, $\ov B+\ti g$ generates a
strongly continuous contraction semigroup $(S^{\ti g}_t)_{t\geq 0}
=(e^{-\la t}S^g_t)_{t\geq 0}$ on $\Ci([0,1]^2)$. Note that
$\Ri(1-\al(B+\ti g))$ is the space of all $v\in\Ci([0,1]^2)$ for which
the Laplace equation $(1-\al(B+\ti g))u=v$ has a solution
$u\in\Ci^{(\infty)}([0,1]^2)$. Therefore, by Lemma~\ref{smoothlap},
for each $\al>0$:
\be\ba{rl}\label{smoothpres}
{\rm (i)}&\mbox{If $g\in\Ci_{1+}$, then $(1-\al(\ov B+\ti g))^{-1}$
maps $\Ri(1-\al(B+\ti g))\cap\Ci_+\cap\Ci_{1+}$ into
$\Ci_+\cap\Ci_{1+}$.}\\[5pt]
{\rm (ii)}&\mbox{If $g\in\Ci_{1+}\cap\Ci_{2+}$, then
$(1-\al(\ov B+\ti g))^{-1}$ maps
$\Ri(1-\al(B+\ti g))\cap\Ci_+\cap\Ci_{1+}\cap\Ci_{2+}$}\\
&\mbox{into $\Ci_+\cap\Ci_{1+}\cap\Ci_{2+}$.}
\ec
By Lemma~\ref{twicegen}, the restriction of the semigroup
$(S^{\ti g}_t)_{t\geq 0}$ to $\Ci^{(2)}([0,1]^2)$ is strongly continuous and
contractive in the norm $\|\cdot\|_{(2)}$. Therefore, by Hille-Yosida
(\ref{HilYos1}), $\Ri(1-\al(B+\ti g))$ is dense in
$\Ci^{(2)}([0,1]^2)$ for each $\al>0$. It follows that
$\Ri(1-\al(B+\ti g))\cap\Ci_+\cap\Ci_{1+}$ is dense in
$\Ci_+\cap\Ci_{1+}$ and likewise $\Ri(1-\al(B+\ti
g))\cap\Ci_+\cap\Ci_{1+}\cap\Ci_{2+}$ is dense in
$\Ci_+\cap\Ci_{1+}\cap\Ci_{2+}$, both in the norm $\|\cdot\|_{(2)}$.
Note that we need density in the norm $\|\cdot\|_{(2)}$ here: if we
would only know that $\Ri(1-\al(B+\ti g))$ is a dense subset of
$\Ci([0,1]^2)$ in the norm $\|\cdot\|$, then $\Ri(1-\al(B+\ti
g))\cap\Ci_+\cap\Ci_{1+}$ might be empty. By approximation in the norm
$\|\cdot\|_{(2)}$ it follows from (\ref{smoothpres}) that:
\be\ba{rl}\label{C2pres}
{\rm (i)}&\mbox{If $g\in\Ci_{1+}$, then $(1-\al(\ov B+\ti g))^{-1}$
maps $\Ci_+\cap\Ci_{1+}$ into itself.}\\[5pt]
{\rm (ii)}&\mbox{If $g\in\Ci_{1+}\cap\Ci_{2+}$, then
$(1-\al(\ov B+\ti g))^{-1}$ maps $\Ci_+\cap\Ci_{1+}\cap\Ci_{2+}$ into itself.}
\ec
Using also continuity in the norm $\|\cdot\|$ we find that:
\be\ba{rl}\label{ovCpres}
{\rm (i)}&\mbox{If $g\in\Ci_{1+}$, then $(1-\al(\ov B+\ti g))^{-1}$
maps $\ov{\Ci_+\cap\Ci_{1+}}$ into itself.}\\[5pt]
{\rm (ii)}&\mbox{If $g\in\Ci_{1+}\cap\Ci_{2+}$, then
$(1-\al(\ov B+\ti g))^{-1}$ maps $\ov{\Ci_+\cap\Ci_{1+}\cap\Ci_{2+}}$
into itself.}
\ec
For $\eps>0$ let
\be
G_\eps:=\eps^{-1}\big((1-\eps(\ov B+\ti g))^{-1}-1\big)
\ee
be the Yosida approximation to $\ov B+\ti g$. Then
\be
e^{G_\eps t}=e^{-\eps^{-1}t}\sum_{n=0}^\infty\frac{t^n}{n!}
(1-\eps(\ov B+\ti g))^{-n}\qquad(t\geq 0),
\ee
and therefore, by (\ref{ovCpres}), for each $t\geq 0$:
\be\ba{rl}\label{epGep}
{\rm (i)}&\mbox{If $g\in\Ci_{1+}$, then $e^{G_\eps t}$ maps
$\ov{\Ci_+\cap\Ci_{1+}}$ into itself.}\\[5pt]
{\rm (ii)}&\mbox{If $g\in\Ci_{1+}\cap\Ci_{2+}$, then
$e^{G_\eps t}$ maps $\ov{\Ci_+\cap\Ci_{1+}\cap\Ci_{2+}}$ into itself.}
\ec
Finally
\be
e^{-\la t}S^g_tu=S^{\ti g}_tu=\lim_{\eps\to 0}\ex{G_\eps t}u
\qquad(t\geq 0,\ u\in\Ci([0,1]^2)),
\ee
so (\ref{epGep}) implies that for each $t\geq 0$:
\be\ba{rl}
{\rm (i)}&\mbox{If $g\in\Ci_{1+}$, then $S^g_t$ maps
$\ov{\Ci_+\cap\Ci_{1+}}$ into itself.}\\[5pt]
{\rm (ii)}&\mbox{If $g\in\Ci_{1+}\cap\Ci_{2+}$, then
$S^g_t$ maps $\ov{\Ci_+\cap\Ci_{1+}\cap\Ci_{2+}}$ into itself.}
\ec
Using the continuity of $S^g_t$ in $g$ (which follows from
Feynman-Kac (\ref{FeyKac})) we arrive at the statements in
Lemma~\ref{Lapcon}.\qed

\section{Convergence to a time-homogeneous process}\label{convsec}

\subsection{Convergence of certain Markov chains}\label{Markcon}

Section~\ref{convsec} is devoted to the proof of
Theorem~\ref{supercon}. In the present subsection, we start by
formulating a theorem about the convergence of certain Markov chains
to con\-tin\-uous-time processes. In Section~\ref{Brasuper} we
specialize to Poisson-cluster branching processes and superprocesses.
In Section~\ref{superproof}, finally, we carry out the necessary
calculations for the specific processes from Theorem~\ref{supercon}.

Let $E$ be a {\em compact} metrizable space. Let $A:\Di(A)\to\Ci(E)$
be an operator defined on a domain $\Di(A)\sub\Ci(E)$. We say that a
process $\y=(\y_t)_{t\geq 0}$ solves the martingale problem for $A$ if
$\y$ has sample paths in $\Di_E\half$ and for each $f\in\Di(A)$, the
process $(M^f_t)_{t\geq 0}$ given by
\be
M^f_t:=f(\y_t)-\int_0^tAf(\y_s)\di s\qquad(t\geq 0)
\ee
is a martingale with respect to the filtration generated by $\y$. We
say that existence (uniqueness) holds for the martingale problem for
$A$ if for each probability measure $\mu$ on $E$ there is at least one
(at most one (in law)) solution $\y$ to the martingale problem for $A$
with initial law $\Li(\y_0)=\mu$. If both existence and uniqueness
hold we say that the martingale problem is well-posed. For each $n\geq
0$, let $X^{(n)}=(X^{(n)}_0,\ldots,X^{(n)}_{m(n)})$ (with $1\leq
m(n)<\infty$) be a (time-inhomogeneous) Markov process in $E$ with
$k$-th step transition probabilities
\be
P_k(x,\di y)=P\big[X^{(n)}_k\in\di y\big|X^{(n)}_{k-1}=x\big]
\qquad(1\leq k\leq m(n)).
\ee
We assume that the $P_k$ are continuous probability kernels on $E$.
Let $(\eps^{(n)}_k)_{1\leq k\leq m(n)}$ be positive constants. Set
\be
A^{(n)}_kf(x):=(\eps^{(n)}_k)^{-1}\Big(\int_EP_k(x,\di y)f(y)-f(x)\Big)
\qquad(1\leq k\leq m(n),\ f\in\Ci(E)).
\ee
Define $t^{(n)}_0:=0$ and
\be\label{tdef}
t^{(n)}_k:=\sum_{l=1}^k\eps^{(n)}_l\qquad(1\leq k\leq m(n)),
\ee
and put
\be\label{kdef}
k^{(n)}(t):=\max\big\{k\;:\;0\leq k\leq m(n),\ t^{(n)}_k\leq t\big\}
\qquad(t\geq 0).
\ee
Define processes $\y^{(n)}=(\y^{(n)}_t)_{t\geq 0}$ with sample paths
in $\Di_E\half$ by
\be
\y^{(n)}_t:=X^{(n)}_{k^{(n)}(t)}\qquad(t\geq 0).
\ee
By definition, a space $\Ai$ of real functions is called an algebra
if $\Ai$ is a linear space and $f,g\in\Ai$ implies $fg\in\Ai$.
\bt{\bf(Convergence of Markov chains)}\label{genth}
Assume that $\Li(X^{(n)}_0)\Rightarrow\mu$ as $n\to\infty$ for some
probability law $\mu$ on $E$. Suppose that there exists at most one
(in law) solution to the martingale problem for $A$ with initial law
$\mu$. Assume that the linear span of $\Di(A)$ contains an algebra
that separates points. Assume that
\be\label{epscon}
{\rm (i)}\ \lim_{n\to\infty}\sum_{k=1}^{m(n)}\eps^{(n)}_k=\infty,
\qquad{\rm (ii)}\ \lim_{n\to\infty}\;\sup_{k:\ t^{(n)}_k\leq T}\eps^{(n)}_k=0,
\ee
and
\be\label{Acon}
\lim_{n\to\infty}\;\sup_{k:\ t^{(n)}_k\leq T}\big\|A^{(n)}_kf-Af\|_\infty=0
\qquad(f\in\Di(A))
\ee
for each $T>0$. Then there exists a unique solution $\y$ to the martingale
problem for $A$ with initial law $\mu$ and moreover
$\Li(\y^{(n)})\Rightarrow\Li(\y)$, where $\Rightarrow$
denotes weak convergence of probability measures on $\Di_E\half$.
\et
{\bf Proof} We apply \cite[Corollary~4.8.15]{EK}. Fix $f\in\Di(A)$.
We start by observing that
\be
f(X^{(n)}_k)-\sum_{i=1}^k\eps^{(n)}_iA^{(n)}_if(X^{(n)}_{i-1})
\qquad(0\leq k\leq m(n))
\ee
is a martingale with respect to the filtration generated by $X^{(n)}$
and therefore,
\be\label{Ymart}
f(\y^{(n)}_t)\;-\!\!\!\!\sum_{i=1}^{k^{(n)}(t)}
\eps^{(n)}_iA^{(n)}_if(\y^{(n)}_{t^{(n)}_{i-1}})\qquad(t\geq 0)
\ee
is a martingale with respect to the filtration generated by $\y^{(n)}$. Put
\be
\lfloor t\rfloor^{(n)}:=t^{(n)}_{k^{(n)}(t)}\qquad(t\geq 0)
\ee
and set
\be
\phi^{(n)}_t:=A^{(n)}_{k^{(n)}(t)+1}f(\y^{(n)}_{\lfloor t\rfloor^{(n)}})
1_{\{t<t^{(n)}_{m(n)}\}}\qquad(t\geq 0)
\ee
and
\be
\xi^{(n)}_t:=f(\y^{(n)}_t)+\int_{\lfloor t\rfloor^{(n)}}^t\phi^{(n)}_s\di s
\qquad(t\geq 0).
\ee
Then we can rewrite the martingale in (\ref{Ymart}) as
\be
\xi^{(n)}_t-\int_0^t\phi^{(n)}_s\di s.
\ee
By \cite[Corollary~4.8.15]{EK} and the compactness of the state space, it
suffices to check the following conditions on $\phi^{(n)}$ and $\xi^{(n)}$:
\be\ba{rl}\label{EKcondi}
{\rm (i)}&\dis\sup_{n\geq N}\;\sup_{t\leq T}E\big[|\xi^{(n)}_t|\big]
<\infty,\\[5pt]
{\rm (ii)}&\dis\sup_{n\geq N}\;\sup_{t\leq T}E\big[|\phi^{(n)}_t|\big]
<\infty,\\[5pt]
{\rm (iii)}&\dis\lim_{n\to\infty}
E\Big[\big(\xi^{(n)}_T-f(\y^{(n)}_T)\big)\prod_{i=1}^rh_i(\y^{(n)}_{s_i})\Big]
=0,\\[5pt]
{\rm (iv)}&\dis\lim_{n\to\infty}E\Big[\big(\phi^{(n)}_T-Af(\y^{(n)}_T)\big)
\prod_{i=1}^rh_i(\y^{(n)}_{s_i})\Big]=0,\\[5pt]
{\rm (v)}&\dis\lim_{n\to\infty}E\Big[\sup_{t\in\Q\cap[0,T]}\big|\xi^{(n)}_t
-f(\y^{(n)}_t)\big|\Big]=0,\\[5pt]
{\rm (vi)}&\dis\sup_{n\geq N}E\big[\|\phi^{(n)}\|_{p,T}\big]<\infty
\qquad\mbox{for some }p\in(1,\infty],
\ec
for some $N\geq 0$ and for each $T>0$, $r\geq 1$,
$0\leq s_1<\cdots<s_r\leq T$, and $h_1,\ldots,h_r\in\Hi\sub\Ci(E)$.
Here $\Hi$ is separating, i.e., $\int h\di\mu=\int h\di\nu$ for all
$h\in\Hi$ implies $\mu=\nu$ whenever $\mu,\nu$ are probability measures
on $E$. In (vi):
\be
\|g\|_{p,T}:=\Big(\int_0^T|g(t)|^p\di t\Big)^{1/p}\qquad(1\leq p<\infty)
\ee
and $\|g\|_{\infty,T}$ denotes the essential supremum of $g$ over $[0,T]$.

The conditions (\ref{EKcondi})~(i)--(vi) are implied by the stronger conditions
\be\ba{rl}\label{strongEK}
{\rm (i)}&\dis\lim_{n\to\infty}\,\sup_{0\leq t\leq T}\big\|\xi^{(n)}_t
-f(\y^{(n)}_t)\big\|_\infty=0,\\[5pt]
{\rm (ii)}&\dis\lim_{n\to\infty}\,\sup_{0\leq t\leq T}\big\|\phi^{(n)}_t
-Af(\y^{(n)}_t)\big\|_\infty=0,
\ec
where we denote the essential supremumnorm of a real-valued random variable
$X$ by $\|X\|_\infty:=\inf\{K\geq 0:|X|\leq K\ \as\}$. Condition
(\ref{strongEK})~(ii) is implied by (\ref{epscon})~(i) and (\ref{Acon}).
To see that also (\ref{strongEK})~(i) holds, set
\be
M_n:=\sup_{0\leq t\leq T}\big\|\phi^{(n)}_t\big\|_\infty,
\ee
and estimate
\be\label{Mnsup}
\sup_{0\leq t\leq T}\big\|\xi^{(n)}_t-f(\y^{(n)}_t)\big\|_\infty
\leq M_n\sup\{\eps^{(n)}_k\,:\,1\leq k\leq m(n),\ t^{(n)}_k\leq T\}.
\ee
Condition (\ref{strongEK})~(ii) implies that $\limsup_nM_n<\infty$
and therefore the right-hand side of (\ref{Mnsup}) tends to zero by
assumption (\ref{epscon})~(ii).\qed

\subsection{Convergence of certain branching processes}\label{Brasuper}

In this section we apply Theorem~\ref{genth} to certain branching
processes and superprocesses.

Throughout this section, $E$ is a compact metrizable space and
$A:\Di(A)\to\Ci(E)$ is a linear operator on $\Ci(E)$ such that the
closure $\ov A$ of $A$ generates a Feller process $\xi=(\xi_t)_{t\geq
  0}$ in $E$ with Feller semigroup $(P_t)_{t\geq 0}$ given by
$P_tf(x):=E^x[f(\xi_t)]$ ($t\geq 0,\ f\in\Ci(E)$).

Let $\al\in\Ci_+(E)$ and $\bet,f\in\Ci(E)$. By definition, a function
$t\mapsto u_t$ from $\half$ into $\Ci(E)$ is a {\em classical} solution
to the semilinear Cauchy problem
\be\label{gencauchy}
\left\{\ba{r@{\,}c@{\,}l}
\dif{t}u_t&=&\ov Au_t+\bet u_t-\al u_t^2\qquad(t\geq 0),\\[5pt]
u_0&=&f
\ea\right.
\ee
if $t\mapsto u_t$ is continuously differentiable (in $\Ci(E)$),
$u_t\in\Di(\ov A)$ for all $t\geq 0$, and (\ref{gencauchy}) holds. We
say that $u$ is a {\em mild} solution to (\ref{gencauchy}) if
$t\mapsto u_t$ is continuous and
\be\label{mild}
u_t=P_tf+\int_0^tP_{t-s}(\bet u_s-\al u_s^2)\di s\qquad(t\geq 0).
\ee
\bl{\bf(Mild and classical solutions)}\label{L:semlin}
Equation (\ref{gencauchy}) has a unique $\Ci_+(E)$-valued mild
solution $u$ for each $f\in\Ci_+(E)$, and $f>0$ implies that $u_t>0$
for all $t\geq 0$. If moreover $f\in\Di(\ov A)$ then $u$ is a
classical solution. For each $t\geq 0$, $u_t$ depends continuously on
$f\in\Ci_+(E)$.
\el
{\bf Proof} It follows from \cite[Theorems~6.1.2, 6.1.4, and
6.1.5]{Paz83} that for each $f\in\Ci(E)$, (\ref{gencauchy}) has a
unique solution $(u_t)_{0\leq t<T}$ up to an explosion time $T$, and
that this is a classical solution if $f\in\Di(\ov A)$. Moreover, $u_t$
depends continuously on $f$. Using comparison arguments based on the
fact that $\ov A$ satisfies the positive maximum principle (which
follows from Hille-Yosida (\ref{HilYos2})) one easily proves the other
statements; compare \cite[Lemmas~23 and 24]{FStrim}.\qed

\noi
We denote the (mild or classical) solution of (\ref{gencauchy}) by
$\Ui_tf:=u_t$; then $\Ui_t:\Ci_+(E)\to\Ci_+(E)$ are continuous
operators and $\Ui=(\Ui_t)_{t\geq 0}$ is a (nonlinear) semigroup on
$\Ci_+(E)$.

Since $E$ is compact, the spaces $\{\mu\in\Mi(E):\mu(E)\leq M\}$ are
compact for each $M\geq 0$. In particular, $\Mi(E)$ is locally
compact. We denote its one-point compactification by
$\Mi(E)_\infty=\Mi(E)\cup\{\infty\}$. We define functions
$F_f\in\Ci(\Mi(E)_{\infty})$ by $F_f(\infty):=0$ and
\be\label{Ffdef}
F_f(\mu):=\ex{-\li\mu,f\re}\qquad(f\in\Ci_+(E),\ f>0,\ \mu\in\Mi(E)).
\ee
We introduce an operator $\Gi$ with domain
\be
\Di(\Gi):=\{F_f:f\in\Di(A),\ f>0\},
\ee
given by $\Gi F_f(\infty):=0$ and
\be\label{Gidef}
\Gi F_f(\mu):=-\li\mu,Af+\bet f-\al f^2\re\,\ex{-\li\mu,f\re}
\qquad(\mu\in\Mi(E)).
\ee
Note that $\Gi F_f\in\Ci(\Mi(E)_{\infty})$ for all $F_f\in\Di(\Gi)$.
\bp{\bf($(\ov A,\al,\bet)$-superprocesses)}\label{superdef}
The martingale problem for the operator $\Gi$ is well-posed. The
solutions to this martingale problem define a Feller process
$\Yi=(\Yi_t)_{t\geq 0}$ in $\Mi(E)_\infty$ with continuous sample
paths, called the $(\ov A,\al,\bet)$-superprocess. If $\Yi_0=\infty$
then $\Yi_t=\infty$ for all $t\geq 0$. If $\Yi_0=\mu\in\Mi(E)$ then
\be\label{logLa}
E^\mu\big[\ex{-\li\Yi_t,f\re}\big]=\ex{-\li\mu,\Ui_tf\re}\qquad(f\in\Ci_+(E)).
\ee
\ep
{\bf Proof} Results of this type are well-known, see for example
\cite[Theorem~9.4.3]{EK}, \cite{Fit88}, and
\cite[Th\'eor\`eme~7]{ER91}. Since, however, it is not completely
straightforward to derive the proposition above from these references,
we give a concise autonomous proof of most of our statements. Only for
the continuity of sample paths we refer the reader to
\cite[Corollary~(4.7)]{Fit88} or \cite[Corollaire~9]{ER91}.

We are going to extend $\Gi$ to an operator $\hat\Gi$ that is linear
and satisfies the conditions of the Hille-Yosida Theorem
(\ref{HilYos2}). For any $\ga\in\Ci_+(E)$ and $\mu\in\Mi(E)$, let
${\rm Clust}_\ga(\mu)$ denote a random measure such that on
$\{\ga=0\}$, ${\rm Clust}_\ga(\mu)$ is equal to $\mu$, and on
$\{\ga>0\}$, ${\rm Clust}_\ga(\mu)$ is a Poisson cluster measure with
intensity $\frac{1}{\ga}\mu$ and cluster mechanism
$\Qi(x,\cdot)=\Li(\tau_{\ga(x)}\de_x)$, where $\tau_{\ga(x)}$ is
exponentially distributed with mean $\ga(x)$. It is not hard to see
that
\be
E\big[\ex{-\li{\rm Clust}_\ga(\mu),f\re}\big]=\ex{-\li\mu,\Vi_\ga f\re}
\qquad(f\in\Ci(E),\ f>0),
\ee
where $\Vi_\ga f(x):=(\frac{1}{f(x)}+\ga(x))^{-1}$.
\detail{Indeed, by (\ref{randclust}),
\[\ba{l}
\dis\Vi_\ga f(x)=\ga(x)^{-1}E\big[1-\ex{-\li\tau_{\ga(x)}\de_x,f\re}\big]
=\ga(x)^{-1}\Big(1-\int_0^\infty e^{-f(x)s}\ffrac{1}{\ga(x)}
e^{s/\ga(x)}\di s\Big)\\[5pt]
\quad\dis=\ga(x)^{-1}\Big(1-\ffrac{1}{\ga(x)}\int_0^\infty 
e^{-(f(x)+1/\ga(x))s}\di s\Big)=\ga(x)^{-1}
\Big(1-\frac{\frac{1}{\ga(x)}}{f(x)+\frac{1}{\ga(x)}}\Big)\\[5pt]
\quad\dis=\frac{f(x)\ga(x)^{-1}}{f(x)+\ga(x)^{-1}}.
\ea\]
}
Note that since $\Vi_\ga 1$ is bounded, the previously mentioned
Poisson cluster measure mentioned above is well-defined. By definition,
we put ${\rm Clust}_\ga(\infty):=\infty$.

Define a linear operator $\Gi_\al$ on $\Ci(\Mi(E))_\infty)$ by
\be\label{G1}
\Gi_\al F(\mu):=\lim_{\eps\to 0}\eps^{-1}
\big(E[F({\rm Clust}_{\eps\al}(\mu))]-F(\mu)\big)
\ee
with as domain $\Di(\Gi_\al)$ the space of all $F\in\Ci(\Mi(E)_\infty)$
for which the limit exists. Define a linear operator $\Gi_\bet$ by
\be\label{G2}
\Gi_\bet F(\mu):=\lim_{\eps\to 0}\eps^{-1}\big(F((1+\eps\bet)\mu)-F(\mu)\big)
\ee
with domain $\Di(\Gi_\bet):=\Ci(\Mi(E))_\infty)$. Define
$P^\ast_t:\Mi(E)_\infty\to\Mi(E)_\infty$ by
$\li P^\ast_t\mu,f\re:=\li\mu,P_tf\re$ $(t\geq 0,\ f\in\Ci(E),\ \mu\in\Mi(E))$
and $P^\ast_t\infty:=\infty$ ($t\geq 0$). Finally, let $\Gi_{\ov A}$
be the linear operator on $\Ci(\Mi(E))_\infty)$ defined by
\be\label{G3}
\Gi_{\ov A}F(\mu):=\lim_{\eps\to 0}\eps^{-1}\big(F(P^\ast_\eps\mu)-F(\mu)\big),
\ee
with as domain $\Di(\Gi_{\ov A})$ the space of all $F$ for which the
limit exists. Define an operator $\hat\Gi$ by
\be
\hat\Gi:=\Gi_\al+\Gi_\bet+\Gi_{\ov A},
\ee
with domain $\Di(\hat\Gi):=\Di(\Gi_\al)\cap\Di(\Gi_{\ov A})$.
If $f\in\Di(\ov A)$, $f>0$, and $F_f$ is as in (\ref{Ffdef}),
then it is not hard to see that $\hat\Gi F_f(\infty)=0$ and
\be
\hat\Gi F_f(\mu):=-\li\mu,\ov Af+\bet f-\al f^2\re\,\ex{-\li\mu,f\re}
\qquad(\mu\in\Mi(E)).
\ee
\detail{$\dif{t}\Vi_{t\al}f(x)=\dif{t}(\frac{1}{f(x)}+t\al(x))^{-1}
=-\al(x)(\frac{1}{f(x)}+t\al(x))^{-2}=-\al(x)\Vi_{t\al}f(x)$.}
In particular, $\hat\Gi$ extends the operator $\Gi$ from
(\ref{Gidef}). Since $\Di(\ov A)$ is dense in $\Ci(E)$, it is easy to
see that $\{F_f:f\in\Di(\ov A),\ f>0\}$ is dense in
$\Ci(\Mi(E)_\infty)$. Hence $\Di(\hat\Gi)$ is dense. Using
(\ref{G1})--(\ref{G3}) it is not hard to show that $\hat\Gi$ satisfies
the positive maximum principle. Moreover, by Lemma~\ref{L:semlin}, for
$f\in\Di(\ov A)$ with $f>0$, the function $t\mapsto F_{\Ui_tf}$ from
$\half$ into $\Ci(\Mi(E)_\infty)$ is continuously differentiable,
satisfies $F_{\Ui_tf}\in\Di(\hat\Gi)$ for all $t\geq 0$, and
\be\label{FC}
\dif{t}F_{\Ui_tf}=\hat\Gi F_{\Ui_tf}\qquad(t\geq 0).
\ee
{F}rom this it is not hard to see that $\hat\Gi$ also satisfies
condition~(\ref{HilYos2})~(ii), so the closure of $\hat\Gi$ generates
a Feller semigroup $(S_t)_{t\geq 0}$ on $\Ci(\Mi(E)_\infty)$. It is
easy to see that $S_tF_f=F_{\Ui_tf}$ $(t\geq 0)$. By
\cite[Theorem~4.2.7]{EK}, this semigroup corresponds to a Feller
process $\Yi$ with cadlag sample paths in $\Mi(E)_\infty$. This means
that $E^\mu[F_f(\Yi_t)]=F_{\Ui_tf}(\mu)$ for all $f\in\Di(\ov A)$ with
$f>0$. If $\mu=\infty$ this shows that $\Yi_t=\infty$ for all $t\geq
0$. If $\mu\in\Mi(E)$ we obtain (\ref{logLa}) for $f\in\Di(\ov A)$,
$f>0$; the general case follows by approximation.\qed

\noi
Now let $(q_\eps)_{\eps>0}$ be continuous weight functions and let
$(\Qi_\eps)_{\eps>0}$ be continuous cluster mechanisms on $E$. Assume that
\be
Z_\eps(x):=\int\Qi_\eps(x,\di\chi)\li\chi,1\re<\infty\qquad(x\in E)
\ee
and define probability kernels $K_\eps$ on $E$ by
\be
\int K_\eps(x,\di y)f(y):=\frac{1}{Z_\eps(x)}
\int\Qi_\eps(x,\di\chi)\li\chi,f\re\qquad(f\in B(E)).
\ee
For each $n\geq 0$, let $(\eps^{(n)}_k)_{1\leq k\leq m(n)}$ (with
$1\leq m(n)<\infty$) be positive constants. Let
$\Xc^{(n)}=(\Xc^{(n)}_0,\ldots,\Xc^{(n)}_{m(n)})$ be a Poisson-cluster
branching process with weight functions
$q_{\eps^{(n)}_1},\ldots,q_{\eps^{(n)}_{m(n)}}$ and cluster mechanisms
$\Qi_{\eps^{(n)}_1},\ldots,\Qi_{\eps^{(n)}_{m(n)}}$. Define
$t^{(n)}_k$ and $k^{(n)}(t)$ as in (\ref{tdef})--(\ref{kdef}). Define
processes $\Yi^{(n)}$ by
\be
\Yi^{(n)}_t:=\Xc^{(n)}_{k^{(n)}(t)}\qquad(t\geq 0).
\ee
\bt{\bf(Convergence of Poisson-cluster branching processes)}\label{bratosup}
Assume that $\Li(\Xc^{(n)}_0)\Rightarrow\rho$ as $n\to\infty$ for some
probability law $\rho$ on $\Mi(E)$. Suppose that the constants
$\eps^{(n)}_k$ fulfill (\ref{epscon}). Assume that
\be\ba{rr@{\,}c@{\,}l}\label{Qlim}
{\rm (i)}&\dis q_\eps(x)\int\Qi_\eps(x,\di\chi)\li\chi,1\re&
=&\dis 1+\eps\bet(x)+o(\eps),\\[5pt]
{\rm (ii)}&\dis q_\eps(x)\int\Qi_\eps(x,\di\chi)\li\chi,1\re^2&
=&\dis\eps\,2\al(x)+o(\eps),\\[5pt]
{\rm (iii)}&\dis q_\eps(x)\int
\Qi_\eps(x,\di\chi)\li\chi,1\re^21_{\{\li\chi,1\re>\de\}}&=&o(\eps)
\ec
for each $\de>0$, and
\be\label{Klim}
\int K_\eps(x,\di y)f(y)=f(x)+\eps Af(x)+o(\eps)
\ee
for each $f\in\Di(A)$, uniformly in $x$ as $\eps\to 0$.
Then $\Li(\Yi^{(n)})\Rightarrow\Li(\Yi)$, where $\Yi$ is the
$(\ov A,\al,\bet)$-superprocess with initial law $\rho$.
\et
Here $\Rightarrow$ denotes weak convergence of probability
measures on $\Di_{\Mi(E)}\half$.\med

\noi
{\bf Proof} We apply Theorem~\ref{genth} to the operator $\Gi$, where
we use the fact that if we view $\Mi_1(\Di_{\Mi(E)}\half)$ as a
subspace of $\Mi_1(\Di_{\Mi(E)_\infty}\half)$ (note the
compactification), equipped with the topology of weak convergence,
then the induced topology on $\Mi_1(\Di_{\Mi(E)}\half)$ is again the
topology of weak convergence.

By Proposition~\ref{superdef}, solutions to the martingale problem for
$\Gi$ are unique. Since $F_fF_g=F_{f+g}$ and $\Di(A)$ is a linear
space, the linear span of the domain of $\Gi$ is an algebra. Using the
fact that $\Di(A)$ is dense in $\Ci(E)$ we see that this algebra
separates points. Therefore, we are left with the task to check
(\ref{Acon}).

Define $\Ui_\eps:\Ci_+(E)\to\Ci_+(E)$ by
\be
\Ui_\eps f(x):=q_\eps(x)\int\Qi_\eps(x,\di\chi)
\big(1-\ex{-\li\chi,f\re}\big)\qquad(x\in E,\ f\in\Ci_+[0,1],\ f>0,\ \eps>0),
\ee
and define transition probabilities $P_\eps(\mu,\di\nu)$ on
$\Mi(E)_\infty$ by $P_\eps(\infty,\,\cdot\,):=\de_\infty$ and
\be\label{pepstr}
\int P_\eps(\mu,\di\nu)\ex{-\li\nu,f\re}=\ex{-\li\mu,\Ui_\eps f\re}.
\ee
We will show that 
\be\label{Uicon}
\lim_{\eps\to 0}\big\|\eps^{-1}(\Ui_\eps f-f)
-(Af+\bet f-\al f^2)\big\|_\infty=0\qquad(f\in\Di(A),\ f>0).
\ee
Together with (\ref{pepstr}) this implies that
\be\label{Gicon}
\int P_\eps(\mu,\di\nu)F_f(\nu)=F_f(\mu)+\eps\Gi F_f(\mu)+o(\eps)
\qquad(f\in\Di(A),\ f>0),
\ee
uniformly in $\mu\in\Mi(E)_\infty$ as $\eps\to 0$. Therefore, the
result follows from Theorem~\ref{genth}.

It remains to prove (\ref{Uicon}). Set $g(z):=1-z+\frac{1}{2}z^2-e^{-z}$
$(z\geq 0)$ and write
\be\label{123}
\Ui_\eps f(x)=q_\eps(x)\int\Qi_\eps(x,\di\chi)\big(\li\chi,f\re
-\ffrac{1}{2}\li\chi,f\re^2+g(\li\chi,f\re)\big).
\ee
Since
\be
g(z)=\int_0^z\di y\int_0^y\di x\int_0^x\di t\,e^{-t}\qquad(z\geq 0),
\ee
it is easy to see that $g$ is nondecreasing on $\half$ and
(since $0\leq e^{-t}\leq 1$ and $\int_0^x\di t\,e^{-t}\leq 1$)
\be
0\leq g(z)\leq\ffrac{1}{2}z^2\wedge\ffrac{1}{6}z^3\qquad(z\geq 0).
\ee
Using these facts and (\ref{Qlim})~(ii) and (iii), we find that
\be\ba{l}
\dis q_\eps(x)\int\Qi_\eps(x,\di\chi)g(\li\chi,f\re)\\[5pt]
\dis\quad\leq\|f\|_\infty q_\eps(x)\Big\{\int\Qi_\eps(x,\di\chi)
g(\li\chi,1\re)1_{\{\li\chi,1\re\leq\de\}}+\int\Qi_\eps(x,\di\chi)
g(\li\chi,1\re)1_{\{\li\chi,1\re>\de\}}\Big\}\\[5pt]
\dis\qquad\leq\|f\|_\infty q_\eps(x)\Big\{\ffrac{1}{6}\de
\int\Qi_\eps(x,\di\chi)\li\chi,1\re^21_{\{\li\chi,1\re\leq\de\}}
+\ffrac{1}{2}\int\Qi_\eps(x,\di\chi)
\li\chi,1\re^21_{\{\li\chi,1\re>\de\}}\Big\}\\[10pt]
\dis\qquad=\ffrac{1}{6}\de\|f\|_\infty\big(\eps\,2\al(x)+o(\eps)\big)+o(\eps).
\ec
Since this holds for any $\de>0$, we conclude that
\be\label{restest}
q_\eps(x)\int\Qi_\eps(x,\di\chi)g(\li\chi,f\re)=o(\eps)
\ee
uniformly in $x$ as $\eps\to 0$. By (\ref{Qlim})~(i) and (\ref{Klim}),
\be\ba{l}\label{1est}
\dis q_\eps(x)\int\Qi_\eps(x,\di\chi)\li\chi,f\re=\Big(q_\eps(x)
\int\Qi_\eps(x,\di\chi)\li\chi,1\re\Big)
\Big(\int K_\eps(x,\di y)f(y)\Big)\\[5pt]
\dis\qquad=\big(1+\eps\bet(x)+o(\eps)\big)\big(f(x)
+\eps Af(x)+o(\eps)\big)\\[5pt]
\dis\qquad=f(x)+\eps\bet(x)f(x)+\eps Af(x)+o(\eps).
\ec
Finally, write
\be\ba{l}\label{pre2est}
\dis q_\eps(x)\int\Qi_\eps(x,\di\chi)\li\chi,f\re^2\\[5pt]
\dis\qquad=q_\eps(x)\int\Qi_\eps(x,\di\chi)\big(\li\chi,f(x)\re^2
+2\li\chi,f(x)\re\li\chi,f-f(x)\re+\li\chi,f-f(x)\re^2\big).
\ec
Then, by (\ref{Qlim})~(ii),
\be\label{fkwa}
q_\eps(x)\int\Qi_\eps(x,\di\chi)\li\chi,f(x)\re^2=f(x)^2
\big(\eps\,2\al(x)+o(\eps)\big).
\ee
We will prove that
\be\label{fmin1}
q_\eps(x)\int\Qi_\eps(x,\di\chi)\li\chi,f-f(x)\re^2=o(\eps).
\ee
Then, by H\"older's inequality, (\ref{Qlim})~(ii), and (\ref{fmin1}),
\be\ba{l}\label{f1f}
\dis\big|q_\eps(x)\int\Qi_\eps(x,\di\chi)\li\chi,f
-f(x)\re\li\chi,f(x)\re\big|\\[5pt]
\dis\qquad\leq\Big(q_\eps(x)
\int\Qi_\eps(x,\di\chi)\li\chi,f-f(x)\re^2\Big)^{1/2}
\Big(q_\eps(x)\int\Qi_\eps(x,\di\chi)\li\chi,f(x)\re^2\Big)^{1/2}\\[5pt]
\dis\qquad\leq\big(o(\eps)(2\al(x)\eps+o(\eps))\big)^{1/2}=o(\eps).
\ec
Inserting (\ref{fkwa}), (\ref{fmin1}) and (\ref{f1f}) into (\ref{pre2est})
we find that
\be\label{2est}
q_\eps(x)\int\Qi_\eps(x,\di\chi)\li\chi,f\re^2=\eps\,2\al(x)f(x)^2+o(\eps).
\ee
Inserting (\ref{restest}), (\ref{1est}) and (\ref{2est}) into (\ref{123}),
we arrive at (\ref{Uicon}). We still need to prove (\ref{fmin1}). To this
aim, we estimate, using (\ref{1est}),
\be\ba{l}
\dis q_\eps(x)\int\Qi_\eps(x,\di\chi)\li\chi,f-f(x)\re^2
1_{\{\li\chi,1\re\leq\de\}}\\[5pt]
\dis\qquad\leq\de\|f-f(x)\|_\infty q_\eps(x)
\int\Qi_\eps(x,\di\chi)\li\chi,f-f(x)\re\\[5pt]
\dis\qquad=\de\|f-f(x)\|_\infty\big(\eps Af(x)+o(\eps)\big)
\ec
and, using (\ref{Qlim})~(iii),
\be\ba{l}
\dis q_\eps(x)\int\Qi_\eps(x,\di\chi)\li\chi,f-f(x)\re^2
1_{\{\li\chi,1\re>\de\}}\\[5pt]
\dis\qquad\leq\|f-f(x)\|_\infty q_\eps(x)
\int\Qi_\eps(x,\di\chi)\li\chi,1\re^21_{\{\li\chi,1\re>\de\}}=o(\eps).
\ec
It follows that
\be
q_\eps(x)\int\Qi_\eps(x,\di\chi)\li\chi,f-f(x)\re^2\leq
\de\eps\|f-f(x)\|_\infty Af(x)+o(\eps)
\ee
for any $\de>0$. This implies (\ref{fmin1}) and completes
the proof of (\ref{Uicon}).\qed

\subsection{Application to the renormalization branching
process}\label{superproof}

{\bf Proof of Theorem~\ref{supercon}~(a)} For any
$f_0,\ldots,f_k\in\Ci_+[0,1]$ one has
\be\ba{l}\label{exexpl}
E\big[\ex{-\li\Xc_{-n},f_0\re}\cdots\ex{-\li\Xc_{-n+k},f_k\re}\big]\\[5pt]
\dis\quad=E\big[\ex{-\li\Xc_{-n},f_0\re}\cdots
\ex{-\li\Xc_{-n+k-1},f_{k-1}+\Ui_{\ga_{n-k}}f_k\re}\big]\\[5pt]
\dis\quad=\cdots=\dis\quad E\big[\ex{-\li\Xc_{-n},g_k\re}\big],
\ec
where we define inductively
\be
g_0:=f_k\quad\mbox{and}\quad g_{m+1}:=f_{k-m-1}+\Ui_{\ga_{n-k+m}}g_m.
\ee
By the compactness of $[0,1]$ and Corollary~\ref{Ugacont}, the map
$(\ga,f)\mapsto\Ui_\ga f$ from
$(0,\infty)\times\Ci_+[0,1]$ to $\Ci_+[0,1]$ (equipped with the
supremumnorm) is continuous. Using this fact and
(\ref{exexpl}) we find that
\be
E\big[\ex{-\li\Xc_{-n},f_0\re}\cdots\ex{-\li\Xc_{-n+k},f_k\re}\big]\asto{n}
E\big[\ex{-\li\Yi^{\ga^\ast}_{-n},f_0\re}\cdots
\ex{-\li\Yi^{\ga^\ast}_{-n+k},f_k\re}\big].
\ee
Since $f_1,\ldots,f_k$ are arbitrary, (\ref{disco}) follows.\qed

\noi
{\bf Proof of Theorem~\ref{supercon}~(b)} We apply
Theorem~\ref{bratosup} to the weight functions $q_\ga$ and cluster
mechanisms $\Qi_\ga$ from (\ref{Qqdef}) and to $A_{\rm WF}
=x(1-x)\diff{x}$ with domain $\Di(A_{\rm WF})=\Ci^{(2)}[0,1]$,
and $\al=\bet=1$. It is well-known that $\ov A_{\rm WF}$ generates a
Feller semigroup \cite[Theorem~8.2.8]{EK}. We observe that
\be
\int\Qi_\ga(x,\di\chi)\li\chi,f\re
=E\big[2\int_0^{\tau_\ga}\!f(\y^\ga_x(-t))\big]
=2E[\tau_\ga]E\big[f(\y^\ga_x(0))\big]
=\ga\int\Ga^\ga_x(\di y)f(y),
\ee
where $\Ga^\ga_x$ is the equilibrium law of the process $\y^\ga_x$
from Corollary~\ref{ergo}. It follows from (\ref{WFmoments}) that
\be\ba{rr@{\,}c@{\,}l}
{\rm (i)}&\dis\int\Ga^\ga_x(\di y)(y-x)&=&0,\\[5pt]
{\rm (ii)}&\dis\int\Ga^\ga_x(\di y)(y-x)^2&=&
\dis\frac{\ga x(1-x)}{1+\ga},\\[5pt]
{\rm (iii)}&\dis\int\Ga^\ga_x(\di y)(y-x)^4&=&\dis O(\ga^2),
\ec
uniformly in $x$ as $\ga\to 0$.
\detail{Indeed, writing $(y-x)^4=y^4-4xy^3+6x^2y^2-4x^3y+x^4$, we find that
\[\ba{r@{\,}c@{\,}l}
\dis\int\Ga^\ga_x(\di y)(y-x)^4
&=&\dis\big\{x(x+\ga)(x+2\ga)(x+3\ga)-4x^2(x+\ga)(x+2\ga)(1+3\ga)\\
&&\dis\phantom{\big\{}+6x^3(x+\ga)(1+2\ga)(1+3\ga)
-3x^4(1+\ga)(1+2\ga)(1+3\ga)\big\}\\
&&\dis\phantom{\big\{}\cdot\big((1+\ga)(1+2\ga)(1+3\ga)\big)^{-1}.
\ea\]
Collecting the terms with different powers of $\ga$, we find $0\cdot\ga^0$ and
\[\ba{l}
\dis\ga\big\{6x^3-4(3x^3+3x^4)+6(x^3+5x^4)-3\cdot 6x^4\big\}\\
\dis\quad=\ga\big\{(6-12+6)x^3+(-12+30-24+6)x^4\big\}=0\cdot\ga^1.
\ea\]}
Therefore, for any $\de>0$,
\be\ba{rr@{\,}c@{\,}l}
{\rm (i)}&\dis\int\Ga^\ga_x(\di y)(y-x)&=&0,\\[8pt]
{\rm (ii)}&\dis\int\Ga^\ga_x(\di y)(y-x)^2&=&\dis\ga x(1-x)+o(\ga),\\[8pt]
{\rm (iii)}&\dis\int\Ga^\ga_x(\di y)1_{\{|y-x|>\de\}}&=&\dis o(\ga),
\ec
uniformly in $x$ as $\ga\to 0$. Consequently, a Taylor expansion of
$f$ around $x$ yields
\be
\int\Ga^\ga_x(\di y)f(x)=f(x)+\ga\ffrac{1}{2}x(1-x)\diff{x}f(x)
+o(\ga)\qquad(f\in\Ci^{(2)}[0,1]),
\ee
uniformly in $x$ as $\ga\to 0$. (For details, in particular the
uniformity in $x$, see for example Proposition~\cite[B.1.1]{Swa99}.)
This shows that condition (\ref{Klim}) is satisfied. Moreover,
\be\ba{l}
\dis\int\Qi_\ga(x,\di\chi)\li\chi,1\re=E[2\tau_\ga]=\ga,\\[5pt]
\dis\int\Qi_\ga(x,\di\chi)\li\chi,1\re^2=E[(2\tau_\ga)^2]
=\int_0^\infty z^2\ffrac{1}{\ga}e^{-z/\ga}\di z=2\ga^2,\\[5pt]
\dis\int\Qi_\ga(x,\di\chi)\li\chi,1\re^3=E[(2\tau_\ga)^3]
=\int_0^\infty z^3\ffrac{1}{\ga}e^{-z/\ga}\di z=6\ga^3,
\ec
which, using the fact that $q_\ga=(\frac{1}{\ga}+1)$, gives
\be\ba{l}
\dis q_\ga\int\Qi_\ga(x,\di\chi)\li\chi,1\re=1+\ga,\\[10pt]
\dis q_\ga\int\Qi_\ga(x,\di\chi)\li\chi,1\re^2=2\ga+o(\ga),\\[10pt]
\dis q_\ga\int\Qi_\ga(x,\di\chi)\li\chi,1\re^3=o(\ga).
\ec
This shows that (\ref{Qlim}) is fulfilled. In particular,
\be
q_\ga\int\Qi_\ga(x,\di\chi)\li\chi,1\re^21_{\{\li\chi,1\re>\de\}}
\leq\de^{-1}q_\ga\int\Qi_\ga(x,\di\chi)\li\chi,1\re^3=o(\ga)
\ee
for all $\de>0$.\qed

\section{Embedded particle systems}\label{embsec}

In this section we use embedded particle systems to prove
Proposition~\ref{Uit}. An essential ingredient in the proofs
is Proposition~\ref{00hom}~(a), which will be proved in the Section~\ref{00sec}.

\subsection{Weighting and Poissonization}\label{Poissec}

{\bf Proof of Proposition~\ref{weightprop}} Obviously
$q^h_k\in\Ci_+(E^h)$ for each $k=1,\ldots,n$. Since $h\in\Ci_+(E)$ and
$h$ is bounded, it is easy to see that the map $\mu\mapsto h\mu$ from
$\Mi(E)$ into $\Mi(E^h)$ is continuous, and therefore the cluster
mechanisms defined in (\ref{Qih}) are continuous. Since
\be
\Ui^h_kf(x)=\frac{q_k(x)}{h(x)}E\big[1-\ex{-\li h\Zi_x,f\re}\big]
=\frac{\Ui_k(hf)(x)}{h(x)}\qquad(x\in E^h,\ f\in B_+(E^h)),
\ee
formula (\ref{htrafo1}) holds on $E^h$. To see that (\ref{htrafo1})
holds on $E\beh E^h$, note that by assumption $\Ui_kh\leq Kh$ for some
$K<\infty$, so if $x\in E\beh E^h$, then $\Ui_kh(x)=0$. By
monotonicity also $\Ui_k(hf)(x)=0$, while $h\Ui^h_kf(x)=0$ by
definition. Since $\sup_{x\in E^h}\Ui^h_k1(x)=\sup_{x\in E^h}
\frac{\Ui_kh(x)}{h(x)}\leq K<\infty$, the log-Laplace operators
$\Ui^h_k$ satisfy (\ref{fincon}). If $\Xc$ is started in an initial
state $\Xc_0$, then the Poisson-cluster branching process $\Xc^h$ with
log-Laplace operators $\Ui^h_1,\ldots,\Ui^h_n$ started in
$\Xc^h_0=h\Xc_0$ satisfies
\bc
\dis E\big[\ex{-\li h\Xc_k,f\re}\big]&
=&\dis E\big[\ex{-\li\Xc_0,\Ui_1\circ\cdots\circ\Ui_k(hf)\re}\big]\\[5pt]
&=&\dis E\big[\ex{-\li\Xc_0,h\Ui^h_1\circ\cdots\circ\Ui^h_k(f)\re}\big]
=E\big[\ex{-\li\Xc^h_k,f\re}\big]\qquad(f\in B_+(E^h)),
\ec
which proves (\ref{weighting}).\qed

\noi
{\bf Proof of Proposition~\ref{Poisprop}} We start by noting that by
(\ref{Vk}),
\be\label{Uiintu}
\Ui_kf(x)=q(x)E\big[1-\ex{-\li\Zi^k_x,f\re}\big]
=q_k(x)P[\Pois(f\Zi^k_x)\neq 0]\qquad(x\in E,\ f\in B_+(E)).
\ee
Into (\ref{Zxdef}), we insert
\be\ba{l}\label{condpois}
\dis P\big[\Pois(h\Zi^k_x)\in\cdot\,\big]\\[5pt]
\dis\quad=P\big[\Pois(h\Zi^k_x)\in\cdot\,\big|\,\Pois(h\Zi^k_x)\neq 0\big]
P[\Pois(h\Zi^k_x)\neq 0]+\de_0P[\Pois(h\Zi^k_x)=0].
\ec
Here and in similar formulas below, if in a conditional probability
the symbol $\Pois(\,\cdot\,)$ occurs twice with the same argument,
then it always refers to the same random variable (and not to
independent Poisson point measures with the same intensity, for
example). Using moreover (\ref{Uiintu}) we can rewrite (\ref{Zxdef})
as
\be\label{Zxfo}
Q^h_k(x,\,\cdot\,)=\frac{\Ui_k h(x)}{h(x)}
P\big[\Pois(h\Zi^k_x)\in\cdot\,\big|\,\Pois(h\Zi^k_x)\neq 0\big]
+\frac{h(x)-\Ui_k h(x)}{h(x)}\de_0(\,\cdot\,).
\ee
In particular, since we are assuming that $h$ is $\Ui_k$-subharmonic,
this shows that $Q^h_k(x,\,\cdot\,)$ is a probability measure. Let
$X^h$ be the branching particle system with offspring mechanisms
$Q^h_1,\ldots,Q^h_k$. Let $Z^{h,k}_x$ be random variables such that
$\Li(Z^{h,k}_x)=Q^h_k(x,\,\cdot\,)$. Then, by (\ref{VV}),
(\ref{Zxdef}), (\ref{elrel}), and (\ref{Uiintu}),
\be\ba{l}
\dis U^h_kf(x)=P[\Thin_f(Z^{h,k}_x)\neq 0]
=\frac{q_k(x)}{h(x)}P[\Thin_f(\Pois(h\Zi^k_x))\neq 0]\\[5pt]
\dis\qquad=\frac{q_k(x)}{h(x)}P[\Pois(hf\Zi^k_x)\neq 0]
=\frac{1}{h(x)}\Ui_k(hf)(x)\qquad(x\in E^h).
\ec
If $x\in E\beh E^h$, then $\Ui_k(hf)(x)\leq\Ui_k(h)(x)\leq
h(x)=0=:h\Ui^h(f)(x)$. This proves (\ref{htrafo}). To see that $Q^h_k$
is a {\em continuous} offspring mechanism, by
\cite[Theorem~4.2]{Kal76} it suffices to show that $x\mapsto\int
Q^h_k(x,\di\nu)\ex{-\li\nu,g\re}$ is continuous for all bounded
$g\in\Ci_+(E^h)$. Indeed, setting $f:=1-e^{-g}$, one has $\int
Q^h_k(x,\di\nu)\ex{-\li\nu,g\re}=\int
Q^h_k(x,\di\nu)(1-f)^\nu=1-\Ui^h_kf(x)=1-\Ui_k(hf)(x)/h(x)$ which is
continuous on $E^h$ by the continuity of $q_k$ and $\Qi_k$.

To see that also (\ref{Poissonization}) holds, just note that by
(\ref{V2}), (\ref{htrafo}), and (\ref{Vi2}),
\be\ba{l}
\dis P^{\Li(\Pois(h\mu))}[\Thin_f(X^h_n)=0]
=P[\Thin_{U^h_1\circ\cdots\circ U^h_nf}(\Pois(h\mu))=0]\\[5pt]
\dis\quad=P[\Pois((hU^h_1\circ\cdots\circ U^h_nf)\mu)=0]
=P[\Pois((\Ui_1\circ\cdots\circ\Ui_n(hf))\mu)=0]\\[5pt]
\dis\quad=P^\mu[\Pois(hf\Xc_n)=0]=P^\mu[\Thin_f(\Pois(h\Xc_n))=0].
\ec
Since this formula holds for all $f\in B_{[0,1]}(E^h)$, formula
(\ref{Poissonization}) follows.\qed

\brm{\bf(Boundedness of $h$)}
Propositions~\ref{weightprop} and \ref{Poisprop} generalize to the
case that $h$ is unbounded, except that in this case the cluster
mechanism in (\ref{Qih}) and the offspring mechanism in (\ref{Zxdef})
need in general not be continuous. Here, in order for (\ref{htrafo1})
and (\ref{htrafo}) to be well-defined, one needs to extend the
definition of $\Ui_kf$ to unbounded functions $f$, but this can always
be done unambiguously \cite[Lemma~9]{FSsup}.
\erm

\subsection{Sub- and superharmonic functions}

This section contains a number of pivotal calculations involving the
log-Laplace operators $\Ui_\ga$ from (\ref{Uiga}).  In particular, we
will prove that the functions $h_{1,1}$, $h_{0,0}$, and $h_{0,1}$ from
Lemmas~\ref{11lem}, \ref{00lem}, and \ref{01lem}, respectively, are
$\Ui_\ga$-superharmonic.

We start with an observation that holds for general log-Laplace operators.
\bl{\bf(Constant multiples)}\label{rhlem}
Let $\Ui$ be a log-Laplace operator of the form (\ref{Vk}) satisfying
(\ref{fincon}) and let $f\in B_+(E)$. Then $\Ui(rf)\leq r\Ui f$ for all
$r\geq 1$, and $\Ui(rf)\geq r\Ui f$ for all $0\leq r\leq 1$. In particular,
if $f$ is $\Ui$-superharmonic then $rf$ is $\Ui$-superharmonic for each
$r\geq 1$, and if $f$ is $\Ui$-subharmonic then $rf$ is $\Ui$-superharmonic
for each $0\leq r\leq 1$.
\el
{\bf Proof} If $\Xc$ is a branching process and $\Ui$ is the log-Laplace
operator of the transition law from $\Xc_0$ to $\Xc_1$ then, using Jensen's
inequality, for all $r\geq 1$,
\be\label{rhsup}
\ex{-\li\mu,\Ui(rf)\re}=E^\mu\big[\ex{-\li\Xc_1,rf\re}\big]
=E^\mu\big[\big(\ex{-\li\Xc_1,f\re}\big)^r\big]
\geq\big(E^\mu\big[\ex{-\li\Xc_1,f\re}\big]\big)^r=\ex{-\li\mu,r\Ui f\re}.
\ee
Since this holds for all $\mu\in\Mi(E)$, it follows that
$\Ui(rf)\leq r\Ui f$. The proof of the statements for $0\leq r\leq 1$
is the same but with the inequality signs reversed.\qed

\noi
We next turn our attention to the functions $h_{1,1}$ and $h_{0,0}$.
\bl{\bf(The catalyzing function $h_{1,1}$)}\label{sup11}
One has
\be\label{cofo}
\Ui_\ga(rh_{1,1})(x)=\frac{1+\ga}{\frac{1}{r}+\ga}\qquad(\ga,r>0,\ x\in[0,1]).
\ee
In particular, $h_{1,1}$ is $\Ui_\ga$-harmonic for each $\ga>0$.
\el
{\bf Proof} Recall (\ref{Zidef})--(\ref{Uiga}). Let $\sig_{1/r}$ be an
exponentially distributed random variable with mean $1/r$, independent
of $\tau_\ga$. Then
\be
\Ui_\ga(rh_{1,1})(x)=(\ffrac{1}{\ga}+1)
E\big[1-\ex{-\int_0^{\tau_\ga}r\di t}\big]
=(\ffrac{1}{\ga}+1)P[\sig_{1/r}<\tau_\ga]
=(\ffrac{1}{\ga}+1)\frac{\ga}{\frac{1}{r}+\ga},
\ee
which yields (\ref{cofo}).\qed

\bl{\bf(The catalyzing function $h_{0,0}$)}\label{sup00}
One has $\Ui_\ga(rh_{0,0})\leq rh_{0,0}$ for each $\ga,r>0$.
\el
{\bf Proof} Let $\Ga^\ga_x$ be the invariant law from Corollary~\ref{ergo}.
Then, for any $\ga>0$ and $f\in B_+[0,1]$,
\bc\label{linest}
\Ui_\ga f(x)&=&\dis(\ffrac{1}{\ga}+1)E\big[1-\ex{-\li\Zi^\ga_x,f\re}\big]
\leq(\ffrac{1}{\ga}+1)E[\li\Zi^\ga_x,f\re]\\[5pt]
&=&\dis(\ffrac{1}{\ga}+1)E\big[\int_0^{\tau_\ga}\!\!\!
f(\y^\ga_x(-t/2))\,\di t\big]=(1+\ga)\li\Ga^\ga_x,f\re\qquad(x\in[0,1]),
\ec
where we have used that $\tau_\ga$ is independent of $\y^\ga_x$ and has
mean $\ga$. In particular, setting $f=rh_{0,0}$ and using (\ref{WFfix})
we find that $\Ui_\ga(rh_{0,0})\leq rh_{0,0}$.\qed

\noi
The aim of the remainder of this section is to derive various bounds on
$\Ui_\ga f$ for $f\in\Hi_{0,1}$. We start with a formula for $\Ui_\ga f$
that holds for general $[0,1]$-valued functions $f$.
\bl{\bf(Action of $\Ui_\ga$ on $[0,1]$-valued functions)}\label{Ui01}
 Let $\y^\ga_x$ be the stationary solution to (\ref{Yclx}) and let
$\tau_{\ga/2}$ be an independent exponentially distributed random
variable with mean $\ga/2$. Let $(\bet_i)_{i\geq 1}$ be independent
exponentially distributed random variables with mean $\frac{1}{2}$,
independent of $\y^\ga_x$ and $\tau_{\ga/2}$, and let
$\sig_k:=\sum_{i=1}^k\bet_i$ $(k\geq 0)$. Then
\be\label{U1}
1-\Ui_\ga f(x)=E\Big[\prod_{k\geq 0:\;\sig_k<\tau_\ga}
\big(1-f(\y^\ga_x(-\sig_k))\big)\Big]
\qquad(\ga>0,\ f\in B_{[0,1]}[0,1],\ x\in[0,1]).
\ee 
\el
{\bf Proof} By Lemma~\ref{sup11}, the constant function
$h_{1,1}(x):=1$ satisfies $\Ui_\ga h_{1,1}=h_{1,1}$ for all $\ga>0$.
Therefore, by Proposition~\ref{Poisprop}, Poissonizing the
Poisson-cluster branching process $\Xc$ with the density $h_{1,1}$
yields a branching particle system
$X^{h_{1,1}}=(X^{h_{1,1}}_{-n},\ldots,X^{h_{1,1}}_0)$ with generating
operators $U^{h_{1,1}}_{\ga_{n-1}},\ldots,U^{h_{1,1}}_{\ga_0}$, where
\be
U^{h_{1,1}}_\ga f=\Ui_\ga f\qquad(f\in B_{[0,1]}[0,1],\ \ga>0).
\ee
By (\ref{VV}) and (\ref{Zxfo}),
\be
U^{h_{1,1}}_\ga f(x)=1-E\big[(1-f)^{\txt\Pois(\Zi^\ga_x)}\,\big|
\,\Pois(\Zi^\ga_x)\neq 0\big]\quad(f\in B_{[0,1]}[0,1],\ x\in[0,1],\ \ga>0).
\ee
Therefore, (\ref{U1}) will follow provided that
\be\label{conpos}
P\big[\Pois(\Zi^\ga_x)\in\cdot\,\big|\,\Pois(\Zi^\ga_x)\neq 0\big]
=\Li\Big(\sum_{k\geq 0:\;\sig_k<\tau_{\ga/2}}\de_{\y^\ga_x(-\sig_k)}\Big).
\ee
Indeed, it is not hard to see that
\be\label{unconpos}
\Pois(\Zi^\ga_x)\isd\sum_{k>0:\;\sig_k<\tau_{\ga/2}}\de_{\y^\ga_x(-\sig_k)}.
\ee
This follows from the facts that $\Zi^\ga_x=2\int_0^{\tau_{\ga/2}}
\de_{\y^\ga_x(-s)}\di s$ and
\be
\sum_{k>0:\;\sig_k<\tau_{\ga/2}}\de_{-\sig_k}
\isd\Pois(2\,1_{(-\tau_{\ga/2},0]}).
\ee
Conditioning $\Pois(2\,1_{(-\tau_{\ga/2},0]})$ on being nonzero means
conditioning on $\tau_{\ga/2}>\sig_1$. Since $\tau_{\ga/2}-\sig_1$,
conditioned on being nonnegative, is exponentially distributed with
mean $\ga/2$, using the stationarity of $\y^\ga_x$, we arrive at
(\ref{conpos}).\qed

\noi
The next lemma generalizes the duality (\ref{WFdual}) to mixed moments
of the Wright-Fisher diffusion $\y$ at multiple times. We can
interpret the left-hand side of (\ref{mixdua}) as the probability that
$m_1,\ldots,m_n$ organisms sampled from the population at times
$t_1,\ldots,t_n$ are all of the genetic type~I.
\bl{\bf(Sampling at multiple times)}\label{multidual}
Fix $0\leq t_1<\cdots<t_n=t$ and nonnegative integers $m_1,\ldots,m_n$.
Let $\y$ be the diffusion in (\ref{WFsde}). Then
\be\label{mixdua}
E^y\Big[\prod_{k=1}^n\y_{t_k}^{m_k}\Big]=E\big[y^{\phi_t}x^{\psi_t}\big],
\ee
where $(\phi_s,\psi_s)_{s\in[0,t]}$ is a Markov process in $\N^2$ started
in $(\phi_0,\psi_0)=(m_n,0)$, that jumps deterministically as
\be
(\phi_s,\psi_s)\to(\phi_s+m_k,\psi_s)\quad\mbox{at time}\quad t-t_k\quad(k<n),
\ee
and between these deterministic times jumps with rates as in (\ref{phipsi}).
\el
{\bf Proof} Induction, with repeated application of (\ref{WFdual}).\qed

\noi
For any $m\geq 1$, we put
\be
h_m(x):=1-(1-x)^m\qquad(x\in[0,1]).
\ee
The next lemma shows that we have particular good control on the action
of $\Ui_\ga$ on the functions $h_m$.
\bl{\bf(Action of $\Ui_\ga$ on the functions $h_m$)}\label{L:Uihm}
Let $m\geq 1$ and let $\tau_\ga$ be an exponentially distributed random
variable with mean $\ga$. Conditional on $\tau_\ga$, let
$(\phi'_t,\psi'_t)_{t\geq 0}$ be a Markov process in $\N^2$, started in
$(\phi'_0,\psi'_0)=(m,0)$ that jumps at time $t$ as:
\be\ba{r@{\,}c@{\,}l@{\qquad}l}\label{1phipsi}
(\phi'_t,\psi'_t)&\to&(\phi'_t-1,\psi'_t)
&\mbox{with rate}\ \phi'_t(\phi'_t-1),\\
(\phi'_t,\psi'_t)&\to&(\phi'_t-1,\psi'_t+1)
&\mbox{with rate}\ \ffrac{1}{\ga}\phi'_t,\\
(\phi'_t,\psi'_t)&\to&(\phi'_t+m,\psi'_t)
&\mbox{with rate}\ 1_{\{\tau_{\ga/2}<t\}}.
\ec
Then the limit $\lim_{t\to\infty}\psi'_t=:\psi'_\infty$ exists a.s., and
\be\label{suphest}
\Ui_\ga h_m(x)=E^{(m,0)}\big[1-(1-x)^{\psi'_\infty}\big]
\qquad(m\geq 1,\ x\in[0,1]).
\ee
\el
{\bf Proof} Let $\y^\ga_x$, $\tau_{\ga/2}$, and $(\sig_k)_{k\geq 0}$
be as in Lemma~\ref{Ui01}. Then, by (\ref{U1}),
\be\label{preUh01}
\Ui_\ga h_m(x)=1-E\Big[\prod_{k\geq 0:\;\sig_k<\tau_{\ga/2}}
\big(1-\y^\ga_x(-\sig_k)\big)^m\Big].
\ee
Let $(\phi',\psi')=(\phi'_t,\psi'_t)_{t\geq 0}$ be a $\N^2$-valued
process started in $(\phi'_0,\psi'_0)=(m,0)$ such that conditioned on
$\tau_\ga$ and $(\sig_k)_{k\geq 0}$, $(\phi',\psi')$ is a Markov
process that jumps deterministically as
\be
(\phi'_t,\psi'_t)\to(\phi'_t+m,\psi'_s)\quad\mbox{at time}
\quad\sig_k\quad(k\geq 1:\;\sig_k<\tau_{\ga/2})
\ee
and between these times jumps with rates as in (\ref{phipsi}). Then
$(\phi'_t,\psi'_t)\to(0,\psi'_\infty)$ as $t\to\infty$ a.s.\ for some
$\N$-valued random variable $\psi'_\infty$, and (\ref{suphest})
follows from Lemma~\ref{multidual}, using the symmetry
$y\leftrightarrow 1-y$. Since $\sig_{k+1}-\sig_k$ are independent
exponentially distributed random variables with mean one,
$(\phi',\psi')$ is the Markov process with jump rates as in
(\ref{1phipsi}).\qed

\noi
The next result is a simple application of Lemma~\ref{L:Uihm}.
\bl{\bf(The catalyzing function $h_1$)}\label{subx}
The function $h_1(x):=x$ $(x\in[0,1])$ is $\Ui_\ga$-subharmonic
for each $\ga>0$.
\el
{\bf Proof} Since $\psi'_\infty\geq 1$ a.s., one has
$1-(1-x)^{\psi'_\infty}\geq x$ a.s.\ $(x\in[0,1])$ in
(\ref{suphest}). In particular, setting $m=1$ yields
$\Ui_\ga h_1\geq h_1$.\qed

\noi
We now set out to prove that $h_7$, which is the function $h_{0,1}$
from Lemma~\ref{01lem}, is $\Ui_\ga$-super\-harmonic. In order to do
so, we will derive upper bounds on the expectation of $\psi'_\infty$.
We derive two estimates: one that is good for small $\ga$ and one that
is good for large $\ga$.

In order to avoid tedious formal arguments, it will be convenient to
recall the interpretation of the process $(\phi',\psi')$ and
Lemma~\ref{multidual}. Recall from the discussion following
(\ref{WFdual}) that $(\y^\ga_x(t))_{t\in\R}$ describes the equilibrium
frequency of genetic type $I$ as a function of time in a population
that is in genetic exchange with an infinite reservoir. From this
population we sample at times $-\sig_k$ ($k\geq 0$,
$\sig_k<\tau_{\ga/2}$) each time $m$ individuals, and ask for the
probability that they are not all of the genetic type~II. In order to
find this probability, we follow the ancestors of the sampled
individuals back in time. Then $\phi'_t$ and $\psi'_t$ are the number
of ancestors that lived at time $-t$ in the population and the
reservoir, respectively, and $E[1-(1-x)^{\psi'_\infty}]$ is the
probability that at least one ancestor is of type~I.
\bl{\bf(Bound for small $\ga$)}\label{smaga}
For each $\ga\in(0,\infty)$ and $m\geq 1$,
\be\label{smag}
\frac{1}{m}E^{(m,0)}[\psi'_\infty]\leq\frac{1}{m}
\sum_{i=0}^{m-1}\frac{1+\ga}{1+i\ga}=:\chi_m(\ga).
\ee
The function $\chi_m$ is concave and satisfies $\chi_m(0)=1$
for each $m\geq 1$.
\el
{\bf Proof} Note that
\be
E\big[\big|\{k\geq 0:\;\sig_k<\tau_{\ga/2}\}\big|\big]=1+\ga.
\ee
We can estimate $(\phi',\psi')$ from above by a process where
ancestors from individuals sampled at different times cannot
coalesce. Therefore,
\be\label{aces}
E^{(m,0)}[\psi'_\infty]\leq(1+\ga)E^{(m,0)}[\psi_\infty],
\ee
where $(\phi,\psi)$ is the Markov process in (\ref{phipsi}).
Note that if $(\phi,\psi)$ is in the state $(m+1,0)$, then the
next jump is to $(m,1)$ with probability
\be
\frac{\frac{1}{\ga}(m+1)}{\frac{1}{\ga}(m+1)+m(m+1)}=\frac{1}{1+m\ga}
\ee
and to $(m,0)$ with one minus this probability. Therefore,
\bc
\dis E^{(m+1,0)}[\psi_\infty]&=&\dis\frac{1}{1+m\ga}
E^{(m,1)}[\psi_\infty]+\Big(1-\frac{1}{1+m\ga}\Big)
E^{(m,0)}[\psi_\infty]\\[5pt]
&=&\dis\frac{1}{1+m\ga}\Big(E^{(m,0)}[\psi_\infty]+1\Big)
+\Big(1-\frac{1}{1+m\ga}\Big)E^{(m,0)}[\psi_\infty]\\[5pt]
&=&\dis E^{(m,0)}[\psi_\infty]+\frac{1}{1+m\ga}.
\ec
By induction, it follows that
\be\label{psim}
E^{(m,0)}[\psi_\infty]=\sum_{i=0}^{m-1}\frac{1}{1+i\ga}.
\ee
Inserting this into (\ref{aces}) we arrive at (\ref{smag}). Finally, since
\be
\Diff{\ga}\frac{1+\ga}{1+i\ga}=\frac{2i(i-1)}{(1+i\ga)^3}\geq 0
\qquad(i\geq 0,\ \ga\geq 0),
\ee
the function $\chi_m$ is convex.\qed

\bl{\bf(Bound for large $\ga$)}\label{laga}
For each $\ga\in(0,\infty)$ and $m\geq 1$,
\be\label{lag}
E^{(m,0)}[\psi'_\infty]\leq(\ffrac{1}{\ga}+1)\sum_{k=1}^m\frac{1}{k}
+\frac{3}{2}.
\ee
\el
{\bf Proof} We start by observing that $\dif{t}E[\psi'_t]
=\frac{1}{\ga}E[\phi'_t]$, and therefore
\be\label{Epsi}
E[\psi'_\infty]=\ffrac{1}{\ga}\int_0^\infty E[\phi'_t]\di t.
\ee
Unlike in the proof of the last lemma, this time we cannot fully
ignore the coalescence of ancestors sampled at different times. In
order to deal with this we use a trick: at time zero we introduce an
extra ancestor that can only jump to the reservoir when
$t\geq\tau_\ga$ and there are no other ancestors left in the
population. We further assume that all other ancestors do not jump to
the reservoir on their own. Let $\xi_t$ be one as long as this extra
ancestor is in the population and zero otherwise, and let $\phi''_t$
be the number of other ancestors in the population according to these
new rules. Then we have at a Markov process $(\xi,\phi'')$ started in
$(\xi_0,\phi''_0)=(1,m)$ that jumps as:
\be\ba{r@{\,}c@{\,}l@{\qquad}l}\label{2phipsi}
(\xi_t,\phi''_t)&\to&(\xi_t,\phi''_t-1)
&\mbox{with rate}\ (\phi''_t+1)\phi''_t,\\
(\xi_t,\phi''_t)&\to&(\xi_t,\phi''_t+m)
&\mbox{with rate}\ 1_{\{\tau_{\ga/2}<t\}},\\
(\xi_t,\phi''_t)&\to&(\xi_t-1,\phi''_t)
&\mbox{with rate}\ \frac{1}{\ga}1_{\{\tau_{\ga/2}\geq t\}}1_{\{\phi''_t=0\}}.
\ec
It is not hard to show that $(\xi,\phi'')$ and $\phi'$ can be coupled
such that $\xi_t+\phi''_t\geq\phi'_t$ for all $t\geq 0$. We now
simplify even further and ignore all coalescence between ancestors
belonging to the process $\phi''$ that are introduced at different
times. Let $\phi^{(k)}_t$ be the number of ancestors in the population
that were introduced at the time $\sig_k$ $(k\geq 0)$. Thus, for
$t<\sig_k$ one has $\phi^{(k)}_t=0$, for $t=\sig_k$ one has
$\phi^{(k)}_t=m$, while for $t>\sig_k$, the process $\phi^{(k)}_t$
jumps from $n$ to $n-1$ with rate $(n+1)n$. Then it is not hard to see
that, for an appropriate coupling, $\phi''_t\leq\sum_{k\geq 0
:\sig_k<\tau_{\ga/2}}\phi^{(k)}_t$ for all $t\geq 0$. We let $\xi'$
be a process such that $\xi'_0=1$ and $\xi'_t$ jumps to zero with rate
\be
\frac{1}{\ga}1_{\{\tau_{\ga/2}\geq t\}}
\prod_{k\geq 0:\sig_k<\tau_{\ga/2}}1_{\{\phi^{(k)}_t=0\}}.
\ee
Then for an appropriate coupling $\xi'_t\geq\xi_t$ $(t\geq 0)$.
Thus, we can estimate
\be\label{phint}
\int_0^\infty E[\phi'_t]\di t\leq\int_0^\infty E[\xi'_t]\di t
+\int_0^\infty E\Big[\sum_{k\geq 0:\sig_k<\tau_{\ga/2}}\phi^{(k)}_t\Big]\di t.
\ee
Set $\rho:=\inf\{t\geq\tau_{\ga/2}:\phi^{(k)}_t=0\ \forall k\geq 0
\mbox{ with }\sig_k<\tau_{\ga/2}\}$ and $\pi:=\inf\{t\geq 0:\xi'_t=0\}$.
Then
\be\label{xint}
\int_0^\infty E[\xi'_t]\di t=E[\tau_{\ga/2}]+E[\rho-\tau_{\ga/2}]
+E[\pi-\rho]=\frac{3}{2}\ga+E[\rho-\tau_{\ga/2}].
\ee
Since
\bc
\dis E[\rho-\tau_{\ga/2}]&\leq&\dis\int_0^\infty
E\Big[1_{\{\sum_{k\geq 0:\sig_k<\tau_{\ga/2}}
\phi^{(k)}_t\neq 0\}}\Big]\di t\\[10pt]
&\leq&\dis\int_0^\infty E\Big[\sum_{k\geq 0:\sig_k<\tau_{\ga/2}}
1_{\{\phi^{(k)}_t\neq 0\}}\Big]\di t,
\ec
using moreover (\ref{phint}) and (\ref{xint}), we can estimate
\be
\int_0^\infty E[\phi'_t]\di t\leq\frac{3}{2}\ga
+\int_0^\infty E\Big[\sum_{k\geq 0:\sig_k<\tau_{\ga/2}}
(\phi^{(k)}_t+1_{\{\phi^{(k)}_t\neq 0\}})\Big]\di t.
\ee
Since $E\big[\big|\{k\geq 0:\;\sig_k<\tau_{\ga/2}\}\big|\big]=1+\ga$, we obtain
\be\label{obta}
\int_0^\infty E[\phi'_t]\di t\leq\frac{3}{2}\ga
+(1+\ga)\int_0^\infty E[\phi^{(0)}_t+1_{\{\phi^{(0)}_t\neq 0\}}]\di t.
\ee
Since $\phi^{(0)}_t$ jumps from $n$ to $n-1$ with rate $(n+1)n$, the
expected total time
that $\phi^{(0)}_t=n$ equals $1/((n+1)n)$, and therefore
\be
\int_0^\infty E[\phi^{(0)}_t+1_{\{\phi^{(0)}_t\neq 0\}}]\di t
=\sum_{n=1}^m\frac{1}{(n+1)n}(n+1_{\{n\neq 0\}})
=\sum_{n=1}^m\frac{1}{n}.
\ee
Inserting this into (\ref{obta}), using (\ref{Epsi}), we arrive
at (\ref{lag}).\qed

\bl{\bf(The catalyzing function $h_{0,1}$)}\label{sup01}
One has $\Ui_\ga(h_{0,1})\leq h_{0,1}$ for each $\ga>0$. Moreover,
for each $r>1$ and $\ga>0$,
\be\label{strongsup}
\sup_{x\in(0,1]}\frac{\Ui_\ga(rh_{0,1})(x)}{rh_{0,1}(x)}<1.
\ee
\el
{\bf Proof} Recall that $h_{0,1}(x)=h_7(x)=1-(1-x)^7$ $(x\in[0,1])$.
We will show that
\be
E^{(7,0)}[\psi'_\infty]<7
\ee
for each $\ga\in(0,\infty)$. The function $\chi_m(\ga)$ from
Lemma~\ref{smaga} satisfies
\be
\chi_m(1)=\frac{1}{m}\sum_{n=1}^m\frac{2}{n}<1\qquad(m\geq 5).
\ee
Since $\chi_m(\ga)$ is concave in $\ga$ and satisfies $\chi_m(0)=1$,
it follows that $\chi_m(\ga)<1$ for all $0<\ga\leq 1$ and $m\geq 5$.
By Lemma~\ref{laga}, for all $\ga\geq 1$,
\be
E^{(m,0)}[\psi'_\infty]\leq 2\sum_{k=1}^m\frac{1}{k}+\frac{3}{2}<m
\qquad(m\geq 7).
\ee
Therefore, if $m\geq 7$, then $m':=E^{(m,0)}[\psi'_\infty]<m$. It
follows by (\ref{suphest}) and Jensen's inequality applied to the
concave function $z\mapsto 1-(1-x)^z$ that
\be\label{hmjens}
\Ui_\ga h_m(x)\leq 1-(1-x)^{E^{(m,0)}[\psi'_\infty]}=1-(1-x)^{m'}
\leq h_m(x)\qquad(x\in[0,1],\ \ga>0).
\ee
This shows that $h_m$ is $\Ui_\ga$-superharmonic for each $\ga>0$.
By Lemma~\ref{rhlem}, for each $r>1$,
\be\label{leftb}
\frac{\Ui_\ga(rh_m)(x)}{rh_m(x)}\leq\frac{r\Ui_\ga(h_m)(x)}{rh_m(x)}
\leq\frac{1-(1-x)^{m'}}{1-(1-x)^m}\qquad(x\in(0,1]).
\ee
By Lemma~\ref{sup11} and the monotonicity of $\Ui_\ga$,
\be\label{rightb}
\frac{\Ui_\ga(rh_m)(x)}{rh_m(x)}\leq\frac{\Ui_\ga(r)(x)}{rh_m(x)}
\leq\frac{1+\ga}{1+r\ga}\frac{1}{(1-(1-x)^m)}\qquad(x\in(0,1]).
\ee
Since the right-hand side of (\ref{leftb}) is smaller than $1$ for
$x\in(0,1)$ and tends to $m'/m<1$ as $x\to 0$, since the right-hand
side of (\ref{rightb}) is smaller than $1$ for $x$ in an open
neighborhood of $1$, and since both bounds are continuous,
(\ref{strongsup}) follows.\qed

\subsection{Extinction versus unbounded growth}\label{versus}

In this section we show that Lemmas~\ref{11lem}--\ref{01lem} are
equivalent to Proposition~\ref{P:exgr}. (This follows from the
equivalence of conditions (i) and (ii) in Lemma~\ref{exgro} below.) We
moreover prove Lemmas~\ref{11lem} and \ref{01lem} and prepare for the
proof of Lemma~\ref{00lem}. We start with some general facts about
log-Laplace operators and branching processes.

For the next lemma, let $E$ be a separable, locally compact,
metrizable space. For $n\geq 0$, let $q_n\in\Ci_+(E)$ be continuous
weight functions, let $\Qi_n$ be continuous cluster mechanisms on $E$,
and assume that the associated log-Laplace operators $\Ui_n$ defined
in (\ref{Vk}) satisfy (\ref{fincon}). Assume that $0\neq h\in\Ci_+(E)$
is bounded and $\Ui_n$-superharmonic for all $n$, let $E^h:=\{x\in
E:h(x)>0\}$, and define generating operators $U^h_n:B_{[0,1]}(E^h)\to
B_{[0,1]}(E)$ as in (\ref{htrafo}). For each $n\geq 0$, let
$(\Xc^{(n)}_0,\Xc^{(n)}_1)$ be a one-step Poisson cluster branching
process with log-Laplace operator $\Ui_n$, and let
$(X^{(n),h}_0,X^{(n),h}_1)$ be the one-step branching particle system
with generating operator $U^h_n$. (In a typical application of this
lemma, the operators $\Ui_n$ will be iterates of other log-Laplace
operators, and $\Xc^{(n)}_0,\Xc^{(n)}_1$ will be the initial and final
state, respectively, of a Poisson cluster branching process with many
time steps.)
\bl{\bf(Extinction versus unbounded growth)}\label{exgro}
Assume that $\rho\in\Ci_{[0,1]}(E^h)$ and put
\be
p(x):=\left\{\ba{ll}h(x)\rho(x)\quad&\mbox{if}\ x\in E^h,\\
0\quad&\mbox{if}\ x\in E\beh E^h.\ea\right.
\ee
Then the following statements are equivalent:
\[\ba{rl}
{\rm (i)}&\dis P^{\de_x}\big[|X^{(n),h}_1|\in\cdot\,\big]\Asto{n}
\rho(x)\de_\infty+(1-\rho(x))\de_0\\
&\hspace{1.5cm}\mbox{locally uniformly for }x\in E^h,\\[5pt]
{\rm (ii)}&\dis P^{\de_x}\big[\li\Xc^{(n)}_1,h\re\in\cdot\,\big]\Asto{n}
\ex{-p(x)}\de_0+\big(1-\ex{-p(x)}\big)\de_\infty\\
&\hspace{1.5cm}\mbox{locally uniformly for }x\in E,\\[5pt]
{\rm (iii)}&\Ui_n(\la h)(x)\asto{n}p(x)\\
&\hspace{1.5cm}\mbox{locally uniformly for }x\in E\quad\forall\la>0,\\[5pt]
{\rm (iv)}&\exists0<\la_1<\la_2<\infty:\quad\Ui_n(\la_i h)(x)\asto{n}p(x)\\
&\hspace{1.5cm}\mbox{locally uniformly for }x\in E\qquad(i=1,2).
\ea\]
\el
{\bf Proof of Lemma~\ref{exgro}} It is not hard to see that (i) is
equivalent to
\be
P^{\de_x}[\Thin_\la(X^{(n),h}_1)\neq 0]\asto{n}\rho(x)
\ee
locally uniformly for $x\in E^h$, for all $0<\la\leq 1$. It follows from
(\ref{V2}) and (\ref{htrafo}) that $h(x)P^{\de_x}[\Thin_\la(X^{(n),h}_1)\neq 0]
=hU^h(\la)(x)=\Ui(\la h)(x)$ $(x\in E)$, so (i) is equivalent to
\[\ba{l}
{\rm (i)'}\quad\Ui_n(\la h)(x)\asto{n}p(x)\\
\hspace{1.5cm}\mbox{locally uniformly for }x\in E\quad\forall0<\la\leq 1.
\ea\]
By (\ref{Vi}), condition~(ii) implies that
\be
\ex{-\Ui_n(\la h)(x)}=E^{\de_x}\big[\ex{-\la\li\Xc_1,h\re}\big]\asto{n}
\ex{-p(x)}
\ee
locally uniformly for $x\in E$ for all $\la>0$, and therefore (ii)
implies (iii). Obviously (iii)$\volgt{\rm (i)'}\volgt$(iv) so we are
done if we show that (iv)$\volgt$(ii). Indeed, (iv) implies that
\be
E^{\de_x}\big[\ex{-\la_1\li\Xc^{(n)}_1,h\re}
-\ex{-\la_2\li\Xc^{(n)}_1,h\re}\big]\asto{n}0
\ee
locally uniformly for $x\in E$, which shows that
\be
P^{\de_x}\big[c<\li\Xc^{(n)}_1,h\re<C\big]\asto{n}0
\ee
for all $0<c<C<\infty$. Using (iv) once more we arive at (ii).\qed

\noi
Our next lemma gives sufficient conditions for the $n$-th iterates of
a single log-Laplace operator $\Ui$ to satisfy the equivalent
conditions of Lemma~\ref{exgro}. Let $E$ (again) be separable, locally
compact, and metrizable. Let $q\in\Ci_+(E)$ be a weight function,
$\Qi$ a continuous cluster mechanism on $E$, and assume that the
associated log-Laplace operator $\Ui$ defined in (\ref{Vk}) satisfies
(\ref{fincon}). Let $\Xc=(\Xc_0,\Xc_1,\ldots)$ be the Poisson-cluster
branching process with log-Laplace operator $\Ui$ in each step, let
$0\neq h\in\Ci_+(E)$ be bounded and $\Ui$-superharmonic, and let
$X^h=(X^h_0,X^h_1,\ldots)$ denote the branching particle system on
$E^h$ obtained from $\Xc$ by Poissonization with a $\Ui$-superharmonic
function $h$, in the sense of Proposition~\ref{Poisprop}.
\bl{\bf(Sufficient condition for extinction versus unbounded
growth)}\label{L:sufexgro}
Assume that
\be\label{strsup}
\sup_{x\in E^h}\frac{\Ui h(x)}{h(x)}<1.
\ee
Then the process $X^h$ started in any initial law $\Li(X^h_0)\in\Mi_1(E^h)$
satisfies
\be\label{extgrow}
\lim_{k\to\infty}|X^h_k|=\infty\quad\mbox{or}\quad\exists k\geq 0
\mbox{ s.t.\ }X^h_k=0\qquad\as
\ee
Moreover, if the function $\rho:E^h\to[0,1]$ defined by
\be
\rho(x):=P^{\de_x}[X^h_n\neq 0\quad\forall n\geq 0]\qquad(x\in E^h)
\ee
satisfies $\inf_{x\in E^h}\rho(x)>0$, then $\rho$ is continuous.
\el
{\bf Proof of Lemma~\ref{L:sufexgro}} Let $\Ai$ denote the tail event
$\Ai=\{X^h_n\neq 0\ \forall n\geq 0\}$ and let $(\Fi_k)_{k\geq 0}$ be
the filtration generated by $X^h$. Then, by the Markov property and
continuity of the conditional expectation with respect to increasing
limits of \si-fields (see Complement 10(b) from
\cite[Section~29]{Loe63} or \cite[Section~32]{Loe78})
\be
P[X^h_n\neq 0\ \forall n\geq 0|X_k]=P(\Ai|\Fi_k)\asto{k}1_\Ai\qquad\as
\ee
In particular, this implies that a.s.\ on the event $\Ai$ one must
have $P[X^h_{k+1}=0|X^h_k]\to 0$ a.s. By (\ref{V2}) and
(\ref{htrafo}), $P^{\de_x}[X^h_1\neq 0]=U^h 1(x)=(\Ui h(x))/h(x)$,
which is uniformly bounded away from one by (\ref{strsup}). Therefore,
$P[X^h_{k+1}=0|X^h_k]\to 0$ a.s.\ on $\Ai$ is only possible if the
number of particles tends to infinity.

The continuity of $\rho$ can be proved by a straightforward adaptation
of the proof of \cite[Proposition~5 (d)]{FStrim} to the present
setting with discrete time and noncompact space $E$. An essential
ingredient in the proof, apart from (\ref{strsup}), is the fact that
the map $\nu\mapsto P^\nu[X^h_n\in\cdot\,]$ from $\Ni(E)$ to
$\Mi_1(\Ni(E))$ is continuous, which follows from the continuity of
$Q^h$.\qed

\noi
We now turn our attention more specifically to the renormalization
branching process $\Xc$. In the remainder of this section,
$(\ga_k)_{k\geq 0}$ is a sequence of positive constants such that
$\sum_n\ga_n=\infty$ and $\ga_n\to\ga^\ast$ for some
$\ga^\ast\in\half$, and $\Xc=(\Xc_{-n},\ldots,\Xc_0)$ is the Poisson
cluster branching process on $[0,1]$ defined in Section~\ref{RBP}. We
put $\Ui^{(n)}:=\Ui_{\ga_{n-1}}\circ\cdots\circ\Ui_{\ga_0}$. If $0\neq
h\in\Ci[0,1]$ is $\Ui_{\ga_k}$-superharmonic for all $k\geq 0$, then
$\Xc^h$ and $X^h$ denote the branching process and the branching
particle system on $\{x\in[0,1]:h(x)>0\}$ obtained from $\Xc$ by
weighting and Poissonizing with $h$ in the sense of
Propositions~\ref{weightprop} and \ref{Poisprop}, respectively.\med

\noi
{\bf Proof of Lemma~\ref{11lem}} By induction, it follows from
Lemma~\ref{sup11} that
\be
\Ui^{(n)}(\la h_{1,1})=\frac{\prod_{k=0}^{n-1}(1+\ga_k)}
{\prod_{k=0}^{n-1}(1+\ga_k)-1+\frac{1}{\la}}\qquad(\la>0).
\ee
It is not hard to see (compare the footnote at (\ref{sumga})) that
\be
\prod_{k=0}^\infty(1+\ga_k)=\infty\quad\mbox{if and only if}
\quad\sum_{k=0}^\infty\ga_k=\infty.
\ee
Therefore, since we are assuming that $\sum_n\ga_n=\infty$,
\be
\Ui^{(n)}(\la h_{1,1})\asto{n}h_{1,1},
\ee
uniformly on $[0,1]$ for all $\la>0$. The result now follows from
Lemma~\ref{exgro} (with $h=h_{1,1}$ and $\rho(x)=1$ $(x\in[0,1])$).\qed

\brm{\bf(Conditions on $(\ga_n)_{n\geq 0}$)}
Our proof of Lemma~\ref{11lem} does not use that $\ga_n\to\ga^\ast$
for some $\ga^\ast\in\half$. On the other hand, the proof shows that
$\sum_n\ga_n=\infty$ is a necessary condition for (\ref{expl}).
\erm
We do not know if the assumption that $\ga_n\to\ga^\ast$ for some
$\ga^\ast\in\half$ is needed in Lemma~\ref{00lem}. We guess that
it can be dropped, but it will greatly simplify proofs to have it around.

We will show that in order to prove Lemmas~\ref{00lem} and
\ref{01lem}, it suffices to prove their analogues for embedded
particle systems in the time-homogeneous processes $\Yi^{\ga^\ast}$
($\ga^\ast\in\half$). More precisely, we will derive
Lemmas~\ref{00lem} and \ref{01lem} from the following two results.
Below, $(\Ui^0_t)_{t\geq 0}$ is the log-Laplace semigroup of the
super-Wright-Fisher diffusion $\Yi^0$, defined in (\ref{cau}). The
functions $p^\ast_{0,1,\ga^\ast}$ ($\ga^\ast\in\half$) are defined in
(\ref{pdef}).
\bp{\bf(Time-homogeneous embedded particle system with
$h_{0,0}$)}\label{00hom}\smallskip

\noi
{\bf (a)} For any $\ga^\ast>0$, one has $\dis(\Ui_{\ga^\ast})^n h_{0,0}
\asto{n}0$ uniformly on $[0,1]$.\smallskip

\noi
{\bf (b)} One has $\dis\Ui^0_th_{0,0}\asto{t}0$ uniformly on $[0,1]$.
\ep

\bp{\bf(Time-homogeneous embedded particle system with $h_{0,1}$)}\label{01hom}

\noi
{\bf (a)} For any $\ga^\ast>0$, one has
$(\Ui_{\ga^\ast})^n(\la h_{0,1})\asto{n}p^\ast_{0,1,\ga^\ast}$
uniformly on $[0,1]$, for all $\la>0$.\med

\noi
{\bf (b)} One has
$\Ui^0_t(\la h_{0,1})\asto{t}p^\ast_{0,1,0}$ uniformly on $[0,1]$,
for all $\la>0$.
\ep
Proposition~\ref{00hom}~(a) will be proved in Section~\ref{002sec}.\med

\noi
{\bf Proof of Proposition~\ref{01hom}~(a)} By formula (\ref{strongsup})
from Lemma~\ref{sup01}, for each $r>1$ the function $rh_{0,1}$ satisfies
condition (\ref{strsup}) from Lemma~\ref{L:sufexgro}. Set $\rho(x)
:=P^{\de_x}[Y^{\ga^\ast,rh_{0,1}}_n\neq 0\ \forall n]$. Then, by
(\ref{V2}) and (\ref{htrafo}),
\bc\label{rholow}
\dis\rho(x)&=&\dis\lim_{n\to\infty}P^{\de_x}[Y^{\ga^\ast,rh_{0,1}}_n\neq 0]
=\lim_{n\to\infty}(U^{rh_{0,1}}_{\ga^\ast})^n1(x)\\[5pt]
&=&\dis\lim_{n\to\infty}\frac{(\Ui_{\ga^\ast})^n(rh_{0,1})(x)}{rh_{0,1}(x)}
\geq\frac{h_1(x)}{rh_{0,1}(x)}\qquad(x\in(0,1]),
\ec
where $h_1(x)=x$ $(x\in[0,1])$ is the $\Ui_{\ga^\ast}$-subharmonic function
from Lemma~\ref{subx}. It follows that $\inf_{x\in(0,1]}\rho(x)>0$ and
therefore, by Lemma~\ref{L:sufexgro}, $\rho$ is continuous in $x$.

By Lemma~\ref{L:sufexgro}, we see that the Poissonized particle system
$X^{rh_{0,1}}$ exhibits extinction versus unbounded growth in the sense
of Lemma~\ref{exgro}, which implies the statement in
Proposition~\ref{01hom}~(a).\qed

\noi
{\bf Proof of Propositions~\ref{00hom}~(b) and \ref{01hom}~(b)}
These statements follow from results in \cite{FSsup}. Indeed,
\cite[Proposition~2]{FSsup} implies that for any $f\in B_+[0,1]$
and $x\in[0,1]$,
\be\ba{r@{\,}c@{\,}ll}\label{point}
\Ui^0_tf(x)&\asto{t}&0\quad&\mbox{if }f(0)=f(1)=0,\\[5pt]
\Ui^0_tf(x)&\asto{t}&p^\ast_{0,1,\ga^\ast}(x)\quad&\mbox{if }f(0)=0,\ f(1)>0.
\ec
To see that the convergence in (\ref{point}) is in fact uniform in
$x\in[0,1]$ we use the fact that each function $f\in B_+[0,1]$ with
$f(0)=f(1)=0$ can be bounded as $f\leq r1_{(0,1)}$ for some $r\geq 1$,
and that each function $f\in B_+[0,1]$ with $f(0)=0$ and $f(1)>0$ can
be bounded as $\eps 1_{\{1\}}\leq f\leq r1_{(0,1]}$ for some
$0<\eps\leq 1$ and $r\geq 1$. Therefore, by the monotonity of
$\Ui^0_t$, it suffices to show that $\Ui^0_t(r1_{(0,1)})$,
$\Ui^0_t(r1_{(0,1]})$, and $\Ui^0_t(\eps 1_{\{1\}})$ converge
uniformly on $[0,1]$. By \cite[Lemma~15]{FSsup}, these functions are
continuous for each $t>0$, and since moreover the limit functions are
continuous, it suffices to show that the convergence is monotone.
Thus, we claim that
\be\ba{r@{\,}c@{\,}ll}\label{monco}
\dis\Ui^0_t(r1_{(0,1)})&\down&\dis 0\qquad&(r\geq 1),\\[5pt]
\dis\Ui^0_t(r1_{(0,1]})&\down&\dis p^\ast_{0,1,\ga^\ast}
\qquad&(r\geq 1),\\[5pt]
\dis\Ui^0_t(\eps1_{\{1\}})&\up&\dis p^\ast_{0,1,\ga^\ast}
\qquad&(0<\eps\leq 1).
\ec
By (an obvious analogue of) Lemma~\ref{rhlem}, it suffices to show
that $1_{(0,1)}$ and $1_{(0,1]}$ are $\Ui^0_t$-superharmonic, while
$1_{\{1\}}$ is $\Ui^0_t$-subharmonic for each $t\geq 0$. Let
$(\Yi^{0,h_{1,1}}_t)_{t\geq 0}$ be the branching particle system
obtained from $(\Yi^0_t)_{t\geq 0}$ by Poissonization with the
constant function $h_{1,1}:=1$. Then $\Yi^{0,h_{1,1}}$ is a system of
binary splitting Wright-Fisher diffusions, which was also studied in
\cite{FSsup}. One has (compare (\ref{V2}))
\be
\Ui^0_t1_{(0,1)}(x)=P[\Thin_{\Ui^0_t1_{(0,1)}}(\de_x)\neq 0]
=P^{\de_x}[\Thin_{1_{(0,1)}}(Y^{0,h_{1,1}}_t)\neq 0]
=P^{\de_x}[Y^{0,h_{1,1}}_t((0,1))>0].
\ee
Likewise,
\be
\Ui^0_t1_{(0,1]}(x)=P^{\de_x}[Y^{0,h_{1,1}}_t((0,1])>0]\quad\mbox{and}
\quad\Ui^0_t1_{\{1\}}(x)=P^{\de_x}[Y^{0,h_{1,1}}_t(\{1\})>0].
\ee
Using the fact that the points $0,1$ are traps for the Wright-Fisher
diffusion and that in a binary splitting Wright-Fisher diffusion,
particles never die, it is easy to see that
$P^{\de_x}[Y^{0,h_{1,1}}_t((0,1))>0]$ and
$P^{\de_x}[Y^{0,h_{1,1}}_t((0,1])>0]$ are nonincreasing in $t$, while
$P^{\de_x}[Y^{0,h_{1,1}}_t(\{1\})>0]$ is nondecreasing in $t$.\qed

\noi
We now show that Propositions~\ref{00hom} and \ref{01hom} imply
Lemmas~\ref{00lem} and \ref{01lem}, respectively.\med

\noi
{\bf Proof of Lemma~\ref{00lem}} We start with the proof that the
embedded particle system $X^{h_{0,0}}$ is critical. For any
$f\in B_+[0,1]$ and $k\geq 1$, we have, by Poissonization
(Proposition~\ref{Poisprop}) and the definition of $\Xc$,
\be\ba{l}\label{critcalc}
\dis h_{0,0}(x)E^{-k,\de_x}[\li X^{h_{0,0}}_{-k+1},f\re]
=E^{-k,\Li(\Pois(h_{0,0}\de_x))}[\li X^{h_{0,0}}_{-k+1},f\re]
=E^{-k,\de_x}[\li\Pois(h_{0,0}\Xc_{-k+1}),f\re]\\[5pt]
\dis\qquad=E^{-k,\de_x}[\li\Xc_{-k+1},h_{0,0}f\re]
=(\ffrac{1}{\ga}+1)E[\li\Zi^\ga_x,h_{0,0}f\re]
=(\ffrac{1}{\ga}+1)\li\Ga^{\ga_{k-1}}_x,h_{0,0}f\re,
\ec
where $\Ga^\ga_x$ is the invariant law of $\y^\ga_x$ from
Corollary~\ref{ergo}. In particular, setting $f=1$ gives
$h_{0,0}(x)E^{-k,\de_x}[|X^{h_{0,0}}_{-k+1}|]=h_{0,0}(x)$ by (\ref{WFfix}).

To prove (\ref{ext}), by Lemma~\ref{exgro} it suffices to show that
\be\label{TB}
\Ui^{(n)}(\la h_{0,0})\asto{n}0
\ee
uniformly on $[0,1]$ for all $0<\la\leq 1$. We first treat the case
$\ga^\ast>0$. Then, by Theorem~\ref{supercon}~(a), for each fixed
$l\geq 1$ and $f\in\Ci_+[0,1]$,
\be\label{Uitohom}
\Ui_{\ga_{n-1}}\circ\cdots\circ\Ui_{\ga_{n-l}}f\asto{n}(\Ui_{\ga^\ast})^lf
\ee
uniformly on $[0,1]$. Therefore, by a diagonal argument, we can find
$l(n)\to\infty$ such that
\be
\big\|(\Ui_{\ga^\ast})^{l(n)}h_{0,0}-\Ui_{\ga_{n-1}}\circ\cdots\circ
\Ui_{\ga_{n-l(n)}}h_{0,0}\big\|_\infty\asto{n}0.
\ee
Using the fact that the function $h_{0,0}$ is $\Ui_\ga$-superharmonic
for each $\ga>0$ and the monotonicity of the operators $\Ui_\ga$, we
derive from Proposition~\ref{00hom}~(a) that
\be
\Ui^{(n)}(\la h_{0,0})\leq\Ui_{\ga_{n-1}}\circ\cdots\circ
\Ui_{\ga_{n-l(n)}}h_{0,0}\asto{n}0
\ee
uniformly on $[0,1]$ for all $0<\la\leq 1$. This proves (\ref{TB})
in the case $\ga^\ast>0$.

The proof in the case $\ga^\ast=0$ is similar. In this case, by
Theorem~\ref{supercon}~(b), for each fixed $t>0$ and $f\in\Ci_+[0,1]$,
\be\label{Uitohom2}
\Ui_{\ga_{n-1}}\circ\cdots\circ\Ui_{\ga_{k_n(t)}}f(x_n)\asto{n}
\Ui^0_tf(x)\qquad\forall x_n\to x\in[0,1],
\ee
which shows that $\Ui_{\ga_{n-1}}\circ\cdots\circ\Ui_{\ga_{k_n(t)}}f$
converges to $\Ui^0_tf$ uniformly on $[0,1]$. By a diagonal argument,
we can find $t(n)\to\infty$ such that
\be
\big\|\Ui^0_t(h_{0,0})-\Ui_{\ga_{n-1}}\circ\cdots\circ
\Ui_{\ga_{k_n(t(n))}}(h_{0,0})\big\|_\infty\asto{n}0,
\ee
and the proof proceeds in the same way as before.\qed

\noi
{\bf Proof of Lemma~\ref{01lem}} By Lemma~\ref{exgro} and the monotonicity
of the operators $\Ui_\ga$ it suffices to show that
\be\ba{rl}\label{TB2}
{\rm (i)}&\dis\limsup_{n\to\infty}\Ui^{(n)}(h_{0,1})
\leq p^\ast_{0,1,\ga^\ast},\\[5pt]
{\rm (ii)}&\dis\liminf_{n\to\infty}\Ui^{(n)}(\ffrac{1}{2}h_{0,1})
\geq p^\ast_{0,1,\ga^\ast},
\ec
uniformly on $[0,1]$. We first consider the case $\ga^\ast>0$.
By (\ref{Uitohom}) and a diagonal argument, we can find $l(n)\to\infty$
such that
\be
\big\|(\Ui_{\ga^\ast})^{l(n)}h_{0,1}-\Ui_{\ga_{n-1}}\circ\cdots\circ
\Ui_{\ga_{n-l(n)}}h_{0,1}\big\|_\infty\asto{n}0.
\ee
Therefore, by Proposition~\ref{01hom}~(a), the fact that $h_{0,1}$
is $\Ui_{\ga_k}$-superharmonic for each $k\geq 0$, and the monotonicity
of the operators $\Ui_\ga$, we find that
\be
\Ui^{(n)}h_{0,1}\leq\Ui_{\ga_{n-1}}\circ\cdots\circ
\Ui_{\ga_{n-l(n)}}h_{0,1}\asto{n}p^\ast_{0,1,\ga^\ast},
\ee
uniformly on $[0,1]$. This proves (\ref{TB2})~(i). To prove also
(\ref{TB2})~(ii) we use the $\Ui_\ga$-{\em sub}harmonic (for each
$\ga>0$) function $h_1$ from Lemma~\ref{subx}. By Lemma~\ref{rhlem}
also $\ffrac{1}{2}h_1$ is $\Ui_\ga$-subharmonic. By bounding
$\ffrac{1}{2}h_1$ from above and below with multiples of $h_{0,1}$ it
is easy to derive from Proposition~\ref{01hom}~(a) that
\be\label{Ughl}
(\Ui_{\ga^\ast})^n(\ffrac{1}{2}h_1)\asto{n}p^\ast_{0,1,\ga^\ast}
\ee
uniformly on $[0,1]$. Arguing as before, we can find $l(n)\to\infty$ such that
\be
\big\|(\Ui_{\ga^\ast})^{l(n)}(\ffrac{1}{2}h_1)-\Ui_{\ga_{n-1}}
\circ\cdots\circ\Ui_{\ga_{n-l(n)}}(\ffrac{1}{2}h_1)\big\|_\infty\asto{n}0.
\ee
Therefore, by (\ref{Ughl}) and the facts that $\ffrac{1}{2}h_1$ is
$\Ui_{\ga_k}$-subharmonic for each $k\geq 0$ and
$\ffrac{1}{2}h_1\leq\ffrac{1}{2}h_{0,1}$,
\be
\Ui^{(n)}(\ffrac{1}{2}h_{0,1})\geq\Ui_{\ga_{n-1}}\circ\cdots\circ
\Ui_{\ga_{n-l(n)}}(\ffrac{1}{2}h_1)\asto{n}p^\ast_{0,1,\ga^\ast},
\ee
uniformly on $[0,1]$, which proves (\ref{TB2})~(ii). The proof
of (\ref{TB2}) in case $\ga^\ast=0$ is completely analogous.\qed

\section{Extinction on the interior}\label{00sec}

\subsection{Basic facts}

In this section we prove Proposition~\ref{00hom}~(a). To simplify
notation, throughout this section $h$ denotes the function $h_{0,0}$.
We fix $0<\ga^\ast<\infty$, we let $Y^h:=Y^{\ga^\ast,h}$ denote the
branching particle system on $(0,1)$ obtained from
$\Yi^{\ga^\ast}=(\Yi^{\ga^\ast}_0,\Yi^{\ga^\ast}_1,\ldots)$ by
Poissonization with $h$ in the sense of Proposition~\ref{Poisprop},
and we denote its log-Laplace operator by $U^h_{\ga^\ast}$. We will
prove that
\be\label{extin}
\rho(x):=P^{\de_x}\big[Y^h_n\neq 0\ \forall n\geq 0\big]=0\qquad(x\in(0,1)).
\ee
Since for each $n$ fixed, $x\mapsto\rho_n(x):=P^{\de_x}[Y^h_n\neq 0]$
is a continuous function that decreases to $\rho(x)$, (\ref{extin})
implies that $\rho_n(x)\to 0$ locally uniformly on $(0,1)$, which, by
an obvious analogon of Lemma~\ref{exgro}, yields
Proposition~\ref{00hom}~(a).

As a first step, we prove:
\bl{\bf (Continuous survival probability)}\label{rholem}
One has either $\rho(x)=0$ for all $x\in(0,1)$ or there exists a
continuous function $\ti\rho:(0,1)\to[0,1]$ such that
$\rho(x)\geq\ti\rho(x)>0$ for all $x\in(0,1)$.
\el
{\bf Proof} Put $p(x):=h(x)\rho(x)$. We will show that either
$p=0$ on $(0,1)$ or there exists a continuous function
$\ti p:(0,1)\to(0,1]$ such that $p\geq\ti p$ on $(0,1)$. Indeed,
\be\ba{l}
\dis p(x)=h(x)P^{\de_x}\big[Y^h_n\neq 0\ \forall n\geq 0\big]
=\lim_{n\to\infty}h(x)P^{\de_x}\big[Y^h_n\neq 0\big]\\[5pt]
\dis\quad=h(x)\lim_{n\to\infty}(U^h_{\ga^\ast})^n1(x)
=\lim_{n\to\infty}(\Ui_{\ga^\ast})^nh(x)\qquad(x\in(0,1)),
\ec
where we have used (\ref{V2}) and (\ref{htrafo}) in the last two steps.
Using the continuity of $\Ui_{\ga^\ast}$ with respect to decreasing
sequences, it follows that
\be\label{fixp}
\Ui_{\ga^\ast}p=p.
\ee

We claim that for any $f\in B_{[0,1]}[0,1]$, one has the bounds
\be\label{Gast}
\li\Ga^\ga_x,f\re\leq\Ui_\ga f(x)\leq(1+\ga)\li\Ga^\ga_x,f\re
\qquad(\ga>0,\ x\in[0,1]).
\ee
Indeed, by Lemma~\ref{Ui01}, $\Ui_\ga f(x)\geq 1-E[(1-f(\y^\ga_x(0)))]
=\li\Ga^\ga_x,f\re$, while the upper bound in (\ref{Gast}) follows
from (\ref{linest}).

By Remark~\ref{R:betadis}, $(0,1)\ni x\mapsto\li\Ga^\ga_x,f\re$ is
continuous for all $f\in B_{[0,1]}[0,1]$. Moreover,
$\li\Ga^\ga_x,f\re=0$ for some $x\in(0,1)$ if and only if $f=0$ almost
everywhere with respect to Lebesgue measure.

Applying these facts to $f=p$ and $\ga=\ga^\ast$, using (\ref{fixp}),
we see that there are two possibilities. Either $p=0$ a.s.\ with
respect to Lebesgue measure, and in this case $p=0$ by the upper bound
in (\ref{Gast}), or $p$ is not almost everywhere zero with respect to
Lebesgue measure, and in this case the function $x\mapsto\ti
p(x):=\li\Ga^\ga_x,f\re$ is continuous, positive on $(0,1)$, and
estimates $p$ from below by the lower bound in (\ref{Gast}).\qed

\subsection{A representation for the Campbell law}\label{002sec}

(Local) extinction properties of critical branching processes are
usually studied using Palm laws. Our proof of formula (\ref{extin}) is
no exception, except that we will use the closely related Campbell
laws. Loosely speaking, Palm laws describe a population that is
size-biased at a given position, plus `typical' particle sampled from
that position, while Campbell laws describe a population that is
size-biased as a whole, plus a `typical' particle sampled from a
random position.

Let $\Pc$ be a probability law on $\Ni(0,1)$ with
$\int_{\Ni(0,1)}\Pc(\di\nu)|\nu|=1$. Then the {\em size-biased law}
$\Pc_{\rm size}$ associated with $\Pc$ is the probability law on
$\Ni(0,1)$ defined by
\be
\Pc_{\rm size}(\,\cdot\,)
:=\int_{\Ni(0,1)}\Pc(\di\nu)\,|\nu|1_{\txt\{\nu\in\cdot\,\}}.
\ee
The {\em Campbell law} associated with $\Pc$ is the probability
law on $(0,1)\times\Ni(0,1)$ defined by
\be\label{Campdef}
\Pc_{\rm Camp}(A\times B):=\int_{\Ni(0,1)}\Pc(\di\nu)\,\nu(A)
1_{\txt\{\nu\in B\}}
\ee
for all Borel-measurable $A\sub(0,1)$ and $B\sub\Ni(0,1)$. If $(v,V)$
is a $(0,1)\times\Ni(0,1)$-valued random variable with law
$\Pc_{\rm Camp}$, then $\Li(V)=\Pc_{\rm size}$, and $v$ is the position
of a `typical' particle chosen from $V$.

Let
\be
\Pc^{x,n}(\,\cdot\,):=P^{\de_x}\big[Y^h_n\in\cdot\,]
\ee
denote the law of $Y^h$ at time $n$, started at time $0$ with one
particle at position $x\in(0,1)$. Note that by criticality,
$\int_{\Ni(0,1)}\Pc^{x,n}(\di\nu)|\nu|=1$. Using again criticality, it
is easy to see that in order to prove the extinction formula
(\ref{extin}), it suffices to show that
\be\label{Sexpl}
\lim_{n\to\infty}\Pc^{x,n}_{\rm size}\big(\{1,\ldots,N\}\big)=0
\qquad(x\in(0,1),\ N\geq 1).
\ee
In order to prove (\ref{Sexpl}), we will write down an expression for
$\Pc^{x,n}_{\rm Camp}$. Let $Q^h$ denote the offspring mechanism of
$Y^h$, and, for fixed $x\in(0,1)$, let $Q^h_{\rm Camp}(x,\,\cdot\,)$
denote the Campbell law associated with $Q^h(x,\,\cdot\,)$. The next
proposition is a time-inhomogeneous version of Kallenberg's famous
backward tree technique; see \cite[Satz~8.2]{Lie81}.
\bp{\bf(Representation of Campbell law)}\label{P:Camprep}
Let $(\vb_k,V_k)_{k\geq 0}$ be the Markov process in $(0,1)\times\Ni(0,1)$
with transition laws
\be
P\big[(\vb_{k+1},V_{k+1})\in\cdot\,\big|\,(\vb_k,V_k)=(x,\nu)\big]
=Q^h_{\rm Camp}(x,\,\cdot\,)\qquad((x,\nu)\in(0,1)\times\Ni(0,1)),
\ee
started in $(\vb_0,V_0)=(\de_x,0)$. Let $(Y^{h,(k)})^{k\geq 1}$ be
branching particle systems with offspring mechanism $Q^h$, conditionally
independent given $(\vb_k,V_k)_{k\geq 0}$, started in
$Y^{h,(k)}_0=V_k-\de_{\vb_k}$. Then
\be\label{Camprep}
\Pc^{x,n}_{\rm Camp}=\Li\Big(\vb_n,\de_{\vb_n}
+\sum_{k=1}^n Y^{h,(k)}_{n-k}\Big).
\ee
\ep
Formula (\ref{Camprep}) says that the Campbell law at time $n$ arises
in such a way, that an `immortal' particle at positions
$\vb_0,\ldots,\vb_n$ sheds off offspring
$V_1-\de_{\vb_1},\ldots,V_n-\de_{\vb_n}$, distributed according to the
size-biased law with one `typical' particle taken out, and this
offspring then evolve under the usual forward dynamics till time $n$.
Note that the position of the immortal particle $(\vb_k)_{k\geq 0}$ is
an autonomous Markov chain.

We need a bit of explicit control on $Q^h_{\rm Camp}$.
\bl{\bf(Campbell law)}\label{L:Campexpl}
One has
\be\label{Campexpl}
Q^h_{\rm Camp}(x,A\times B)=\frac{\frac{1}{\ga^\ast}+1}{h(x)}
\int P[\Pois(h\Zi^{\ga^\ast}_x)\in\di\chi]\chi(A)1_{\{\chi\in A\}},
\ee
where the random measures $\Zi^{\ga^\ast}_x$ are defined in (\ref{Zidef}).
\el
{\bf Proof} By the definition of the Campbell law (\ref{Campdef}),
and (\ref{Zxdef}),
\bc
\dis Q^h_{\rm Camp}(x,A\times B)
&=&\dis\int Q^h(x,\di\chi)\chi(A)1_{\{\chi\in B\}}\\[5pt]
&=&\dis\frac{\frac{1}{\ga^\ast}+1}{h(x)}
\int P[\Pois(h\Zi^{\ga^\ast}_x)\in\di\chi]\chi(A)1_{\{\chi\in B\}}
+\Big(1-\frac{\frac{1}{\ga^\ast}+1}{h(x)}\Big)\cdot 0.
\ec
\qed

\noi
Recall that by (\ref{Zidef}),
\be
\Zi^{\ga^\ast}_x:=\int_0^{\tau_{\ga_\ast}}\de_{\y^{\ga^\ast}_x(-t/2)}\di t,
\ee
where $(\y^{\ga^\ast}_x(t))_{t\in\R}$ is a stationary solution to the
SDE (\ref{Yclx}) with $\ga=\ga^\ast$. By Lemma~\ref{L:Campexpl}, the
transition law of the Markov chain $(\vb_k)_{k\geq 0}$ from
Proposition~\ref{P:Camprep} is given by
\be\label{vbtr}
P[\vb_{k+1}\in\di y|\vb_k=x]=\frac{\frac{1}{\ga^\ast}+1}{h(x)}
E[\Pois(h\Zi^{\ga^\ast}_x)(\di y)]
=\frac{1+\ga^\ast}{h(x)}h(y)\Ga^{\ga^\ast}_x(\di y),
\ee
where $\Ga^{\ga^\ast}_x$ is the invariant law of $\y^{\ga^\ast}_x$
from Corollary~\ref{ergo}. In the next section we will prove the
following lemma.
\bl{\bf(Immortal particle stays in interior)}\label{stayin}
The Markov chain $(\vb_k)_{k\geq 0}$ started in any $\vb_0=x\in(0,1)$ satisfies
\be
(\vb_k)_{k\geq 0}\mbox{ has a cluster point in }(0,1)\quad\as
\ee
\el
We now show that Lemma~\ref{stayin}, together with our previous results,
implies Proposition~\ref{00hom}~(a).\med

\noi
{\bf Proof of Proposition~\ref{00hom}~(a)} We need to prove (\ref{extin}).
By our previous analysis, it suffices to prove (\ref{Sexpl}) under the
assumption that $\rho\neq 0$. By Proposition~\ref{P:Camprep},
\be\label{sizerep}
\Pc^{x,n}_{\rm size}=\Li\Big(\de_{\vb_n}+\sum_{k=1}^n Y^{h,(k)}_{n-k}\Big).
\ee
Conditioned on $(\vb_k,V_k)_{k\geq 0}$, the $(Y^{h,(k)}_{n-k})_{k=1,\ldots,n}$
are independent random variables with
\be
P\big[Y^{h,(k)}_{n-k}\neq 0\big]
\geq P\big[Y^{h,(k)}_m\neq 0\ \forall m\geq 0\big]
=P[\Thin_\rho(V_k-\de_{\vb_k})\neq 0].
\ee
Therefore, (\ref{Sexpl}) will follow by Borel-Cantelli provided
that we can show that
\be\label{BC}
\sum_{k=1}^\infty P[\Thin_\rho(V_k-\de_{\vb_k})\neq 0|\vb_{k-1}]=\infty\quad\as
\ee
Define $f(x):=P[\Thin_\rho(V_k-\de_{\vb_k})\neq 0|\vb_{k-1}=x]$
$(x\in(0,1))$. We need to show that $\sum_{k=1}^\infty f(x)=\infty$
a.s. Using Lemma~\ref{rholem} and Lemma~\ref{L:Campexpl} we can estimate
\be\label{notnil}
f(x)\geq P[\Thin_{\ti\rho}(V_k-\de_{\vb_k})\neq 0|\vb_{k-1}=x]
=\int_{\Ni(0,1)}Q^h_{\rm Camp}(x,\di y,\di\nu)
\{1-(1-\ti\rho)^{\nu-\de_y}\big\}>0
\ee
for all $x\in(0,1)$. Since $\Qi_{\ga^\ast}$, defined in (\ref{Qqdef}),
is a continuous cluster mechanism, also $Q^h_{\rm Camp}(x,\cdot)$ is
continuous as a function of $x$, hence the bound in (\ref{notnil}) is
locally uniform on $(0,1)$, hence Lemma~\ref{stayin} implies that
there is an $\eps>0$ such that
\be
P[\Thin_\rho(V_k-\de_{\vb_k})\neq 0|\vb_{k-1}]\geq\eps
\ee
at infinitely many times $k-1$, which in turn implies (\ref{BC}).\qed

\subsection{The immortal particle}

{\bf Proof of Lemma~\ref{stayin}} Let $K(x,\di y)$ denote the transition
kernel (on $(0,1)$) of the Markov chain $(\vb_k)_{k\geq 0}$, i.e., by
(\ref{vbtr}),
\be
K(x,\di y)=(1+\ga^\ast)\frac{y(1-y)}{x(1-x)}\Ga^{\ga^\ast}_x(\di y).
\ee
It follows from (\ref{WFmoments}) that
\be\label{Kyy}
\int K(x,\di y)y(1-y)
=\frac{x(1-x)+\ga^\ast(1+\ga^\ast)}{(1+2\ga^\ast)(1+3\ga^\ast)}.
\ee
\detail{Indeed,
\[\ba{l}
\dis\int K(x,\di y)y(1-y)=\frac{1+\ga}{x(1-x)}
\int\Ga^\ga_x(\di y)\{y^2-2y^3+y^4\}\\[5pt]
\dis\quad=\frac{1+\ga}{x(1-x)}\Big\{\left(\frac{x+0}{1+0}\right)
\left(\frac{x+\ga}{1+\ga}\right)-2\left(\frac{x+0}{1+0}\right)
\left(\frac{x+\ga}{1+\ga}\right)\left(\frac{x+2\ga}{1+2\ga}\right)\\
\dis\qquad\phantom{=}+\left(\frac{x+0}{1+0}\right)
\left(\frac{x+\ga}{1+\ga}\right)\left(\frac{x+2\ga}{1+2\ga}\right)
\left(\frac{x+3\ga}{1+3\ga}\right)\Big\}\\[5pt]
\dis\quad=\frac{x(x+\ga)}{x(1-x)(1+2\ga)(1+3\ga)}
\Big\{(1+2\ga)(1+3\ga)-2(x+2\ga)(1+3\ga)+(x+2\ga)(x+3\ga)\Big\}\\[5pt]
\dis\quad=\frac{x+\ga}{(1-x)(1+2\ga)(1+3\ga)}
\Big\{(1+2\ga)(1+3\ga)-2(((x-1)+1+2\ga)(1+3\ga)\\
\dis\quad\phantom{=\frac{x+\ga}{(1-x)(1+2\ga)(1+3\ga)}
\Big\{}+((x-1)+1+2\ga)((x-1)+1+3\ga)\Big\}\\[5pt]
\dis\quad=\frac{x+\ga}{(1-x)(1+2\ga)(1+3\ga)}
\Big\{(x-1)\big\{-2(1+3\ga)+(1+2\ga)+(1+3\ga)\big\}+(x-1)^2\Big\}\\[5pt]
\dis\quad=\frac{(x+\ga)\{\ga+(1-x)\}}{(1+2\ga)(1+3\ga)}
=\frac{x(1-x)+\ga(1+\ga)}{(1+2\ga)(1+3\ga)}.
\ea\]}
Set
\be
g(x):=\int K(x,\di y)y(1-y)-x(1-x)\qquad(x\in(0,1)).
\ee
Then
\be
M_n:=\vb_n(1-\vb_n)-\sum_{k=0}^{n-1}g(\vb_k)\qquad(n\geq 0)
\ee
defines a martingale $(M_n)_{n\geq 0}$. Since $g>0$ in an open
neighborhood of $\{0,1\}$,
\be
P[(\vb_k)_{k\geq 0}\mbox{ has no cluster point in }(0,1)]\leq
P[\lim_{n\to\infty}M_n=-\infty]=0,
\ee
where in the last equality we have used that $(M_n)_{n\geq 0}$
is a martingale.\qed

\section{Proof of the main result}\label{final}

{\bf Proof of Theorem~\ref{main}} Part~(a) has been proved in
Section~\ref{dualsec}. It follows from (\ref{sumga}), (\ref{gamma}),
(\ref{hutform}), and (\ref{renbra}) that part~(b) is equivalent to
the following statement. Assuming that
\be\label{varsga}
{\rm (i)}\quad\sum_{n=1}^\infty\ga_n=\infty\qquad\mbox{and}
\qquad{\rm (ii)}\quad\ga_n\asto{n}\ga^\ast
\ee
for some $\ga^\ast\in\half$, one has, uniformly on $[0,1]$,
\be\label{varwast}
\Ui_{\ga_{n-1}}\circ\cdots\circ\Ui_{\ga_0}(p)\asto{n}p^\ast_{l,r,\ga^\ast},
\ee
where $p^\ast_{l,r,\ga^\ast}$ is the unique solution in $\Hi_{l,r}$ of
\be\ba{rr@{\,}c@{\,}ll}\label{varwiga}
{\rm (i)}&\Ui_{\ga^\ast}p^\ast&=&p^\ast\quad
&\mbox{if }0<\ga^\ast<\infty,\\[5pt]
{\rm (ii)}&\ffrac{1}{2}x(1-x)\diff{x}p^\ast(x)-p^\ast(x)(1-p^\ast(x))
&=&0\quad(x\in[0,1])\quad&\mbox{if }\ga^\ast=0.
\ec
It follows from Proposition~\ref{Uit} that the left-hand side of
(\ref{varwast}) converges uniformly to a limit $p^\ast_{l,r,\ga^\ast}$
which is given by (\ref{pdef}). We must show $1^\circ$ that
$p^\ast_{l,r,\ga^\ast}\in\Hi_{l,r}$ and $2^\circ$ that
$p^\ast_{l,r,\ga^\ast}$ is the unique solution in this class to
(\ref{varwiga}). We first treat the case $\ga^\ast>0$.

$1^\circ$ Since $p^\ast_{0,0,\ga^\ast}\equiv 0$ and
$p^\ast_{1,1,\ga^\ast}\equiv 1$, it is obvious that
$p^\ast_{0,0,\ga^\ast}\in\Hi_{0,0}$ and
$p^\ast_{1,1,\ga^\ast}\in\Hi_{1,1}$. Therefore, by symmetry, it
suffices to show that $p^\ast_{0,1,\ga^\ast}\in\Hi_{0,1}$. By
Lemmas~\ref{subx} and \ref{sup01}, $x\leq p\leq 1-(1-x)^7$ implies
$x\leq\Ui_{\ga_k}p\leq 1-(1-x)^7$ for each $k$. Iterating this
relation, using (\ref{varwast}), we find that
\be\label{lowup}
x\leq p^\ast_{0,1,\ga^\ast}(x)\leq 1-(1-x)^7.
\ee
By Proposition~\ref{moncon}, the left-hand side of (\ref{varwast}) is
nondecreasing and concave in $x$ if $p$ is, so taking the limit we
find that $p^\ast_{0,1,\ga^\ast}$ is nondecreasing and concave.
Combining this with (\ref{lowup}) we conclude that
$p^\ast_{0,1,\ga^\ast}$ is Lipschitz continuous. Moreover
$p^\ast_{0,1,\ga^\ast}(0)=0$ and $p^\ast_{0,1,\ga^\ast}(1)=1$ so
$p^\ast_{0,1,\ga^\ast}\in\Hi_{0,1}$.

$2^\circ$ Taking the limit $n\to\infty$ in
$(\Ui_{\ga^\ast})^{n}p=\Ui_{\ga^\ast}(\Ui_{\ga^\ast})^{n-1}p$, using
the continuity of $\Ui_{\ga^\ast}$ (Corollary~\ref{Ugacont}) and
(\ref{varwast}), we find that
$\Ui_{\ga^\ast}p^\ast_{l,r,\ga^\ast}=p^\ast_{l,r,\ga^\ast}$. It
follows from (\ref{varwast}) that $p^\ast_{l,r,\ga^\ast}$ is the only
solution in $\Hi_{l,r}$ to this equation.

For $\ga^\ast=0$, it has been shown in \cite[Proposition~3]{FSsup}
that $p^\ast_{l,r,0}$ is the unique solution in $\Hi_{l,r}$ to
(\ref{varwiga})~(ii). In particular, it has been shown there that
$p^\ast_{0,1,0}$ is twice continuously differentiable on $[0,1]$
(including the boundary). This proves parts~(b) and (c) of the
theorem.\qed

\appendix

\section{Appendix: Infinite systems of linearly interacting
diffusions}\label{app}

\subsection{Hierarchically interacting diffusions}\label{hier}

For any $N\geq 2$, the {\em hierarchical group} with freedom $N$ is
the set $\om_N$ of all sequences $\xi=(\xi_1,\xi_2,\ldots)$, with
coordinates $\xi_k$ in the finite set $\{0,\ldots,N-1\}$, which are
different from $0$ only finitely often, equipped with componentwise
addition modulo $N$. Setting
\be
\|\xi\|:=\min\{n\geq 0:\xi_k=0\ \forall k>n\}\qquad(\xi\in\om_N),
\ee
$\|\xi-\eta\|$ is said to be the {\em hierarchical distance}
between two sites $\xi$ and $\eta$ in $\om_N$.

Let $D\sub\R^d$ be open and convex, and let $\Wi$ be a renormalization
class on $\ov D$. Let $\sig$ be a continuous root of a diffusion
matrix $w\in\Wi$ as in Remark~\ref{sderem}. Consider a collection
$\x=(\x_\xi)_{\xi\in\om_N}$ of $\ov D$-valued processes, solving a
system of SDE's of the form
\be\label{infsde}
\di\x_\xi(t)=\sum_{k=0}^\infty\frac{c_k}{N^k}
\Big(\x^{k+1}_\xi(t)-\x_\xi(t)\Big)\di t+\sqrt{2}\sig(\x_\xi(t))\di B_\xi(t)
\qquad(t\geq 0,\ \xi\in\om_N),
\ee
where $(B_\xi)_{\xi\in\om_N}$ is a collection of independent standard
Brownian motions, with initial condition
\be\label{infinit}
\x_\xi(0)=\tet\in D\qquad(\xi\in\om_N).
\ee
Here the $(c_k)_{k\geq 0}$ are positive constants satisfying
$\sum_k c_k/N^k<\infty$, and $\x^k_\xi(t)$ denotes the
{\em $k$-block average} around $\xi$:
\be\label{kblock}
\x^k_\xi(t):=\frac{1}{N^k}\sum_{\eta:\|\xi-\eta\|\leq k}\x_\eta(t)
\qquad\qquad(k\geq 0).
\ee
(Note that $|\{\eta:\|\xi-\eta\|\leq k\}|=N^k$.) Under suitable
additional assumptions on $\sig$, one can show that (\ref{infsde}) has
a unique (weak or strong) solution (see \cite{DG93a,DG96,Swa00}). We
call $\x$ a system of {\em hierarchically interacting $\ov D$-valued
diffusions} with migration constants $(c_k)_{k\geq 0}$ and local
diffusion rate $w_{ij}=\sum_k\sig_{ik}\sig_{jk}$. Such systems are
used to model gene frequencies or population sizes in population
biology \cite{SF83}.

The long-time behavior of the system in (\ref{infsde}) depends
crucially on the recurrence versus transience of the continuous-time
random walk on $\om_N$ which jumps from a point $\xi$ to a point
$\eta\neq\xi$ with rate
\be\label{a}
a(\eta-\xi):=\sum_{k=\|\xi-\eta\|}^\infty\frac{c_{k-1}}{N^{2k-1}}.
\ee
This random walk is recurrent if and only if
\be\label{rec}
\sum_{k=0}^\infty\frac{1}{d_k}=\infty,\qquad\mbox{where}
\quad d_k:=\sum_{n=0}^\infty\frac{c_{k+n}}{N^n}
\ee
(see \cite{DG93a,Kle96}; a similar problem is treated in \cite{DE68}).
Assuming that the law of $\x(t)$ converges weakly as $t\to\infty$ to
the law of some $\ov D^{\om_N}$-valued random variable $\x(\infty)$,
one expects that in the recurrent case $\x(\infty)$ must have the
following properties:
\be\ba{rll}\label{clust}
{\rm (i)}&\x_\xi(\infty)=\x_\eta(\infty)\quad&\mbox{a.s.}
\qquad\forall\xi,\eta\in\om_N,\\
{\rm (ii)}&\x_\xi(\infty)\in\pa_w D\quad&\mbox{a.s.}\qquad\forall\xi\in\om_N.
\ec
Here $\pa_w D$ is the effective boundary of $D$, defined in
(\ref{eff}). If $\x(t)$ converges in law to a limit $\x(\infty)$
satisfying (\ref{clust}), then we say that $\x$ {\em clusters}. In the
transient case, it is believed that solutions of (\ref{infsde}) do not
cluster. (For compact $\ov D$ these facts were proved in
\cite{Swa00}.)

An important tool in the study of solutions to (\ref{infsde}) is the
so-called {\em interaction chain}. This is the chain
$(\x^0_0(t),\x^1_0(t),\ldots)$ of block-averages around the origin.
Heuristic arguments suggest that in the {\em local mean field limit}
$N\to\infty$, the interaction chain converges to a certain
well-defined Markov chain.
\begin{conjecture}\label{Ichain}
Fix $w\in\Wi$, $\tet\in D$, and positive numbers $(c_k)_{k\geq 0}$
such that for $N$ large enough, $\sum_k c_k/N^k<\infty$. For all $N$
large enough, let $\x^N$ be a solution to
(\ref{infsde})--(\ref{infinit}), and assume that $t_N$ are constants
such that, for some $n\geq 1$,
$\lim_{N\to\infty}N^{-n}t_N=T\in\half$. Then
\be\label{toI}
\Big(\x^{N,n}_0(t_N),\ldots,\x^{N,0}_0(t_N)\Big)
\underset{N\to\infty}{\dis\Longrightarrow}(I^{w}_{-n},\ldots,I^{w}_0),
\ee
where $(I^{w}_{-n},\ldots,I^{w}_0)$ is a Markov chain with transition laws
\be\label{transspec}
P[I^{w}_{-k}\in\di y|I^{w}_{-k-1}=x]=\nu^{c_k,F^{(k)}w}_x(\di y)
\qquad(x\in\ov D,\ 0\leq k\leq n-1)
\ee
and initial state
\be\label{tinit}
I^{w}_{-n}=\y_T,\quad\mbox{where}\quad\di\y_t=c_n(\tet-\y_t)\di t
+\sqrt{2}\sig^{(n)}(\y_t)\di B_t,\quad\y_0=\tet,
\ee
and $\sig^{(n)}$ is a root of the diffusion matrix $F^{(n)}w$.
\end{conjecture}
Rigorous versions of conjecture~\ref{Ichain} have been proved for
renormalization classes on $\ov D=[0,1]$ and $\ov D=\half$ in
\cite{DG93a,DG93b}.

Note that the iterated kernels $K^{w,(n)}$ defined in (\ref{Kdef})
are the transition probabilities from time $-n$ to time $0$ of the
interaction chain in the mean-field limit:
\be\label{Kdef2}
K^{w,(n)}_x(\di y)=P[I^{w}_0\in\di y|I^{w}_{-n}=x]
\qquad\quad(x\in\ov D,\ n\geq 0).
\ee
Lemma~\ref{clusK} expresses the fact that the system $\x^N$ clusters
in the local mean-field limit $N\to\infty$. The condition
$s_n\to\infty$ in Lemma~\ref{clusK} means that $\sum_{k\geq
  0}\frac{1}{c_k}=\infty$, which, in a sense, is the $N\to\infty$
limit of condition (\ref{rec}).

\subsection{The clustering distribution of linearly interacting
diffusions}\label{clustapp}

Let $D\sub\R^d$ be open, {\em bounded}, and convex, and let $\Wi$ be a
renormalization class on $\ov D$. Fix migration constants
$(c_k)_{k\geq 0}$ and assume that $s_n\to\infty$ and $s_{n+1}/s_n\to
1+\ga^\ast$ for some $\ga^\ast\in[0,\infty]$. Recall the definition of
the iterated probability kernels $K^{w,(n)}$ in (\ref{Kdef}). Recall
Conjecture~\ref{con}. Assuming that the rescaled renormalized
diffusion matrices $s_nF^{(n)}w$ converge to a limit $w^\ast$, we can
make a guess about the limit of the iterated probability kernels
$K^{w,(n)}$.
\begin{conjecture}{\bf(Limits of iterated probability kernels)}\label{Kcon}
Assume that $s_nF^{(n)}w\to w^\ast$ as $n\to\infty$. Then, for any $w\in\Wi$,
\be
K^{w,(n)}\asto{n}K^\ast,
\ee
where $K^\ast$ has the following description:
\renewcommand{\labelenumi}{\rm(\roman{enumi})}
\begin{enumerate}
\item If $0<\ga^\ast<\infty$, then
\be
K^\ast_x=\lim_{n\to\infty}P^x[I^{\ga^\ast}_n\in\,\cdot\,],
\ee
where $(I^{\ga^\ast}_n)_{n\geq 0}$ is the Markov chain with
transition law $P[I^{\ga^\ast}_{n+1}\in\cdot\,|I^{\ga^\ast}_n=x]
=\nu^{1/\ga^\ast,w^\ast}$. 
\item If $\ga^\ast=0$, then
\be\label{Kcrit}
K^\ast_x=\lim_{t\to\infty}P^x[I^0_t\in\,\cdot\,],
\ee
where $(I^0_s)_{s\geq 0}$ is the diffusion process with generator
$\sum_{i,j=1}^dw^\ast_{ij}(y)\difif{y_i}{y_j}$.
\item If $\ga^\ast=\infty$, then
\be
K^\ast_x=\lim_{\ga\to\infty}\nu^{1/\ga,w^\ast}_x.
\ee
\end{enumerate}
\end{conjecture}
For each $N\geq 2$, let $\x^N=(\x^N_\xi)_{\xi\in\om_N}$ be a system of
hierarchically interacting diffusions as in (\ref{infsde}) and
(\ref{infinit}). If $\ga^\ast=0$, then because of
Conjectures~\ref{Ichain} and \ref{Kcon}, we expect\footnote{For
$\ga^\ast>0$, the situation is more complex. In this case at the
right-hand side of (\ref{dublim}) we expect the law $\int_{\ov D}
P^\tet[\y_T\in\di x]K^\ast_x$, where $\y$ solves the SDE
$\di\y_t=\frac{1}{\ga^\ast}(\tet-\y_t)\di t+\sqrt2\sig^\ast(\y_t)\di
B_t$ and $\sig^\ast$ is a root of the diffusion matrix $w^\ast$.
Note that in this case the right-hand side of (\ref{dublim}) depends
on $T$.} that
\be\label{dublim}
\lim_{n\to\infty}\lim_{N\to\infty}\Li(\x^N_0(N^nT))=K^\ast_\tet\qquad(T>0),
\ee
where $K^\ast$ is the kernel in (\ref{Kcrit}).

In particular, consider the case that the migration constants
$(c_k)_{k\geq 0}$ are of the form $c_k=r^k$ for some $r>0$. In this
case, $s_{n+1}/s_n\to\frac{1}{r}\vee 1$, and $s_n\to\infty$ if and
only if $r\leq 1$. One can check (see (\ref{rec})) that for fixed
$N\geq 2$, the random walk with the kernel $a$ in (\ref{a}) is
recurrent if and only if $r\leq 1$. The critical case $r=1$
corresponds to a {\em critically recurrent} random walk. For a precise
definition of critical recurrence, see \cite[formula~(1.15)]{Kle96}.
For $r=1$, we expect that the double limit in (\ref{dublim}) can be
replaced by a single limit. More precisely, for each fixed $N\geq 2$,
we expect that
\be\label{enklim}
\lim_{t\to\infty}\Li(\x^N_0(t))=K^\ast_\tet.
\ee
In this case, we call $K^\ast_\tet$ the {\em clustering distribution}
of $\x^N$. The clustering distribution of linearly interacting
isotropic diffusions was studied in \cite{Swa00}. We expect
(\ref{enklim}) to hold, even more generally, for all systems of
linearly interacting diffusions with a critically recurrent migration
mechanism. In particular, we expect (\ref{enklim}) to hold for
symmetric nearest-neighbor interaction on $\Z^d$ in the critical
dimension $d=2$. If one is ready to make this enormous leap of faith,
then combining Conjectures~\ref{con} and \ref{Kcon}, one arrives at
the following conjecture.
\begin{conjecture}{\bf (Critical clustering)}\label{Z2con}
Let $D\sub\R^d$ be open, bounded, and convex, and let $\Wi$ be a
renormalization class on $\ov D$. Assume that the asymptotic fixed
point equation (\ref{afp})~(ii) has a unique solution $w^\ast$ in
$\Wi$. Let $\sig$ be a continuous root of a diffusion matrix
$w\in\Wi$. Let $\x=(\x_\xi)_{\xi\in\Z^2}$ be a $\ov D^{\Z^2}$-valued
process, solving the system of SDE's
\be
\di\x_\xi(t)=\sum_{\eta:\,|\eta-\xi|=1}
\!\big(\x_\eta(t)-\x_\xi(t)\big)\,\di t+\sig(\x_\xi(t))\di B_\xi(t),
\ee
with initial condition $\x_\xi(0)=\tet\in\ov D$ $(\xi\in\Z^2)$. Then
\be
\x_\xi(t)\Asto{t}I^\tet_\infty\qquad(\xi\in\Z^2),
\ee
where $(I^\tet_s)_{s\geq 0}$ is the diffusion with generator
$\sum_{i,j}w^\ast_{ij}(y)\difif{y_i}{y_j}$ and initial condition
$I^\tet_0=\tet$.
\end{conjecture}

\subsection*{Acknowledgements}

We thank Anton Wakolbinger and Martin M\"ohle for pointing out
reference \cite{Ewe04} and the fact that the distribution in
(\ref{Gadef}) is a $\bet$-distribution.

\newcommand{\noopsort}[1]{}

\end{document}